\documentclass[11pt,a4paper]{amsart}

\usepackage[T1]{fontenc}
\usepackage[utf8]{inputenc}
\usepackage{lmodern}
\usepackage{microtype}
\usepackage[a4paper,margin=2.65cm]{geometry}
\usepackage{amsmath,amssymb,amsfonts,amsthm,mathtools}
\usepackage{mathrsfs}
\usepackage{enumitem}
\usepackage{xcolor}
\usepackage{cite}
\usepackage{hyperref}

\hypersetup{
  colorlinks=true,
  linkcolor=blue!85!black,
  citecolor=blue!85!black,
  urlcolor=blue!85!black,
  pdftitle={Nonlinear Darcy Laws for Generalized Newtonian Fluids with Orlicz Growth and Oscillating Navier Slip in Periodic Porous Media},
  pdfauthor={Hermann Douanla and Joel Fotso Tachago and Jean Louis Woukeng},
  pdfsubject={Periodic homogenization of generalized Newtonian flows with Orlicz growth and nonlinear Navier slip},
  pdfkeywords={periodic homogenization, generalized Newtonian fluids, Orlicz-Sobolev spaces, surface two-scale convergence, nonlinear Navier slip, nonlinear Darcy law}
}

\allowdisplaybreaks
\setlist[enumerate]{leftmargin=2.2em,itemsep=0.25em,topsep=0.35em}
\setlist[itemize]{leftmargin=2.0em,itemsep=0.25em,topsep=0.35em}
\numberwithin{equation}{section}
\newtheorem{theorem}{Theorem}[section]
\newtheorem{proposition}[theorem]{Proposition}
\newtheorem{lemma}[theorem]{Lemma}
\newtheorem{corollary}[theorem]{Corollary}
\theoremstyle{definition}
\newtheorem{definition}[theorem]{Definition}
\newtheorem{remark}[theorem]{Remark}
\newtheorem{assumption}[theorem]{Assumption}
\newtheorem*{theoremA}{Theorem A}

\newcommand{\R}{\mathbb{R}}
\newcommand{\Z}{\mathbb{Z}}
\newcommand{\D}{\mathbb{D}}
\newcommand{\eps}{\varepsilon}
\newcommand{\norm}[1]{\left\lVert #1\right\rVert}
\newcommand{\abs}[1]{\left|#1\right|}
\newcommand{\dd}{\,\mathrm{d}}

\newcommand{\diver}{\operatorname{div}}
\newcommand{\supp}{\operatorname{supp}}

\title[Nonlinear Darcy laws for generalized Newtonian fluids]
{Nonlinear Darcy Laws for Generalized Newtonian Fluids with Orlicz Growth and Oscillating Navier Slip in Periodic Porous Media}

\author[H. Douanla]{Hermann Douanla$^{\dagger}$}
\thanks{$^{\dagger}$Department of Mathematics, Faculty of Science, University of Yaound\'e I, P.O. Box 812, Yaound\'e, Cameroon; \texttt{hermann.douanla@facsciences-uy1.cm}.}

\author[J. Fotso Tachago]{Joel Fotso Tachago$^{\ddagger}$}
\thanks{$^{\ddagger}$Department of Mathematics, Higher Teachers Training College, University of Bamenda, P.O. Box 39, Bambili, Cameroon; \texttt{fotsotachago@yahoo.fr}.}

\author[J. L. Woukeng]{Jean Louis Woukeng$^{*}$}
\thanks{$^{*}$Corresponding author. Department of Mathematics and Computer Science, University of Dschang, P.O. Box 67, Dschang, Cameroon; \texttt{jeanlouis.woukeng@univ-dschang.org}.}

\begin{document}

\begin{abstract}
We study the periodic homogenization of a stationary generalized Navier--Stokes system in a periodically perforated domain $\Omega\subset\mathbb R^N$, $N\in\{2,3\}$, with a heterogeneous bulk stress of general Orlicz growth and an oscillating nonlinear Navier friction law on the pore walls, both at the Darcy scale. If the lower growth exponent at infinity satisfies $p^-_\infty>3N/(N+2)$ and the inertial exponent exceeds $\kappa_c:=1+N(2/p^-_\infty-1)_+$, inertia vanishes and, under strict monotonicity of the bulk stress, the whole family converges to the nonlinear Darcy law $U=\mathcal K(f-\nabla p_0)$, $\operatorname{div}U=0$, $U\cdot\nu=0$. The permeability map $\mathcal K$ is continuous, strictly monotone, obeys two-sided Orlicz bounds and incorporates the wall slip through a nonlinear cell problem, and the extended microscopic pressure converges strongly to $p_0$. The proof relies on a uniform Orlicz--Korn inequality, a Darcy-scaled two-scale trace theorem and a coupled bulk--surface Minty argument.


\end{abstract}

\maketitle

\section{Introduction}
\label{sec:introduction}

The rigorous derivation of Darcy's law by homogenization of the Stokes system in periodically perforated
domains is classical under linear viscosity and no-slip conditions on the pore walls
\cite{Allaire1989,Hornung1997,SanchezPalencia1980,Tartar1980}.  The present work replaces both assumptions
simultaneously by a generalized Newtonian constitutive law with Orlicz growth and an oscillating nonlinear
Navier friction law.  As a consequence, the effective permeability is in general a nonlinear map.

\subsection{Model and motivation}
\label{subsec:motivation-objectives}

We consider porous transport of non-Newtonian fluids whose stress--strain relation is not accurately
represented by a single power law and for which wall slip is significant.  Shear-thinning and
Carreau-type rheologies are standard in polymer and suspension flows through porous media
\cite{OrgeasGeindreauAuriaultBloch2007}, and, for Newtonian flow, wall slip raises
the permeability from its no-slip value towards the perfect-slip value, which bound it from below and
above \cite{LasseuxValdesParada2024}.  General Orlicz growth describes several rheological regimes within one
constitutive law, while the nonlinear Navier law models tangential dissipation at the solid--fluid
interface.

Let $N\in\{2,3\}$ and $Y=(0,1)^N$.  The geometry satisfies the following assumption.

\begin{assumption}\label{ass:geometry}
The domain $\Omega\subset\mathbb R^N$ is bounded, connected and of class $C^{1,1}$, and
$T\Subset Y$ is a nonempty connected open set with $\partial T\in C^{1,1}$ such that
$Y^*:=Y\setminus\overline T$ is connected.
\end{assumption}

For $0<\eps<1$, define the set of complete cells, the union of holes, the perforated domain and the union
of pore walls by
\[
K_\eps=\{k\in\Z^N:\ \eps(k+\overline Y)\subset\Omega\},
\qquad
T_\eps=\bigcup_{k\in K_\eps}\eps(k+T),
\]
\[
\Omega_\eps=\Omega\setminus\overline{T_\eps},
\qquad
\Gamma_\eps=\bigcup_{k\in K_\eps}\eps(k+\partial T).
\]
Only complete cells are perforated, so that no hole meets $\partial\Omega$.  Let $n(y)$ be the unit
normal to $\partial T$ pointing out of $Y^*$ and into $T$; the outward unit normal $n_\eps$ to $\Omega_\eps$
on $\Gamma_\eps$ is its rescaling, $n_\eps(x)=n(x/\eps)$.  For $z\in\R^N$ we set
$P_\tau(y):=I-n(y)\otimes n(y)$ and $z_\tau=P_\tau(y)z$, and for a vector field $v$ we use the symmetric
gradient $\D v:=\frac12\big(\nabla v+(\nabla v)^T\big)$.  For $A,B\in\mathbb R^{N\times N}$ we write
$A:B:=\sum_{i=1}^{N}\sum_{j=1}^{N}A_{ij}B_{ij}$ and $|A|:=(A:A)^{1/2}$, and for $a\in\R^m$ and
$b\in\R^N$, $a\otimes b\in\R^{m\times N}$ is the matrix with entries $(a\otimes b)_{ij}:=a_ib_j$, that is,
$a\otimes b=ab^T$.  For a vector field $v$, $(\nabla v)_{ij}:=D_jv_i$, and the divergence of a
matrix field is taken row by row; a subscript $x$ or $y$, as in $\nabla_y$, $\D_y$, $\diver_y$ or $\diver_x$,
restricts the derivatives to that variable.  We denote by $\sigma_x$ and $\sigma_y$ the surface measures on
$\Gamma_\eps$ and $\partial T$, and by $\nu$ the outer unit normal to $\partial\Omega$.  For a Sobolev function $v$ on a Lipschitz domain,
$\gamma_0v$ denotes its trace on the boundary; for a function $v(x,y)$, $\gamma_0^yv$ denotes the trace of
$v(x,\cdot)$ on $\partial T$ with respect to $y$, the variable $x$ being fixed.

For a fixed real exponent $\kappa$, the microscopic problem is
\begin{equation}
\label{Peps}
\left\{
\begin{aligned}
-\eps\,\diver\!\left[
\mathbb S\!\left(\frac{x}{\eps},\eps\D u_\eps\right)
\right]
+\eps^\kappa\diver(u_\eps\otimes u_\eps)+\nabla p_\eps
&=f &&\text{in }\Omega_\eps,\\
\diver u_\eps&=0 &&\text{in }\Omega_\eps,\\
u_\eps&=0 &&\text{on }\partial\Omega,\\
u_\eps\cdot n_\eps&=0 &&\text{on }\Gamma_\eps,\\
\left[
\eps\mathbb S\!\left(\frac{x}{\eps},\eps\D u_\eps\right)n_\eps
\right]_\tau
+\eps\mathcal B\!\left(\frac{x}{\eps},(u_\eps)_\tau\right)
&=0 &&\text{on }\Gamma_\eps,
\end{aligned}
\right.
\end{equation}
where $\mathbb S$ is the bulk stress and $\mathcal B$ the tangential wall friction, both $Y$-periodic in
their first variable.  The scaling is that of Darcy flow.  In pores of size $\eps$ carrying a velocity of
order one, the strain rate is of order $\eps^{-1}$; the constitutive law is therefore evaluated at the
scaled strain rate $\eps\D u_\eps$, and the factor $\eps$ in front of the divergence makes the viscous force
of order one, balanced by a pressure gradient of order one.  For a power law
$\mathbb S(y,E)=\mu|E|^{p-2}E$ this amounts, after rescaling of the velocity, to the classical choice of a viscosity of order
$\eps^{p}$; compare \cite[Theorem~9]{BourgeatMikelic1996}.  For a non-homogeneous law, this scaling fixes
the rheological scales, such as the time constant of a Carreau law, at the pore scale.  The factor $\eps^\kappa$ measures the strength of inertia: since the viscous term is of order one, the ratio
of inertial to viscous forces at the pore scale is $\eps^{\kappa-1}$, which is the local Reynolds number of
\cite{MarusicPalokaMikelic2000}, and larger values of $\kappa$ correspond to weaker inertia
\cite{AnguianoBonnivardSuarezGrau2026,MarusicPalokaMikelic2000}.

For smooth solutions, testing \eqref{Peps} by $u_\eps$ and using $\diver u_\eps=0$, $u_\eps=0$ on
$\partial\Omega$ and $u_\eps\cdot n_\eps=0$ on $\Gamma_\eps$ gives the energy identity
\begin{equation}
\label{introenergyidentity}
\int_{\Omega_\eps}\mathbb S\!\left(\frac{x}{\eps},\eps\D u_\eps\right):\eps\D u_\eps\,\dd x
+\eps\int_{\Gamma_\eps}\mathcal B\!\left(\frac{x}{\eps},(u_\eps)_\tau\right)\cdot(u_\eps)_\tau\,\dd\sigma_x
=\int_{\Omega_\eps}f\cdot u_\eps\,\dd x,
\end{equation}
in which the convection does not contribute.  Its weak counterpart is given in Definition~\ref{def:weakvelocity}.
Since there are of order $\eps^{-N}$ holes, each with surface area of order $\eps^{N-1}$, the total pore-wall
area is of order $\eps^{-1}$, and the factor $\eps$ in the slip law places the bulk and wall dissipations in
\eqref{introenergyidentity} at the same order \cite{AllaireDamlamianHornung1996,NeussRadu1996}; physically, it
corresponds to a slip length comparable to the pore size.

The central question is whether the two microscopic dissipation mechanisms in
\eqref{introenergyidentity} combine into a closed macroscopic filtration law, and how the wall friction
enters the effective permeability.

\subsection{Literature review}
\label{subsec:literature}

Nonlinear Darcy laws for power-law and Carreau fluids were derived by Bourgeat--Mikeli\'c
\cite{BourgeatMikelic1996} and decoupled into well-posed macroscopic filtration laws by
Bourgeat--Gipouloux--Maru\v{s}i\'c-Paloka \cite{BourgeatGipoulouxMarusicPaloka2003}; effective permeation laws
for generalized-Newtonian flow through anisotropic porous media are derived in
\cite{OrgeasGeindreauAuriaultBloch2007}, and lower-dimensional linear and nonlinear Carreau laws in thin
porous media in \cite{AnguianoBonnivardSuarezGrau2022,AnguianoBonnivardSuarezGrau2025}.

The works closest to the present bulk rheology are the following.  Kalousek \cite{Kalousek2016} homogenized
steady and unsteady generalized Stokes systems with Orlicz-growth stresses under the same Darcy scaling,
using two-scale convergence in the Orlicz setting, and obtained a nonlinear two-scale limit system without
inertia.  Lu--Qian \cite{LuQian2024} and H\"ofer--Lu--Oschmann \cite{HoferLuOschmann2026} treated
Carreau--Yasuda flows, the latter through holes much smaller than the period, in scalings for which the
nonlinear viscosity disappears in the limit, and obtained linear Darcy laws, with quantitative rates in
\cite{HoferLuOschmann2026}.  All three works impose the no-slip condition on the pore walls.  Here, by
contrast, the walls carry a nonlinear Navier friction, the convective term is retained in the microscopic
problem and shown to vanish above an explicit threshold, and the nonlinear permeability map is analysed
together with the strong convergence of the microscopic pressure.

For Newtonian flow with no-slip conditions, Maru\v{s}i\'c-Paloka--Mikeli\'c \cite{MarusicPalokaMikelic2000}
showed that inertia does not affect the Darcy limit when the local Reynolds number tends to zero, and derived,
for a local Reynolds number of order one and small forces, a nonlinear filtration law with a Navier--Stokes
cell problem, which for small velocities yields a non-classical Dupuit--Ergun--Forchheimer law.  For Carreau fluids in thin porous media,
Anguiano--Bonnivard--Su\'arez-Grau \cite{AnguianoBonnivardSuarezGrau2026} showed that inertia does not affect
the limit up to a critical Reynolds number of order $1/\eps$, and solved numerically the nonlinear Darcy law
of the critical case.

Slip conditions in porous-flow homogenization were studied by Allaire \cite{AllaireSlip1991} for holes smaller
than the period, and by Cioranescu--Donato--Ene \cite{CioranescuDonatoEne1996} and Capatina--Ene
\cite{CapatinaEne2011} for inclusions of the size of the period.  The limit law depends on the order of the
slip coefficient (Darcy, Brinkman or Stokes type in \cite{CioranescuDonatoEne1996}), and the friction enters
the local problem when this coefficient is of order $\eps^{-1}$ relative to the viscous stress
\cite{AnguianoSuarezGrau2026}; the present scaling is precisely this critical order, which lies between the
regimes treated in \cite{CioranescuDonatoEne1996}.  Effective Navier laws over oscillating walls were derived in
\cite{MikelicNecasovaNeussRadu2013}.  In thin perforated layers, Navier slip
\cite{FabriciusGahn2023} and non-homogeneous slip on the obstacles
\cite{AnguianoSuarezGrau2019,AnguianoSuarezGrau2026} were combined with homogenization and dimension
reduction, pore-filling extensions and uniform Korn inequalities in $L^p$ were developed
\cite{GahnJagerNeussRadu2026}, and Navier-slip-type effective interface laws were derived for a thin elastic
layer between two fluid regions \cite{GahnNeussRadu2025}.  In these works the layer reduces to a
lower-dimensional limit; here the limit is a bulk Darcy law, and the symmetric-gradient estimates and
oscillating tests are developed in Orlicz spaces with a nonlinear friction.  At the macroscopic level,
S\'anchez-Vargas--Vald\'es-Parada--Lasseux \cite{SanchezVargasValdesParadaLasseux2022} derive an unsteady
generalized-Newtonian model with inertia and slip.

At fixed scale, existence for non-Newtonian fluids with generalized Orlicz growth is developed in
\cite{GwiazdaSwierczewskaGwiazdaWroblewska2010}, and under frictional boundary conditions in reflexive Orlicz
spaces in \cite{MigorskiPaczka2018}.  Our analysis uses two-scale convergence \cite{Allaire1992,Nguetseng1989}
and its periodic-surface version \cite{AllaireDamlamianHornung1996,NeussRadu1992,NeussRadu1996}; the periodic
unfolding method in perforated domains, with its boundary unfolding operator,
\cite{CioranescuDamlamianDonatoGrisoZaki2012} provides an equivalent framework (Remark~\ref{rem:unfolding}).  In
nonreflexive Musielak--Orlicz spaces, without $\Delta_2$ or $\nabla_2$ conditions, the two-scale homogenization of
nonlinear elliptic systems was carried out by Bul\'i\v{c}ek--Gwiazda--Kalousek--\'Swierczewska-Gwiazda
\cite{BulicekGwiazdaKalousekSwierczewskaGwiazda2019} under a log-H\"older-type continuity of the $N$-function; here
the $\Delta_2$ conditions on $\Phi$ and $\widetilde\Phi$ keep the spaces reflexive, as required by the Korn,
Bogovskii and Helmholtz tools used below.  A
Darcy-scaled trace mechanism appears, in the Hilbert setting, in the perforated-domain Steklov analyses
\cite{Douanla2010,Douanla2013Steklov}, and a nonlinear Orlicz--Steklov surface-to-volume problem is studied in
\cite{DouanlaNnangTachagoWoukeng2026}, where the surface terms are approximated by cell averages, whereas here
the velocity keeps its microscopic dependence and its surface limit is the microscopic trace of the two-scale
limit, which enters a nonlinear cell problem with friction.

\subsection{Hypotheses}
\label{subsec:hypotheses}

In addition to Assumption~\ref{ass:geometry}, we impose the following assumptions on the growth, the
constitutive laws, the force and the scaling.  A measurable function $g:\mathbb R^N\to\mathbb R$ is
$Y$-periodic if $g(y+k)=g(y)$ for a.e. $y\in\mathbb R^N$ and every $k\in\mathbb Z^N$; periodic functions are
identified with their restrictions to $Y$, and $C_{\rm per}(Y)$, $C^\infty_{\rm per}(Y)$ and
$L^\infty_{\rm per}(Y)$ denote the continuous, smooth and essentially bounded $Y$-periodic functions.  Since
$\overline T\subset Y$, functions on $\partial T$ extend periodically to $\partial T+\mathbb Z^N$ without
any compatibility condition, and we simply write $C(\partial T)$ and $L^\infty(\partial T)$.  The subscript
``${\rm per},0$'' denotes periodicity together with zero mean.  We write $\R_{\rm sym}^{N\times N}$ for the
symmetric matrices and $\mathbb T_y:=\{z\in\R^N:\ z\cdot n(y)=0\}$ for the tangent space at
$y\in\partial T$.  For $1\le s<N$ we write $s^*:=Ns/(N-s)$ for the Sobolev exponent, and $s^*:=\infty$ if
$s\ge N$.

\begin{assumption}\label{ass:orlicz}
The function $\Phi$ is given by $\Phi(t)=\int_0^t\varphi(s)\,\dd s$, where $\varphi:[0,\infty)\to[0,\infty)$ is nondecreasing,
right-continuous and positive on $(0,\infty)$, and there exist $1<p^-\le p^+<\infty$ such that
\begin{equation}\label{indices}
 p^-\le \frac{t\varphi(t)}{\Phi(t)}\le p^+\qquad(t>0).
\end{equation}
\end{assumption}

Condition \eqref{indices} implies that $\Phi$ is an $N$-function, with $\varphi(0^+)=0$ and
$\varphi(t)\to\infty$, and that $\Phi$ and its complementary function
$\widetilde\Phi(s):=\sup_{t\ge0}\{st-\Phi(t)\}$ satisfy the global $\Delta_2$ condition
\cite[Theorem~II.3.3]{RaoRen1991}:
\begin{equation}\label{Delta2}
\exists\, C_{\Delta_2}\ge1:\qquad \Upsilon(2t)\le C_{\Delta_2}\Upsilon(t)\qquad\forall\,t\ge0,\ \Upsilon\in\{\Phi,\widetilde\Phi\}.
\end{equation}
The upper bound in \eqref{indices} gives the doubling property of $\Phi$ with $C_{\Delta_2}=2^{p^+}$, and
$p^->1$ gives that of $\widetilde\Phi$.  Property \eqref{Delta2} is the one used throughout: it makes
$L^\Phi(D)$ and $L^{\widetilde\Phi}(D)$ reflexive and mutually dual for every measurable $D\subset\mathbb R^N$ of finite measure.

On sets of finite measure, the integrability provided by $\Phi$ depends only on its behaviour at infinity.
We record it through the lower growth exponent at infinity
\begin{equation}\label{pinfty}
 p^-_\infty:=\sup\big\{s\ge1:\ \exists\,C>0,\ t^s\le C\big(1+\Phi(t)\big)\ \ \forall\,t\ge0\big\}.
\end{equation}
The set in \eqref{pinfty} is an interval.  By \eqref{indices}, $\Phi(1)t^{p^-}\le\Phi(t)\le\Phi(1)t^{p^+}$ for
$t\ge1$, so that $p^-\le p^-_\infty\le p^+$, and the supremum is attained when $p^-_\infty=p^-$.  The exponent
$p^-_\infty$ exceeds $p^-$ when $\Phi$ grows more slowly near the origin than at infinity.

\begin{assumption}
\label{ass:stress}
The mapping $\mathbb S:\R^N\times\R_{\rm sym}^{N\times N}\to \R_{\rm sym}^{N\times N}$ is $Y$-periodic and
measurable in the first variable, continuous in the second, and $\mathbb S(y,0)=0$.  There are constants
$c_S,C_S,a_S>0$ such that, for almost every $y\in Y$ and all $E,G\in\R_{\rm sym}^{N\times N}$,
\begin{align}
\label{Sgrowth}
\widetilde\Phi\!\left(\frac{\abs{\mathbb S(y,E)}}{c_S}\right)
&\le C_S\Phi(\abs E),\\
\label{Scoercive}
\mathbb S(y,E):E&\ge a_S\Phi(\abs E),\\
\label{Smonotone}
\big(\mathbb S(y,E)-\mathbb S(y,G)\big):(E-G)&\ge0.
\end{align}
\end{assumption}

\begin{assumption}
\label{ass:friction}
For $\sigma_y$-almost every $y\in\partial T$, the mapping $z\mapsto\mathcal B(y,z)$ maps $\mathbb T_y$ into
$\mathbb T_y$, is continuous, and satisfies $\mathcal B(y,0)=0$; the dependence on $y$ is measurable.  There
are constants $c_B,C_B,a_B>0$ such that, for $\sigma_y$-almost every $y\in\partial T$ and all $z,w\in\mathbb T_y$,
\begin{align}
\label{Bgrowth}
\widetilde\Phi\!\left(\frac{\abs{\mathcal B(y,z)}}{c_B}\right)
&\le C_B\Phi(\abs z),\\
\label{Bcoercive}
\mathcal B(y,z)\cdot z&\ge a_B\Phi(\abs z),\\
\label{Bmonotone}
\big(\mathcal B(y,z)-\mathcal B(y,w)\big)\cdot(z-w)&\ge0.
\end{align}
\end{assumption}

The measurability in $y$ is understood through the extension
$\widehat{\mathcal B}(y,z):=\mathcal B(y,P_\tau(y)z)$, $(y,z)\in\partial T\times\R^N$, which is required to be a
Carath\'eodory function, that is, measurable in $y$ and continuous in $z$.  Since $P_\tau$ is continuous on
$\partial T$, the compositions $\mathcal B(x/\eps,(u_\eps)_\tau)$ and $\mathcal B(y,\gamma_0^yu_0)$ used below are
then measurable.

\begin{assumption}\label{ass:force-exponent}
The body force satisfies $f\in L^{\widetilde\Phi}(\Omega)^N$, and
\begin{equation}\label{NSexponent}
 p^-_\infty>\frac{3N}{N+2},
\end{equation}
that is, $p^-_\infty>3/2$ if $N=2$ and $p^-_\infty>9/5$ if $N=3$.
\end{assumption}

\begin{assumption}\label{ass:kappa}
\begin{equation}\label{finalkappacondition}
 \kappa>\kappa_c:=1+N\left(\frac{2}{p^-_\infty}-1\right)_+,
\end{equation}
where $(t)_+:=\max\{t,0\}$ for $t\in\R$.
\end{assumption}

Since \eqref{NSexponent} and \eqref{finalkappacondition} are strict inequalities and the set in
\eqref{pinfty} is an interval, we fix once and for all an exponent $\underline p\in[p^-,p^-_\infty]$ in that set
such that
\begin{equation}\label{underlinep}
 t^{\underline p}\le C\big(1+\Phi(t)\big)\quad(t\ge0),\qquad
 \underline p>\frac{3N}{N+2},\qquad
 \kappa>1+N\left(\frac2{\underline p}-1\right)_+;
\end{equation}
the last condition is required only under Assumption~\ref{ass:kappa}, and $\underline p=p^-_\infty$ whenever
the supremum in \eqref{pinfty} is attained.

\begin{assumption}\label{ass:strictmono}
For almost every $y\in Y$ and all $E,G\in\R_{\rm sym}^{N\times N}$ with $E\ne G$,
\begin{equation}\label{strictmonocell}
 [\mathbb S(y,E)-\mathbb S(y,G)]:(E-G)>0.
\end{equation}
\end{assumption}

The roles of these assumptions are distinct.  Microscopic solvability uses
Assumptions~\ref{ass:geometry}--\ref{ass:force-exponent} and holds for every $\kappa\in\R$; the limit passage
uses in addition Assumption~\ref{ass:kappa}; compactness and the identification of the nonlinear limits use
the monotonicity conditions \eqref{Smonotone} and \eqref{Bmonotone}, and not their strict versions; Assumption~\ref{ass:strictmono},
which concerns the bulk stress only, is used for the uniqueness of the cell velocity, of the effective
permeability and of the macroscopic pressure; Subsection~\ref{subsec:relaxations} describes what remains under
weaker hypotheses.

\subsection{Main result}
\label{subsec:main-results}

The effective law is governed by a nonlinear cell problem.  For $\zeta\in\R^N$, let
$(w_\zeta,\pi_\zeta)$ be the solution, with $w_\zeta$ and $\pi_\zeta$ $Y$-periodic and
$\int_{Y^*}\pi_\zeta\,\dd y=0$, of
\begin{equation}
\label{introcellproblem}
\left\{
\begin{aligned}
-\diver_y\mathbb S(y,\D_yw_\zeta)+\nabla_y\pi_\zeta&=\zeta &&\text{in }Y^*,\\
\diver_yw_\zeta&=0 &&\text{in }Y^*,\\
w_\zeta\cdot n&=0 &&\text{on }\partial T,\\
[\mathbb S(y,\D_yw_\zeta)n]_\tau+\mathcal B(y,\gamma_0w_\zeta)&=0 &&\text{on }\partial T,
\end{aligned}
\right.
\end{equation}
understood in the weak sense of \eqref{cellfullweak}; it admits a solution, which is unique under
Assumption~\ref{ass:strictmono} (Theorem~\ref{thm:cellproblem} and Corollary~\ref{cor:cellpressurezeta}).  The
effective permeability is
\begin{equation}
\label{introK}
 \mathcal K(\zeta):=\int_{Y^*}w_\zeta(y)\,\dd y,\qquad \zeta\in\R^N.
\end{equation}
For $0<\eps<\eps_0$, with $\eps_0$ fixed in Theorem~\ref{thm:boundarylayer}, problem \eqref{Peps} admits at
least one weak solution $u_\eps$ in the sense of Definition~\ref{def:weakvelocity}
(Theorem~\ref{thm:microexistence}).  Below, a tilde denotes zero extension from $\Omega_\eps$ to $\Omega$,
$\chi_\eps$ and $\chi_{Y^*}$ are the characteristic functions of $\Omega_\eps$ in $\Omega$ and of $Y^*$ in $Y$,
$\mathbb W^\Phi_{\rm per}(\Omega\times Y^*)$ is the space \eqref{productOrliczSpace}, and the volume and surface
two-scale convergences $\overset{2s}{\rightharpoonup}$ and $\overset{2s,\Gamma}{\rightharpoonup}$ are those of
Definitions~\ref{def:weak-two-scale} and~\ref{def:surface-two-scale}.  We write
$\mathcal P^{\widetilde\Phi}_0(\Omega):=\{p\in W^1L^{\widetilde\Phi}(\Omega):\ \int_\Omega p\,\dd x=0\}$.  Convergence
``for the whole family'' means convergence along every sequence in $(0,\eps_0)$ tending to zero, without
extraction.  The following theorem is proved in Subsection~\ref{subsec:proofA}.

\begin{theoremA}
Let Assumptions~\ref{ass:geometry}--\ref{ass:strictmono} hold, and let $(u_\eps)_{0<\eps<\eps_0}$ be any
family of weak solutions of \eqref{Peps}.  Then there exist unique
$u_0\in\mathbb W^\Phi_{\rm per}(\Omega\times Y^*)^N$ and $p_0\in\mathcal P^{\widetilde\Phi}_0(\Omega)$ such that,
as $\eps\to0$ and for the whole family,
\begin{align}
 \widetilde u_\eps&\overset{2s}{\rightharpoonup}\chi_{Y^*}u_0,
 \label{finalconvU}\\
 \eps\widetilde{\nabla u_\eps}&\overset{2s}{\rightharpoonup}\chi_{Y^*}\nabla_yu_0,
 \qquad
 \eps\widetilde{\D u_\eps}\overset{2s}{\rightharpoonup}\chi_{Y^*}\D_yu_0,
 \label{finalconvgrad}\\
 u_\eps&\overset{2s,\Gamma}{\rightharpoonup}\gamma_0^y u_0,
 \label{finalconvtrace}\\
 \chi_\eps\,\mathbb S\!\left(\frac{x}{\eps},\eps\D u_\eps\right)
 &\overset{2s}{\rightharpoonup}\chi_{Y^*}\mathbb S(y,\D_yu_0),
 \label{finalconvstress}\\
 \mathcal B\!\left(\frac{x}{\eps},(u_\eps)_\tau\right)
 &\overset{2s,\Gamma}{\rightharpoonup}\mathcal B(y,\gamma_0^y u_0),
 \label{finalconvfriction}
\end{align}
and
\begin{equation}
\label{finalu0cell}
 u_0(x,\cdot)=w_{f(x)-\nabla p_0(x)}\qquad\text{for a.e. }x\in\Omega.
\end{equation}
The permeability map \eqref{introK} is continuous and strictly monotone and satisfies the two-sided Orlicz
bounds $\Phi(c|\mathcal K(\zeta)|)\le C\widetilde\Phi(C|\zeta|)$ and
$\zeta\cdot\mathcal K(\zeta)\ge c\,\widetilde\Phi(c|\zeta|)$ for all $\zeta\in\R^N$.  The filtration velocity
$U:=\int_{Y^*}u_0\,\dd y$ belongs to $L^\Phi(\Omega)^N$ and satisfies
\begin{equation}
\label{finalDarcyboxed}
 U=\mathcal K(f-\nabla p_0),\qquad \diver U=0\ \text{ in }\Omega,\qquad U\cdot\nu=0\ \text{ on }\partial\Omega;
\end{equation}
equivalently, $p_0$ is the unique solution in $\mathcal P^{\widetilde\Phi}_0(\Omega)$ of
\begin{equation}
\label{Darcyweakpressure}
 \int_\Omega
 \mathcal K(f-\nabla p_0)\cdot\nabla\omega\,\dd x=0
 \qquad\forall\,\omega\in W^1L^{\widetilde\Phi}(\Omega).
\end{equation}
Finally, the microscopic pressure $p_\eps$, extended to each hole by its mean over the surrounding annulus and
normalized to zero mean in $\Omega$, converges strongly to $p_0$ in $L^{s}(\Omega)$ for some $s>1$
(Theorem~\ref{thm:pressureconv}).
\end{theoremA}

From the viewpoint of mathematical physics, Theorem~A yields a nonlinear effective transport law in which
the bulk non-Newtonian dissipation and the microscopic wall friction are coupled in the single cell problem
\eqref{introcellproblem}, so that the friction law enters the permeability map itself rather than an
effective boundary condition.  For the Newtonian benchmark $\Phi(t)=t^2/2$, $\mathbb S(y,E)=2\mu(y)E$ and
$\mathcal B(y,z)=\beta(y)z$, the cell problem is linear and $\mathcal K(\zeta)=K_\beta\zeta$ for a symmetric
positive definite matrix $K_\beta$ depending on $\mu$ and $\beta$, a slip-dependent Darcy law in the critical
scaling, to be compared with the upscaled models of \cite{LasseuxValdesParada2024}.  The map $\mathcal K$ also
inherits the symmetries of the pore space and, for variational laws, is the gradient of a convex potential
which decreases as the wall friction increases and lies between the no-slip and perfect-slip potentials, so that
$K_\beta$ lies between the no-slip and perfect-slip permeability matrices (Remarks~\ref{rem:Kpotential}
and~\ref{rem:Kfriction}).

With respect to the classical setting, three features change the compactness mechanism: the constitutive
law has general Orlicz growth with periodic coefficients, the pore walls carry a nonlinear tangential
friction, and, without the no-slip condition, the energy controls only the scaled symmetric gradient
$\eps\D u_\eps$, so that the full gradient must be recovered by a Korn inequality uniform in $\eps$.  As in
the classical Darcy scaling, the leading-order velocity keeps a microscopic $y$-dependence; here, in
addition, its trace on the pore walls survives through the critical surface measure $\eps\,\dd\sigma_x$.

These features are handled by the following new ingredients.  First, a uniform Orlicz--Korn--Poincar\'e inequality
in which impermeability and tangential friction remove the rigid motions on each cell, completed near
$\partial\Omega$ by a hole-filling argument on the fixed domain (Corollary~\ref{cor:globalkorn}).  Second, a
Darcy-scaled two-scale trace theorem in Orlicz spaces, which identifies the surface limit with the
microscopic trace $\gamma_0^yu_0$ through a Green formula combining the volume and scaled-gradient limits
(Theorem~\ref{thm:darcyscaledtrace}).  Third, exactly solenoidal oscillating recovery tests preserving
the pore-wall trace (Theorem~\ref{thm:solenoidaltests}), which yield the explicit inertial threshold
$\kappa_c$; it depends only on the dimension and on the lower growth exponent at infinity $p^-_\infty$ of the rheology, and not on
the specific constitutive law, the friction law or the pore geometry, and it is sharp at the level of the
Orlicz energy estimates (Proposition~\ref{prop:kappathreshold} and Remark~\ref{rem:kappasharp}).  In particular
$\kappa_c=1$ whenever $p^-_\infty\ge2$, for instance for Carreau laws with positive infinite-shear viscosity,
consistently with the critical Reynolds number independent of the flow index found in
\cite{AnguianoBonnivardSuarezGrau2026}.  In the Newtonian case, $\kappa_c=1$ is exactly the critical scaling of
\cite{MarusicPalokaMikelic2000}, for which, with no-slip conditions and small forces, inertia enters the cell
problem; in the critical regime $\kappa=\kappa_c$ an analogous inertial cell problem with friction is therefore
expected.  Fourth, a coupled bulk--surface Minty principle, closed by a single limit
energy equality, which identifies the bulk stress and the wall friction simultaneously
(Theorem~\ref{thm:energyclosure}).  Two Orlicz de Rham arguments then reconstruct the microscopic and
macroscopic pressures from the limit system, and a Tartar-type restriction operator with modular Orlicz bounds
yields the strong convergence of the extended microscopic pressure (Theorem~\ref{thm:pressureconv}).  As
complements, the macroscopic problem \eqref{Darcyweakpressure} is well posed independently of the homogenization
process (Corollary~\ref{cor:Darcywellposed}), uniform monotonicity yields strong two-scale convergence of the
scaled symmetric gradient and of the wall trace (Proposition~\ref{prop:corrector}), and without strict
monotonicity the Darcy law becomes an inclusion for a maximal monotone permeability (Remark~\ref{rem:nonstrict}).

\subsection{Benchmark constitutive models covered by the hypotheses}
\label{subsec:benchmarks}

The hypotheses cover a broad family of generalized Newtonian bulk laws together with nonlinear finite-slip wall laws.  Throughout this subsection, $t\ge0$, $y\in Y$ and $E\in\R^{N\times N}_{\rm sym}$ in the bulk laws, and $y\in\partial T$ and $z\in\mathbb T_y$ in the friction laws.  A natural class is obtained directly from the $N$-function in Assumption~\ref{ass:orlicz}.  Let $a\in L^\infty_{\rm per}(Y)$ and $b\in L^\infty(\partial T)$ satisfy
\[
 0<a_*\le a(y)\le a^*,\qquad 0<b_*\le b(y)\le b^*
\]
almost everywhere, and set, with the value $0$ at the origin,
\[
 \mathbb S_\Phi(y,E)=a(y)\frac{\varphi(|E|)}{|E|}E,
 \qquad
 \mathcal B_\Phi(y,z)=b(y)\frac{\varphi(|z|)}{|z|}z.
\]
The second mapping preserves the tangent space $\mathbb T_y$.  By \eqref{indices},
\[
 \mathbb S_\Phi(y,E):E\ge a_*p^-\Phi(|E|),\qquad
 \mathcal B_\Phi(y,z)\cdot z\ge b_*p^-\Phi(|z|).
\]
Choosing $c_S\ge a^*$ and $c_B\ge b^*$, the complementary-function identity
\[
 \widetilde\Phi(\varphi(t))=t\varphi(t)-\Phi(t)\le (p^+-1)\Phi(t)
\]
gives \eqref{Sgrowth} and \eqref{Bgrowth}.  Moreover, these maps are the gradients of the convex radial potentials $a(y)\Phi(|E|)$ and $b(y)\Phi(|z|)$, respectively; hence they are monotone, and they are strictly monotone whenever the corresponding radial potentials are strictly convex.  Thus every constitutive pair generated by an admissible $\Phi$ in this way satisfies Assumptions~\ref{ass:stress}--\ref{ass:friction}.  Since Assumptions~\ref{ass:force-exponent} and~\ref{ass:kappa} involve only $p^-_\infty$, we record for each law
below the indices $p^-$, $p^+$ and the exponent $p^-_\infty$; Theorem~A then applies whenever $p^-_\infty>3N/(N+2)$ and
$\kappa>1+N(2/p^-_\infty-1)_+$.  The following familiar laws are included.

\begin{enumerate}[label=\textup{(\roman*)}]
\item \emph{Newtonian viscosity with linear Navier slip.}  For
\[
 \Phi(t)=\frac{t^2}{2},\qquad
 \mathbb S(y,E)=2\mu(y)E,\qquad
 \mathcal B(y,z)=\beta(y)z,
\]
with positive bounded periodic coefficients $\mu$ and $\beta$, one has $p^-=p^+=p^-_\infty=2$ and hence $\kappa_c=1$.  Uniformly positive definite periodic linear tensors in the bulk, and uniformly positive definite tangential tensors on $\partial T$, are covered as well.

\item \emph{Power-law fluids with nonlinear power-law slip.}  For $1<p<\infty$,
\[
 \Phi(t)=\frac{t^p}{p},\qquad
 \mathbb S(y,E)=\mu(y)|E|^{p-2}E,\qquad
 \mathcal B(y,z)=\beta(y)|z|^{p-2}z,
\]
one has $p^-=p^+=p^-_\infty=p$; in particular $\kappa_c=1$ for $p\ge2$.
Power-law filtration laws form one of the classical nonlinear Darcy benchmarks.

\item \emph{Mixed-power and Ladyzhenskaya-type rheologies.}  Let $1<p<q<\infty$ and
\[
 \Phi(t)=\frac{t^p}{p}+\frac{t^q}{q}.
\]
Then the constitutive pair
\[
 \mathbb S(y,E)=a(y)\big(|E|^{p-2}+|E|^{q-2}\big)E,
 \qquad
 \mathcal B(y,z)=b(y)\big(|z|^{p-2}+|z|^{q-2}\big)z
\]
satisfies the hypotheses with $p^-=p$ and $p^+=q$, while $p^-_\infty=q$ because $\Phi(t)\ge t^q/q$; the conditions of
Theorem~A thus involve only $q$, for every $p\in(1,q)$, and $\kappa_c=1$ when $q\ge2$.  The choice $p=2<q$ gives the standard Newtonian-plus-power-law, or Ladyzhenskaya-type, structure.  More generally, finite sums of powers are covered by taking $\Phi$ to be the corresponding sum of power $N$-functions.

\item \emph{Carreau--Yasuda-type laws.}  Let $a_C>0$, $\lambda>0$, $r>1$ and
$0\le\mu_\infty<\mu_0<\infty$, and set
\[
 \mu_C(t)=\mu_\infty+(\mu_0-\mu_\infty)\big(1+(\lambda t)^{a_C}\big)^{(r-2)/a_C},
 \qquad
 \Phi_C(t)=\int_0^t \mu_C(s)s\,\dd s.
\]
A direct differentiation gives
\[
 \frac{d}{dt}\bigl[t\mu_C(t)\bigr]
 =\mu_\infty+(\mu_0-\mu_\infty)
 \bigl(1+(\lambda t)^{a_C}\bigr)^{(r-2)/a_C-1}
 \bigl[1+(r-1)(\lambda t)^{a_C}\bigr]>0,
\]
so the associated radial constitutive map is strictly monotone.  Since, for every $\varsigma\ge0$, the number $1+(r-1)\varsigma$ lies between
$\min\{1,r-1\}(1+\varsigma)$ and $\max\{1,r-1\}(1+\varsigma)$, this derivative satisfies
\[
 \min\{1,r-1\}\,\mu_C(t)\le\frac{d}{dt}\bigl[t\mu_C(t)\bigr]\le\max\{1,r-1\}\,\mu_C(t),
\]
also when $\mu_\infty>0$.  With $\varphi_C(t):=t\mu_C(t)$, the identity
$t\varphi_C(t)=\int_0^t\bigl[\varphi_C(s)+s\varphi_C'(s)\bigr]\dd s$ then gives
\[
 \min\{2,r\}\le\frac{t\varphi_C(t)}{\Phi_C(t)}\le\max\{2,r\}\qquad(t>0),
\]
so that one may take $p_C^-=\min\{2,r\}$ and $p_C^+=\max\{2,r\}$; in particular
$\Phi_C,\widetilde\Phi_C\in\Delta_2$.  At infinity, $\Phi_C$ grows like $t^r$ if $\mu_\infty=0$ and like
$t^{\max\{2,r\}}$ if $\mu_\infty>0$, so that $p^-_\infty=r$ in the first case and $p^-_\infty=\max\{2,r\}$ in the
second.  Consequently, for a shear-thinning law with $\mu_\infty>0$, Theorem~A applies for every $r\in(1,2)$,
with $\kappa_c=1$, and for $\mu_\infty=0$ it applies whenever $r>3N/(N+2)$.  Therefore
\[
 \mathbb S(y,E)=a(y)\mu_C(|E|)E,
 \qquad
 \mathcal B(y,z)=b(y)\mu_C(|z|)z
\]
satisfies the structural hypotheses.  Carreau and Carreau--Yasuda porous homogenization provides an important nonlinear Darcy comparison class; see \cite{AnguianoBonnivardSuarezGrau2022,AnguianoBonnivardSuarezGrau2025,BourgeatGipoulouxMarusicPaloka2003,LuQian2024}.

\item \emph{Regularized Bingham--Papanastasiou rheology.}  Let
$\mu,\tau_0\in L^\infty_{\rm per}(Y)$ and $\beta\in L^\infty(\partial T)$ satisfy
\[
 0<\mu_*\le\mu(y)\le\mu^*,\qquad
 0\le\tau_0(y)\le\tau^*,\qquad
 0<\beta_*\le\beta(y)\le\beta^*,
\]
and let $m>0$.  Define
\[
 \mathbb S(y,E)=2\mu(y)E+\tau_0(y)\big(1-e^{-m|E|}\big)\frac{E}{|E|},
 \qquad
 \mathcal B(y,z)=\beta(y)z,
\]
with the second bulk term equal to $0$ at $E=0$.  Since $1-e^{-mt}\le mt$,
$|\mathbb S(y,E)|\le(2\mu^*+m\tau^*)|E|$, while
\[
 \mathbb S(y,E):E\ge 2\mu_*|E|^2,
 \qquad
 \mathcal B(y,z)\cdot z\ge\beta_*|z|^2.
\]
The radial stress magnitude $2\mu(y)t+\tau_0(y)(1-e^{-mt})$ has derivative
$2\mu(y)+m\tau_0(y)e^{-mt}>0$.  Hence both constitutive maps satisfy the growth, coercivity and strict-monotonicity requirements with $\Phi(t)=t^2/2$, so that $p^-_\infty=2$ and $\kappa_c=1$.  Thus regularized viscoplastic models are included.

\item \emph{Genuinely non-polynomial Orlicz laws.}  For $p>1$ and $\alpha>0$, set
\[
 \Phi_{p,\alpha}(t)=\int_0^t s^{p-1}\big(1+\log(1+s)\big)^\alpha\,\dd s.
\]
With $L(t):=1+\log(1+t)$, the bound $\Phi_{p,\alpha}(t)\le t^pL(t)^\alpha/p$, which holds because $L$ is
increasing, and the bound $\frac{d}{dt}[t^pL(t)^\alpha]\le(p+\alpha)t^{p-1}L(t)^\alpha$, integrated from $0$ to $t$,
give $p\le t\Phi'_{p,\alpha}(t)/\Phi_{p,\alpha}(t)\le p+\alpha$ for $t>0$.
Thus one may take $p^-=p$ and $p^+=p+\alpha$, and $p^-_\infty=p$ because $\Phi_{p,\alpha}(t)\le t^pL(t)^\alpha/p$
is $o(t^s)$ as $t\to\infty$ for every $s>p$; in particular
$\Phi_{p,\alpha}$ and its complementary function satisfy $\Delta_2$.  The laws
\[
 \mathbb S(y,E)=a(y)|E|^{p-2}\big(1+\log(1+|E|)\big)^\alpha E,
 \qquad
 \mathcal B(y,z)=b(y)|z|^{p-2}\big(1+\log(1+|z|)\big)^\alpha z
\]
therefore provide an Orlicz, non-power benchmark.
\end{enumerate}

The preceding list concerns finite coercive wall friction.  The no-slip condition is an infinite-friction limit of Assumption~\ref{ass:friction}; perfect slip $\mathcal B\equiv0$ can be included under the relaxation described in Remark~\ref{rem:perfectslip}; ideal Bingham and Herschel--Bulkley laws, whose yield graphs are multivalued, are approximated by the regularized model (v).

\subsection{Organization of the paper}
\label{subsec:organization}
Section~\ref{sec:preliminaries} collects the functional setting, two-scale convergence in volume and on
surfaces, and the cell and perforated-domain estimates.  Section~\ref{sec:microscopic-estimates} proves
coercivity, existence and the uniform Darcy-scale bounds.  Section~\ref{sec:homogenization} carries out the
homogenization and identifies the nonlinear limits.  Section~\ref{sec:effectivedarcy} treats the pressures, the
cell problem and the permeability map, proves Theorem~A, and discusses relaxations of the hypotheses.

\section{Functional setting and preliminary results}
\label{sec:preliminaries}

This section collects the Orlicz tools, the volume and surface two-scale convergences and the
reference-cell estimates used in the sequel.  Two-scale limits, written $\lim_{\eps\to0}$, are taken along an arbitrary sequence
$\eps\to0$.  We denote by $\mathscr E\subset(0,1)$ the set of its terms; uniform bounds are suprema over
$\eps\in\mathscr E$, and ``after extraction'' refers to a subsequence of $\mathscr E$.  A number
$\eps_0\in(0,1]$ is fixed in Section~\ref{sec:microscopic-estimates}, and from
Section~\ref{sec:homogenization} on we take $\mathscr E\subset(0,\eps_0)$.

\subsection{Functional setting}
\label{subsec:functional-setting}

Let $(Q,\mathcal A,\mu)$ be a finite measure space.  Under Assumption~\ref{ass:orlicz}, define
\begin{equation}\label{OrliczHeartDefinition}
L^\Phi(Q,\mu)=\left\{v:Q\to\mathbb R\text{ measurable}:\
\lim_{\delta\to0^+}\int_Q\Phi(\delta|v|)\,\dd\mu=0\right\}.
\end{equation}
Since $\Phi\in\Delta_2$, this agrees with the usual Orlicz class.  The Luxemburg norm is
\begin{equation}\label{LuxemburgNorm}
 \|v\|_{\Phi,Q,\mu}:=
 \inf\left\{a>0:\int_Q\Phi\!\left(\frac{|v|}{a}\right)\dd\mu\le1\right\},
\end{equation}
and the generalized H\"older inequality reads
\begin{equation}\label{OrliczHolder}
 \int_Q|vw|\,\dd\mu\le2\|v\|_{\Phi,Q,\mu}\|w\|_{\widetilde\Phi,Q,\mu}.
\end{equation}
Because $\Phi,\widetilde\Phi\in\Delta_2$, the complementary spaces are separable and reflexive and
$(L^{\widetilde\Phi}(Q,\mu))'=L^\Phi(Q,\mu)$ under the integral pairing \cite[Corollaries~IV.1.9 and~IV.1.12]{RaoRen1991}.  Modular and norm convergence are equivalent.
For vector- or matrix-valued $v$, $\int_Q\Phi(|v|)\,\dd\mu$ is taken with the Euclidean (Frobenius) norm
$|v|$; by $\Delta_2$ it is equivalent, up to constants inside and outside $\Phi$, to the sum of the
componentwise modulars.  The following elementary consequences of \eqref{indices} are used throughout.

\begin{lemma}\label{lem:orlicztools}
Let $(Q,\mathcal A,\mu)$ be a finite measure space.
\begin{enumerate}[label=\textup{(\roman*)}]
\item For all $t\ge0$,
\begin{equation}\label{scalingPhi}
\lambda^{p^-}\Phi(t)\le\Phi(\lambda t)\le\lambda^{p^+}\Phi(t)\quad(\lambda\ge1),
\qquad
\lambda^{p^+}\Phi(t)\le\Phi(\lambda t)\le\lambda^{p^-}\Phi(t)\quad(0<\lambda\le1).
\end{equation}
\item For all $v\in L^\Phi(Q,\mu)$,
\begin{equation}\label{modularnorm}
\min\big\{\|v\|_{\Phi,Q,\mu}^{p^-},\|v\|_{\Phi,Q,\mu}^{p^+}\big\}
\le\int_Q\Phi(|v|)\,\dd\mu
\le\max\big\{\|v\|_{\Phi,Q,\mu}^{p^-},\|v\|_{\Phi,Q,\mu}^{p^+}\big\}.
\end{equation}
\item For all $s,t\ge0$,
\begin{equation}\label{youngPhi}
st\le\Phi(s)+\widetilde\Phi(t),
\qquad
\widetilde\Phi(\varphi(t))=t\varphi(t)-\Phi(t)\le(p^+-1)\Phi(t).
\end{equation}
\item Let $Q\subset\mathbb R^m$ be a bounded connected Lipschitz domain.  There are $C,c>0$, depending only on
$Q$ and $p^+$, such that, with $\bar v:=|Q|^{-1}\int_Qv$,
\begin{equation}\label{modularPoincare}
\int_Q\Phi(|v-\bar v|)\,\dd x\le C\int_Q\Phi(c|\nabla v|)\,\dd x
\qquad\forall\,v\in W^1L^\Phi(Q).
\end{equation}
If $Q\subset\mathbb R^m$ is a bounded open set with diameter $d_Q$, then
\begin{equation}\label{modularPoincarezero}
\int_Q\Phi(|v|)\,\dd x\le \int_Q\Phi(d_Q|\nabla v|)\,\dd x
\qquad\forall\,v\in W^1_0L^\Phi(Q).
\end{equation}
\end{enumerate}
\end{lemma}

\begin{proof}
By \eqref{indices}, $t\mapsto\Phi(t)t^{-p^-}$ is nondecreasing and $t\mapsto\Phi(t)t^{-p^+}$ is
nonincreasing, which is \eqref{scalingPhi}; \eqref{modularnorm} follows by applying \eqref{scalingPhi} to
$v/\|v\|_{\Phi,Q,\mu}$.  The first inequality in \eqref{youngPhi} is the definition of $\widetilde\Phi$, with
equality for $t=\varphi(s)$, which gives the identity; the bound follows from \eqref{indices}.  For
\eqref{modularPoincare}, by density of
smooth functions in $W^1L^\Phi(Q)$ under $\Delta_2$ \cite[Theorem~2.3]{DonaldsonTrudinger1971}, it suffices to consider $v\in
C^1(\overline Q)$.  Consider first the unit cube $Q_0=(0,1)^m$.  For $1\le j\le m$, let $M_j$ denote the average
with respect to the $j$-th coordinate over $(0,1)$, the other coordinates being fixed.  The operators $M_j$
commute, and by Jensen's inequality $\int_{Q_0}\Phi(M_jg)\,\dd x\le\int_{Q_0}\Phi(g)\,\dd x$ for $g\ge0$.  The
fundamental theorem of calculus in the $j$-th variable gives $|(I-M_j)W|\le M_j|D_jW|$ for
$W\in C^1(\overline{Q_0})$, and the telescoping identity
$W-\int_{Q_0}W\,\dd x=\sum_{j=1}^mM_1\cdots M_{j-1}(I-M_j)W$ yields
$|W-\int_{Q_0}W|\le\sum_{j=1}^mM_1\cdots M_j|D_jW|$.  Convexity of $\Phi$ and Jensen's inequality for the successive
averages give
\begin{equation}\label{cubePoincare}
\int_{Q_0}\Phi\Big(\Big|W-\int_{Q_0}W\,\dd x\Big|\Big)\dd x
\le\frac1m\sum_{j=1}^m\int_{Q_0}\Phi(m|D_jW|)\,\dd x
\le\int_{Q_0}\Phi(m|\nabla W|)\,\dd x.
\end{equation}
For a general $Q$, cover $Q$ by finitely many open sets $V_1,\dots,V_M\subset Q$, each the image of $Q_0$
under a bi-Lipschitz map $\Theta_\ell$: cubes contained in $Q$ and, near $\partial Q$, the intersections of $Q$
with Lipschitz coordinate cylinders, flattened by the graph maps.  Since $Q$ is connected, the graph whose
vertices are $1,\dots,M$ and whose edges join $\ell$ and $j$ when $|V_\ell\cap V_j|>0$ is connected.  With
$m_\ell:=\int_{Q_0}v\circ\Theta_\ell\,\dd t$, estimate \eqref{cubePoincare} applied to $v\circ\Theta_\ell$, the bounds
on the Jacobians of $\Theta_\ell$ and \eqref{scalingPhi} give
\[
\int_{V_\ell}\Phi(|v-m_\ell|)\,\dd x\le C\int_{V_\ell}\Phi(C|\nabla v|)\,\dd x,
\qquad \ell=1,\dots,M.
\]
If $|V_\ell\cap V_j|>0$, integrating $\Phi(|m_\ell-m_j|/2)\le\frac12\Phi(|v-m_\ell|)+\frac12\Phi(|v-m_j|)$ over
$V_\ell\cap V_j$ bounds $|V_\ell\cap V_j|\,\Phi(|m_\ell-m_j|/2)$ by the right-hand sides of the last estimate for
$\ell$ and $j$.  Following a path from $1$ to $\ell$ in the graph and using \eqref{scalingPhi}, we obtain
$\Phi(|m_\ell-m_1|)\le C\int_Q\Phi(C|\nabla v|)\,\dd x$ for every $\ell$.  Since
$\Phi(|v-m_1|)\le\frac12\Phi(2|v-m_\ell|)+\frac12\Phi(2|m_\ell-m_1|)$ on $V_\ell$, summation over $\ell$ gives
$\int_Q\Phi(|v-m_1|)\,\dd x\le C\int_Q\Phi(C|\nabla v|)\,\dd x$.  Finally, Jensen's inequality gives
$|Q|\,\Phi(|\bar v-m_1|)\le\int_Q\Phi(|v-m_1|)\,\dd x$, and a final convexity argument replaces $m_1$ by $\bar v$,
which proves \eqref{modularPoincare}.

For \eqref{modularPoincarezero}, let $v\in C_c^\infty(Q)$ and write $x=(x_1,x')$.  For fixed $x'$, the function
$v(\cdot,x')$ is supported in an interval $I(x')$ of length at most $d_Q$, and
$|v(x)|\le\int_{I(x')}|D_1v(t,x')|\,\dd t$.  Jensen's inequality with respect to $d_Q^{-1}\dd t$ gives
$\Phi(|v(x)|)\le d_Q^{-1}\int_{I(x')}\Phi(d_Q|D_1v(t,x')|)\,\dd t$; integrating in $x_1\in I(x')$ and then in $x'$
proves \eqref{modularPoincarezero} for smooth $v$, and density in $W^1_0L^\Phi(Q)$ concludes.
\end{proof}

For an open set $Q\subset\mathbb R^m$,
\begin{equation}\label{OrliczSobolevSpace}
 W^1L^\Phi(Q):=\{v\in L^\Phi(Q):D_{x_i}v\in L^\Phi(Q),\ i=1,\ldots,m\},
\end{equation}
with the usual sum norm.  The space $W^1_0L^\Phi(Q)$ is the closure of $C_c^\infty(Q)$ in this norm.
On bounded connected Lipschitz domains the Orlicz--Poincar\'e inequality holds, and
\begin{equation}\label{compactOrliczEmbedding}
 W^1L^\Phi(Q)\hookrightarrow\hookrightarrow L^\Phi(Q);
\end{equation}
see \cite[Proposition~IX.3.7]{RaoRen1991}.  Using the periodicity convention fixed in Subsection~\ref{subsec:hypotheses}, for
$\Upsilon\in\{\Phi,\widetilde\Phi\}$ we define
\[
 L^\Upsilon_{\rm per}(\Omega\times Y)
 :=\bigl\{v|_{\Omega\times Y}:v\in L^\Upsilon_{\rm loc}(\Omega\times\mathbb R^N),\
 v(x,y+k)=v(x,y)\ \text{for a.e. }(x,y),\ k\in\mathbb Z^N\bigr\},
\]
and
\[
 W^1_{\rm per}L^\Upsilon(Y^*)
 :=\overline{\,C^\infty_{\rm per}(Y)|_{Y^*}\,}^{\,W^1L^\Upsilon(Y^*)},
\]
\[
 W^1_{\rm per,0}L^\Upsilon(Y^*)
 :=\left\{w\in W^1_{\rm per}L^\Upsilon(Y^*):\int_{Y^*}w(y)\,\dd y=0\right\}.
\]
Connectedness of $Y^*$ makes $\sum_{i=1}^N\|D_{y_i}w\|_{\Phi,Y^*}$ an equivalent norm on
$W^1_{\rm per,0}L^\Phi(Y^*)$.  For the Darcy-scale limits we also use the product space
\begin{equation}\label{productOrliczSpace}
\mathbb W^\Phi_{\rm per}(\Omega\times Y^*)
:=\left\{v\in L^\Phi(\Omega\times Y^*):\nabla_yv\in L^\Phi(\Omega\times Y^*)^N,\
 v(x,\cdot)\in W^1_{\rm per}L^\Phi(Y^*)\text{ a.e. }x\right\},
\end{equation}
normed by $\|v\|_{\Phi,\Omega\times Y^*}+\|\nabla_yv\|_{\Phi,\Omega\times Y^*}$, which is the norm used below;
for $v\in\mathbb W^\Phi_{\rm per}(\Omega\times Y^*)^N$ the fast gradient belongs to
$L^\Phi(\Omega\times Y^*)^{N\times N}$.

\subsection{Two-scale convergence in volume and on surfaces}
\label{subsec:preliminary-results}

Two-scale convergence in Orlicz spaces is defined through oscillating test functions.  For
$\psi\in L^{\widetilde\Phi}(\Omega;C_{\rm per}(Y))$ set
$\psi^\eps(x)=\psi(x,x/\eps)$.

\begin{definition}\label{def:weak-two-scale}
A family $(v_\eps)_{\eps\in\mathscr E}\subset L^\Phi(\Omega)$ converges weakly two-scale to
$v_0\in L^\Phi_{\rm per}(\Omega\times Y)$ if
\begin{equation}\label{weakTwoScaleDefinition}
 \lim_{\eps\to0}\int_\Omega v_\eps(x)\psi\!\left(x,\frac{x}{\eps}\right)\dd x
 =\int_\Omega\int_Yv_0(x,y)\psi(x,y)\,\dd y\,\dd x
 \qquad\forall\,\psi\in L^{\widetilde\Phi}(\Omega;C_{\rm per}(Y)).
\end{equation}
We write $v_\eps\overset{2s}{\rightharpoonup}v_0$.  The same definition, with $\Phi$ and
$\widetilde\Phi$ interchanged, is used in $L^{\widetilde\Phi}(\Omega)$.  Since the tests include the
functions of $L^{\widetilde\Phi}(\Omega)$ independent of $y$, the uniform boundedness principle shows
that every sequence which converges weakly two-scale is bounded in $L^\Phi(\Omega)$.
\end{definition}

\begin{definition}\label{def:strong-two-scale}
A family $(v_\eps)_{\eps\in\mathscr E}\subset L^\Phi(\Omega)$ converges strongly two-scale to
$v_0\in L^\Phi_{\rm per}(\Omega\times Y)$ if, for every $\eta>0$ and every
$\psi\in L^\Phi(\Omega;C_{\rm per}(Y))$ satisfying
$\|v_0-\psi\|_{\Phi,\Omega\times Y}\le\eta/2$, there exists $\eps_\eta>0$ such that
\begin{equation}\label{strongTwoScaleDefinition}
 \left\|v_\eps-\psi\!\left(\cdot,\frac{\cdot}{\eps}\right)\right\|_{\Phi,\Omega}\le\eta
 \qquad\forall\,\eps\in\mathscr E\cap(0,\eps_\eta).
\end{equation}
We write $v_\eps\overset{2s}{\longrightarrow}v_0$ strongly in $L^\Phi(\Omega)$; the weak and strong notions are those of
\cite[Definition~4.1]{FotsoTachagoNnang2012}.  The same definition, with $\Phi$ and $\widetilde\Phi$ interchanged, is used
in $L^{\widetilde\Phi}(\Omega)$.
\end{definition}

\begin{lemma}\label{lem:weak-strong-two-scale-pairing}
Let $(f_\eps)$ be bounded in $L^{\Phi}(\Omega)$ and
$f_\eps\overset{2s}{\rightharpoonup}f_0$ weakly two-scale.  If
$g_\eps\overset{2s}{\longrightarrow}g_0$ strongly two-scale in
$L^{\widetilde\Phi}(\Omega)$, then
\begin{equation}\label{weakStrongTwoScalePairing}
 \int_\Omega f_\eps(x)g_\eps(x)\,\dd x
 \longrightarrow
 \int_\Omega\int_Y f_0(x,y)g_0(x,y)\,\dd y\,\dd x.
\end{equation}
The same conclusion holds with $\Phi$ and $\widetilde\Phi$ interchanged, and componentwise for
Euclidean and Frobenius products of vector and matrix fields.
\end{lemma}

\begin{proof}
Fix $\eta>0$ and choose
$\psi\in L^{\widetilde\Phi}(\Omega;C_{\rm per}(Y))$ as in Definition~\ref{def:strong-two-scale}, so that
$\|g_0-\psi\|_{\widetilde\Phi,\Omega\times Y}\le\eta/2$ and
$\|g_\eps-\psi^\eps\|_{\widetilde\Phi,\Omega}\le\eta$ for all sufficiently small $\eps$.  Decompose
\[
 \int_\Omega f_\eps g_\eps
 -\int_\Omega\int_Yf_0g_0
 =\int_\Omega f_\eps(g_\eps-\psi^\eps)
 +\left(\int_\Omega f_\eps\psi^\eps-\int_\Omega\int_Yf_0\psi\right)
 +\int_\Omega\int_Yf_0(\psi-g_0).
\]
The middle term tends to zero by weak two-scale convergence, and the first and third terms are bounded
by $C\eta$, uniformly in $\eps$, by \eqref{OrliczHolder}.  Letting first $\eps\to0$ and then
$\eta\downarrow0$ proves \eqref{weakStrongTwoScalePairing}.
\end{proof}

We recall the compactness theorem of \cite[Theorem~4.1]{FotsoTachagoNnang2012}.

\begin{theorem}\label{thm:volume-two-scale-compactness}
Every bounded family in $L^\Phi(\Omega)$ possesses a subsequence converging weakly two-scale in the sense
of Definition~\ref{def:weak-two-scale}.
\end{theorem}

The recovery tests of Section~\ref{sec:homogenization} have a fast profile which is only Lebesgue
integrable in $y$.  Their admissibility is recorded next.

\begin{lemma}\label{lem:Lqadmissible}
Let $q>p^+$, let $K\Subset\Omega$ be compact, and let $V:\Omega\times\mathbb R^N\to\mathbb R$ be $Y$-periodic in $y$,
vanish for $x\notin K$, and satisfy $V\in L^q_{\rm per}(Y;C(\overline\Omega))$, that is, $y\mapsto V(\cdot,y)$ is
strongly measurable from $Y$ into $C(\overline\Omega)$ and $\sup_{x\in\overline\Omega}|V(x,\cdot)|\in L^q(Y)$.  Then
$V^\eps(x):=V(x,x/\eps)$ is measurable and $V^\eps\overset{2s}{\longrightarrow}V$ strongly two-scale in
$L^\Phi(\Omega)$.  The same holds in $L^{\widetilde\Phi}(\Omega)$ when $q>\widetilde p^{\,+}:=(p^-)'$, and
componentwise for vector and matrix fields.
\end{lemma}

\begin{proof}
By \eqref{indices}, $t\mapsto\Phi(t)t^{-p^+}$ is nonincreasing and $t\mapsto\Phi(t)t^{-p^-}$ is
nondecreasing, so that $\Phi(t)\le\Phi(1)\max\{t^{p^-},t^{p^+}\}\le\Phi(1)(1+t^q)$ for $t\ge0$.  Hence, for
every $a>0$ and $\psi\in C_c(\Omega;C_{\rm per}(Y))$, the function $h_a(x,y):=\Phi(|V(x,y)-\psi(x,y)|/a)$
satisfies
\[
 \sup_{x\in\overline\Omega}h_a(x,y)\le\Phi(1)\Big(1+a^{-q}\big(\sup_{x\in\overline\Omega}|V(x,y)|+\|\psi\|_{L^\infty(\Omega\times Y)}\big)^q\Big),
 \qquad y\in Y.
\]
Since $\psi$ is uniformly continuous, $y\mapsto\psi(\cdot,y)$ is continuous into $C(\overline\Omega)$, and
$y\mapsto h_a(\cdot,y)$ is strongly measurable as the composition of a strongly measurable map with a continuous
one.  Hence $h_a\in L^1_{\rm per}(Y;C(\overline\Omega))$, so that $h_a(x,x/\eps)$ is measurable and
\[
 \lim_{\eps\to0}\int_\Omega h_a\!\left(x,\frac{x}{\eps}\right)\dd x
 =\int_\Omega\int_Yh_a(x,y)\,\dd y\,\dd x
\]
by \cite[Corollary~5.4]{Allaire1992} and \cite[Theorem~4]{LukkassenNguetsengWall2002}.  Some continuity
in $x$ uniform in $y$, as encoded in $L^1_{\rm per}(Y;C(\overline\Omega))$, cannot be dropped: continuity of
$x\mapsto h(x,\cdot)$ with values in $L^1(Y)$ alone does not ensure this convergence
\cite[Proposition~5.8]{Allaire1992}.  Since this holds for every $a>0$, the Luxemburg norms converge:
$\|V^\eps-\psi^\eps\|_{\Phi,\Omega}\to\|V-\psi\|_{\Phi,\Omega\times Y}$.  Indeed, each Luxemburg norm is the
infimum of the $a>0$ for which the corresponding modular is at most one, and the limit modular
$a\mapsto\int_\Omega\int_Y\Phi(|V-\psi|/a)\,\dd y\dd x$ is continuous and, where positive, strictly
decreasing, because $\Phi$ is strictly increasing and satisfies $\Delta_2$; pointwise convergence of the
modulars in $a$ therefore implies convergence of the infima.  The same convergence of norms
holds for $\psi\in L^\Phi(\Omega;C_{\rm per}(Y))$ in place of $V$, since then $h_a\in L^1(\Omega;C_{\rm per}(Y))$ \cite[Lemma~5.2]{Allaire1992}.

Fix $\eta>0$ and choose $\psi_0\in C_c(\Omega;C_{\rm per}(Y))$ with
$\|V-\psi_0\|_{\Phi,\Omega\times Y}\le\eta/8$.  For any $\psi\in L^\Phi(\Omega;C_{\rm per}(Y))$ with
$\|V-\psi\|_{\Phi,\Omega\times Y}\le\eta/2$,
\[
 \limsup_{\eps\to0}\|V^\eps-\psi^\eps\|_{\Phi,\Omega}
 \le\lim_{\eps\to0}\|V^\eps-\psi_0^\eps\|_{\Phi,\Omega}
 +\lim_{\eps\to0}\|(\psi_0-\psi)^\eps\|_{\Phi,\Omega}
 \le\frac\eta8+\frac{5\eta}8<\eta,
\]
which is Definition~\ref{def:strong-two-scale}.  For $L^{\widetilde\Phi}(\Omega)$ one uses
$\widetilde\Phi(t)\le\widetilde\Phi(1)(1+t^{q})$ for $q>(p^-)'$, which follows from \eqref{indices} in the
same way.
\end{proof}

The surface tests are sampled through the following Riemann-sum formula.

\begin{lemma}\label{lem:riemann}
Let $a:\Omega\times Y\to\mathbb R$, $Y$-periodic in $y$, and $h:\Omega\times\partial T\to\mathbb R$ be bounded,
measurable in $y$, continuous in $x$ for almost every $y$, and vanishing for $x$ outside a compact set
$K\Subset\Omega$.  Then
\[
 \int_\Omega a\!\left(x,\frac{x}{\eps}\right)\dd x\longrightarrow\int_\Omega\int_Ya\,\dd y\dd x,
 \qquad
 \eps\int_{\Gamma_\eps}h\!\left(x,\frac{x}{\eps}\right)\dd\sigma_x\longrightarrow\int_\Omega\int_{\partial T}h\,\dd\sigma_y\dd x.
\]
\end{lemma}

\begin{proof}
The change of variables $x=\eps(k+y)$ gives
\[
 \eps\int_{\Gamma_\eps}h\!\left(x,\frac{x}{\eps}\right)\dd\sigma_x
 =\int_{\partial T}\sum_{k\in K_\eps}\eps^N h(\eps(k+y),y)\,\dd\sigma_y.
\]
For almost every $y$, the inner sum is a Riemann sum of the continuous compactly supported function
$h(\cdot,y)$; since $K\Subset\Omega$, only complete cells contribute for small $\eps$, and the sum converges to
$\int_\Omega h(x,y)\,\dd x$.  Dominated convergence in $y$ gives the surface statement; the volume statement
is identical with $\sum_{k\in\mathbb Z^N}$ and $\int_Y$.
\end{proof}

On $\Gamma_\eps$ set $\dd\mu_\eps:=\eps\,\dd\sigma_x$, and equip $\Omega\times\partial T$ with the product
measure $\dd x\,\dd\sigma_y$.  An equivalence class in
$L^{\widetilde\Phi}(\Omega;C(\partial T))$ has no canonical restriction to the oscillating surface, so the
surface tests use continuous macroscopic representatives, as in the classical periodic-surface theory.

\begin{definition}\label{def:surface-two-scale}
Let $z_\eps\in L^\Phi(\Gamma_\eps,\mu_\eps)$, $\eps\in\mathscr E$.  We say that $z_\eps$ converges weakly
two-scale on the surface to $z_0\in L^\Phi(\Omega\times\partial T)$, and write
$z_\eps\overset{2s,\Gamma}{\rightharpoonup}z_0$, if
\begin{equation}\label{surfaceTwoScaleDefinition}
 \lim_{\eps\to0}\eps\int_{\Gamma_\eps}z_\eps(x)\psi\!\left(x,\frac{x}{\eps}\right)\dd\sigma_x
 =\int_\Omega\int_{\partial T}z_0(x,y)\psi(x,y)\,\dd\sigma_y\,\dd x
 \qquad\forall\,\psi\in C_c(\Omega;C(\partial T)).
\end{equation}
\end{definition}

The surface notion enjoys the analogous compactness under the scaled modular bound.

\begin{theorem}\label{thm:surface-two-scale-compactness}
If
\begin{equation}
\label{surfaceTwoScaleBound}
\sup_{\eps\in\mathscr E}\eps\int_{\Gamma_\eps}\Phi(|z_\eps|)\,\dd\sigma_x<\infty,
\end{equation}
then, after extraction, $z_\eps\overset{2s,\Gamma}{\rightharpoonup}z_0$ for some
$z_0\in L^\Phi(\Omega\times\partial T)$.  The same statement holds with $\Phi$ and $\widetilde\Phi$
interchanged.
\end{theorem}

\begin{proof}
Bound \eqref{surfaceTwoScaleBound} and the modular--norm equivalence give
$\sup_{\eps\in\mathscr E}\|z_\eps\|_{\Phi,\Gamma_\eps,\mu_\eps}<\infty$.
For every $\psi\in C_c(\Omega;C(\partial T))$ and every $a>0$, Lemma~\ref{lem:riemann} yields
\[
 \int_{\Gamma_\eps}\widetilde\Phi\!\left(\frac{|\psi(x,x/\eps)|}{a}\right)\dd\mu_\eps
 \longrightarrow
 \int_\Omega\int_{\partial T}\widetilde\Phi\!\left(\frac{|\psi(x,y)|}{a}\right)\dd\sigma_y\,\dd x,
\]
so the sampled tests are uniformly bounded in the complementary surface Orlicz norm, and
\eqref{OrliczHolder} bounds the functionals
$L_\eps(\psi):=\int_{\Gamma_\eps}z_\eps(x)\psi(x,x/\eps)\,\dd\mu_\eps$ uniformly.
Choose a countable subset of $C_c(\Omega;C(\partial T))$ dense in
$L^{\widetilde\Phi}(\Omega\times\partial T)$ and extract a diagonal subsequence on which $L_\eps$ converges on
it.  The uniform bound extends the limit to a bounded linear functional on
$L^{\widetilde\Phi}(\Omega\times\partial T)$, which reflexive Orlicz duality represents by some
$z_0\in L^\Phi(\Omega\times\partial T)$; density gives \eqref{surfaceTwoScaleDefinition}.  This is the
reflexive-Orlicz counterpart of \cite{AllaireDamlamianHornung1996,NeussRadu1996}.  Since
$\Phi$ and $\widetilde\Phi$ play symmetric roles in this argument, the interchanged statement follows in the
same way.
\end{proof}

\subsection{Reference-cell and perforated-domain estimates}
\label{subsec:cell-estimates}

We collect here the fixed-cell and perforated-domain estimates used below, including extension, trace and
restriction estimates.

\begin{lemma}\label{lem:cellextension}
There exist a linear operator $\mathsf E:W^1L^\Phi(Y^*)\to W^1L^\Phi(Y)$, with $\mathsf Ew=w$ in $Y^*$ and $\mathsf E1=1$, and a
constant $C_{\rm cell}>0$ such that, for every $w\in W^1L^\Phi(Y^*)$,
\begin{align}
 \int_Y\Phi(|\mathsf Ew|)\,\dd y&\le C_{\rm cell}\int_{Y^*}\Phi(C_{\rm cell}|w|)\,\dd y,\label{cellExtensionValue}\\
 \int_Y\Phi(|\nabla \mathsf Ew|)\,\dd y&\le C_{\rm cell}\int_{Y^*}\Phi(C_{\rm cell}|\nabla w|)\,\dd y,\label{cellExtensionGradient}\\
 \int_{\partial T}\Phi(|\gamma_0w|)\,\dd\sigma_y
 &\le C_{\rm cell}\int_{Y^*}\big[\Phi(C_{\rm cell}|w|)+\Phi(C_{\rm cell}|\nabla w|)\big]\dd y.\label{fixedtrace}
\end{align}
\end{lemma}

\begin{proof}
The construction is classical \cite{AcerbiChiadoPiatDalMasoPercivale1992}.  By
density of smooth functions under $\Delta_2$, it suffices to consider
$g\in C^\infty(\overline{Y^*})$.  Choose finitely many open cylinders $U_1,\dots,U_J$ with closures in $Y$,
covering $\partial T$, in which, after a rigid change of variables, $Y^*\cap U_j=\{(y',y_N):\ y_N>g_j(y')\}\cap U_j$
with $g_j$ Lipschitz.  The reflections $R_j(y',y_N)=(y',2g_j(y')-y_N)$ are bi-Lipschitz and equal the identity on
$\partial T$.  Choose $\rho_j\in C_c^\infty(U_j)$, $0\le\rho_j\le1$, with $\sum_{j=1}^J\rho_j=1$ on a neighbourhood
of $\partial T$; after shrinking the supports, $R_j$ maps $T\cap\supp\rho_j$ into $Y^*\cap U_j$.  Set
\[
 \mathsf E_0g:=g\ \text{ in }Y^*,
 \qquad
 \mathsf E_0g:=\sum_{j=1}^J\rho_j\,(g\circ R_j)\ \text{ in }T.
\]
Since $R_j$ is the identity on $\partial T$ and $\sum_{j=1}^{J}\rho_j=1$ there, both definitions have the same trace on
$\partial T$, so $\mathsf E_0g\in W^1L^\Phi(Y)$.  The changes of variables $z=R_j(y)$, the bounds on the Jacobians
of $R_j$, the chain rule
$\nabla(\rho_j\,g\circ R_j)=(g\circ R_j)\nabla\rho_j+\rho_j(\nabla g\circ R_j)DR_j$, convexity and
\eqref{scalingPhi} give
\begin{equation}\label{reflectionbounds}
 \int_T\Phi(|\mathsf E_0g|)\,\dd y\le C\int_{Y^*}\Phi(C|g|)\,\dd y,
 \qquad
 \int_T\Phi(|\nabla\mathsf E_0g|)\,\dd y\le C\int_{Y^*}\big[\Phi(C|g|)+\Phi(C|\nabla g|)\big]\dd y.
\end{equation}
Let $c(w):=|Y^*|^{-1}\int_{Y^*}w\,\dd y$ and define
\[
 \mathsf Ew:=c(w)+\mathsf E_0\big(w-c(w)\big).
\]
The operator $\mathsf E$ is linear, $\mathsf Ew=w$ in $Y^*$, and $\mathsf E1=1$.  Since $\nabla\mathsf Ew=\nabla\mathsf E_0(w-c(w))$,
the second estimate in \eqref{reflectionbounds}, applied to $g=w-c(w)$, and the modular Poincar\'e
inequality \eqref{modularPoincare} on the connected Lipschitz set $Y^*$ give \eqref{cellExtensionGradient}.  On
$T$, $|\mathsf Ew|\le|c(w)|+|\mathsf E_0(w-c(w))|$; Jensen's inequality gives
$\Phi(2|c(w)|)\le|Y^*|^{-1}\int_{Y^*}\Phi(2|w|)\,\dd y$, and the first estimate in \eqref{reflectionbounds} with
$|w-c(w)|\le|w|+|c(w)|$ gives \eqref{cellExtensionValue}.

For \eqref{fixedtrace}, choose a smooth vector field $X$ on $\overline{Y^*}$, vanishing near $\partial Y$, with
$X\cdot n\ge1$ on $\partial T$.  The divergence theorem applied to $\Phi(|w|)X$, the bound $|\nabla|w||\le|\nabla w|$
and \eqref{youngPhi} give
\[
 \int_{\partial T}\Phi(|w|)\,\dd\sigma_y
 \le\int_{Y^*}\Big[\Phi(|w|)\,|\diver X|+\varphi(|w|)\,|\nabla w|\,|X|\Big]\dd y
 \le C\int_{Y^*}\big[\Phi(|w|)+\Phi(|\nabla w|)\big]\dd y
\]
for smooth $w$, and \eqref{fixedtrace} follows by density.
\end{proof}
Set
\[
 X_\eps:=\{v\in W^1L^\Phi(\Omega_\eps):\gamma_0v=0\text{ on }\partial\Omega\},
 \qquad
 R_\eps:=\Omega\setminus\bigcup_{k\in K_\eps}\eps(k+Y).
\]
Inside surface integrals over $\Gamma_\eps$ and surface two-scale pairings, a Sobolev function is identified
with its trace on $\Gamma_\eps$ and $\gamma_0$ is omitted; $\gamma_0$ is kept in the definitions of function
spaces and in boundary conditions, and $\gamma_0^y$ is used only for traces of two-scale limits.  The estimates below are stated for
scalar functions; vector-valued versions follow componentwise.

\begin{proposition}\label{prop:uniform-extension-prelim}
There exist linear maps $P_\eps:X_\eps\to W^1_0L^\Phi(\Omega)$ with $P_\eps v=v$ in $\Omega_\eps$ and a
constant $C>0$, independent of $\eps$, such that
\begin{equation}\label{prelimUniformExtension}
 \int_\Omega\big[\Phi(|P_\eps v|)+\Phi(|\nabla P_\eps v|)\big]\dd x
 \le C\int_{\Omega_\eps}\big[\Phi(C|v|)+\Phi(C|\nabla v|)\big]\dd x
 \qquad\forall\,v\in X_\eps.
\end{equation}
Moreover, for every $a>0$, every $k\in K_\eps$ and every $v\in X_\eps$,
\begin{equation}
\label{holeextensionlocal}
\int_{\eps(k+Y)}
\Phi\big(a\eps\abs{\nabla P_\eps v}\big)\,\dd x
\le
C_{\rm cell}\int_{\eps(k+Y^*)}
\Phi\big(C_{\rm cell}a\eps\abs{\nabla v}\big)\,\dd x.
\end{equation}
\end{proposition}

\begin{proof}
For $k\in K_\eps$ put $V_k(y)=v(\eps(k+y))$, $y\in Y^*$, and define $P_\eps v$ on $\eps(k+Y)$ by
$P_\eps v(\eps(k+y))=(\mathsf EV_k)(y)$; on $R_\eps$ we leave $v$ unchanged.  Since $\mathsf EV_k=V_k$ on $Y^*$, in
particular near the outer faces of the cell, the cellwise extensions patch across neighbouring faces, and
the zero trace of $v$ on $\partial\Omega$ gives $P_\eps v\in W^1_0L^\Phi(\Omega)$.

Because $\mathsf E$ is linear, \eqref{cellExtensionGradient} may be applied to $aV_k$ for any $a>0$.  Since
$\nabla_y(\mathsf EV_k)(y)=\eps\nabla P_\eps v(\eps(k+y))$ and $\nabla_yV_k(y)=\eps\nabla v(\eps(k+y))$, the change of
variables $x=\eps(k+y)$ gives
\[
\begin{aligned}
 \int_{\eps(k+Y)}\Phi\big(a\eps|\nabla P_\eps v|\big)\dd x
 &=\eps^N\int_Y\Phi\big(a|\nabla_y\mathsf EV_k|\big)\dd y\\
 &\le C_{\rm cell}\eps^N\int_{Y^*}\Phi\big(C_{\rm cell}a|\nabla_yV_k|\big)\dd y
 =C_{\rm cell}\int_{\eps(k+Y^*)}\Phi\big(C_{\rm cell}a\eps|\nabla v|\big)\dd x,
\end{aligned}
\]
which is \eqref{holeextensionlocal}.  Taking $a=\eps^{-1}$ gives the unscaled cellwise gradient estimate,
and \eqref{cellExtensionValue} gives
$\int_{\eps(k+Y)}\Phi(|P_\eps v|)\le C_{\rm cell}\int_{\eps(k+Y^*)}\Phi(C_{\rm cell}|v|)$.  Summing over $k\in K_\eps$ and adding
the unchanged contribution of $R_\eps$ yields \eqref{prelimUniformExtension}.
\end{proof}

Rescaling the fixed-cell trace estimate \eqref{fixedtrace} gives the corresponding bound on the pore walls.

\begin{proposition}\label{prop:scaled-trace-prelim}
For every $v\in W^1L^\Phi(\Omega_\eps)$,
\begin{equation}\label{prelimScaledTrace}
 \eps\int_{\Gamma_\eps}\Phi(|v|)\,\dd\sigma_x
 \le C_{\rm cell}\int_{\Omega_\eps}\big[\Phi(C_{\rm cell}|v|)+\Phi(C_{\rm cell}\eps|\nabla v|)\big]\dd x.
\end{equation}
\end{proposition}

\begin{proof}
For $k\in K_\eps$ put $V_k(y)=v(\eps(k+y))$.  Applying \eqref{fixedtrace} to $V_k$ and multiplying by
$\eps^N$ gives, since $\dd\sigma_x=\eps^{N-1}\dd\sigma_y$,
\[
 \eps\int_{\eps(k+\partial T)}\Phi(|v|)\,\dd\sigma_x
 \le C_{\rm cell}\int_{\eps(k+Y^*)}
 \big[\Phi(C_{\rm cell}|v|)+\Phi(C_{\rm cell}\eps|\nabla v|)\big]\dd x.
\]
Since $\Gamma_\eps=\bigcup_{k\in K_\eps}\eps(k+\partial T)$, summation over $k\in K_\eps$ yields
\eqref{prelimScaledTrace}.
\end{proof}

These tools are now applied to the microscopic problem, whose $\eps$-uniform estimates are the subject of
the next section.

The pressure analysis of Section~\ref{sec:effectivedarcy} and the solenoidal corrections of
Subsection~\ref{sec:preliminarylimit} use a restriction operator of Tartar type \cite{Allaire1989,Tartar1980}
from $\Omega$ to the perforated domain, whose Orlicz bounds follow from the modular Bogovskii estimate of
\cite[Theorem~6.6]{DieningRuzickaSchumacher2010}.

\begin{lemma}
\label{lem:tartarrestriction}
Fix a smooth open set $U$ with $\overline T\subset U\Subset Y$, set $A:=U\setminus\overline T$, and choose
$\eta\in C^\infty(\overline U)$ with $\eta=0$ in a neighbourhood of $\overline T$ and $\eta=1$ in a neighbourhood of
$\partial U$.  For $k\in K_\eps$ write $A_{\eps,k}:=\eps(k+A)$, $T_{\eps,k}:=\eps(k+T)$ and
$\eta_{\eps,k}(x):=\eta(x/\eps-k)$.  For $0<\eps<1$ there is a linear operator $\mathcal R_\eps$, independent of the $N$-function, which maps
$W^1_0L^\Theta(\Omega)^N$ into $W^1_0L^\Theta(\Omega_\eps)^N$ for every $N$-function $\Theta$ with
$\Theta,\widetilde\Theta\in\Delta_2$, and such that, for every such $w$:
\begin{enumerate}[label=\textup{(\roman*)}]
\item $\mathcal R_\eps w=w$ in $\Omega_\eps\setminus\bigcup_{k\in K_\eps}\eps(k+U)$;
\item $\diver\mathcal R_\eps w=\diver w+\sum_{k\in K_\eps}m_{\eps,k}(w)\,\mathbf 1_{A_{\eps,k}}$ in $\Omega_\eps$, where
$m_{\eps,k}(w):=|A_{\eps,k}|^{-1}\int_{T_{\eps,k}}\diver w\,\dd x$;
\item there is $C>0$, depending only on $Y^*$, $T$, $U$, $N$ and the $\Delta_2$ data of $\Theta$ and
$\widetilde\Theta$, such that
\begin{equation}
\label{tartarmodular}
\int_{\Omega_\eps}\Big[\Theta(|\mathcal R_\eps w|)+\Theta(\eps|\nabla\mathcal R_\eps w|)\Big]\dd x
\le C\int_\Omega\Big[\Theta(C|w|)+\Theta(C\eps|\nabla w|)\Big]\dd x.
\end{equation}
\end{enumerate}
In particular, $\|\mathcal R_\eps w\|_{\Theta,\Omega_\eps}+\eps\|\nabla\mathcal R_\eps w\|_{\Theta,\Omega_\eps}
\le C\big(\|w\|_{\Theta,\Omega}+\eps\|\nabla w\|_{\Theta,\Omega}\big)$, and this holds for $\Theta=\Phi$ and for
$\Theta(t)=t^q$, $1<q<\infty$.
\end{lemma}

\begin{proof}
The annulus $A$ is a bounded Lipschitz domain, hence a John domain, and satisfies the emanating chain
condition \cite[Theorem~3.8]{DieningRuzickaSchumacher2010}.  Let $\mathfrak B_A$ be the Bogovskii operator of
\cite[Theorem~6.6]{DieningRuzickaSchumacher2010} on $A$; it is one linear operator for all admissible
$N$-functions, satisfies $\diver\mathfrak B_Ag=g$ and the modular bound
$\int_A\Theta(|\nabla\mathfrak B_Ag|)\,\dd y\le C\int_A\Theta(|g|)\,\dd y$ for $g\in L^\Theta_0(A)$.

Fix $w$ and $k\in K_\eps$, and set $h_{\eps,k}:=\diver\big((\eta_{\eps,k}-1)w\big)$ in $A_{\eps,k}$.  On the
outer boundary of $A_{\eps,k}$ the factor $\eta_{\eps,k}-1$ vanishes; on $\partial T_{\eps,k}$ it equals $-1$ and
the outward normal of the annulus is the inner normal of the hole.  The divergence theorem therefore gives
$\int_{A_{\eps,k}}h_{\eps,k}\,\dd x=\int_{\partial T_{\eps,k}}w\cdot\nu\,\dd\sigma_x=\int_{T_{\eps,k}}\diver w\,\dd x$,
where $\nu$ is the outward normal of the hole, so that $g_{\eps,k}:=h_{\eps,k}-m_{\eps,k}(w)$ has zero mean on
$A_{\eps,k}$.  With $\widehat g(y):=g_{\eps,k}(\eps(k+y))$, $y\in A$, set
\[
 z_{\eps,k}(x):=\eps\,(\mathfrak B_A\widehat g)\Big(\frac x\eps-k\Big),\qquad x\in A_{\eps,k}.
\]
Then $z_{\eps,k}\in W^1_0L^\Theta(A_{\eps,k})^N$ and $\diver z_{\eps,k}=g_{\eps,k}$.  Define $\mathcal R_\eps w:=
\eta_{\eps,k}w-z_{\eps,k}$ in $A_{\eps,k}$, $k\in K_\eps$, and $\mathcal R_\eps w:=w$ in the rest of $\Omega_\eps$.
Since $\eta_{\eps,k}=1$ near the outer boundary of $A_{\eps,k}$ and $z_{\eps,k}$ has zero trace, the pieces
match, and since $\eta_{\eps,k}=0$ near $\partial T_{\eps,k}$, $\mathcal R_\eps w$ vanishes on $\Gamma_\eps$; it
vanishes on $\partial\Omega$ because $w$ does and only complete cells are modified.  This gives (i), and in
$A_{\eps,k}$,
\[
 \diver\mathcal R_\eps w=\diver(\eta_{\eps,k}w)-h_{\eps,k}+m_{\eps,k}(w)=\diver w+m_{\eps,k}(w),
\]
which is (ii).  All operations are linear and independent of $\Theta$.

For (iii), the change of variables $x=\eps(k+y)$ and the modular Bogovskii bound, applied to $\eps\widehat g$
by linearity, give
\[
 \int_{A_{\eps,k}}\Theta(\eps|\nabla z_{\eps,k}|)\,\dd x
 =\eps^N\int_A\Theta(\eps|\nabla\mathfrak B_A\widehat g|)\,\dd y
 \le C\int_{A_{\eps,k}}\Theta(\eps|g_{\eps,k}|)\,\dd x.
\]
Since $h_{\eps,k}=(\eta_{\eps,k}-1)\diver w+\nabla\eta_{\eps,k}\cdot w$ and $|\nabla\eta_{\eps,k}|\le C/\eps$, we have
$\eps|h_{\eps,k}|\le C(\eps|\nabla w|+|w|)$.  Jensen's inequality gives
$|A_{\eps,k}|\,\Theta(\eps|m_{\eps,k}(w)|)\le C\int_{T_{\eps,k}}\Theta(C\eps|\nabla w|)\,\dd x$, because
$|T_{\eps,k}|/|A_{\eps,k}|=|T|/|A|$.  Convexity and $\Delta_2$ therefore yield
\[
 \int_{A_{\eps,k}}\Theta(\eps|\nabla z_{\eps,k}|)\,\dd x
 \le C\int_{\eps(k+U)}\Big[\Theta(C|w|)+\Theta(C\eps|\nabla w|)\Big]\dd x.
\]
The scaled form of \eqref{modularPoincarezero} on $A_{\eps,k}$, whose diameter is at most $\eps\sqrt N$, bounds
$\int_{A_{\eps,k}}\Theta(|z_{\eps,k}|)$ by the same quantity, and $\eps|\nabla(\eta_{\eps,k}w)|\le\eps|\nabla w|+C|w|$.
Summing over the pairwise disjoint sets $\eps(k+U)$, $k\in K_\eps$, and adding the contribution of the
unmodified region proves \eqref{tartarmodular}.  Applying \eqref{tartarmodular} to $w/\lambda$ with
$\lambda=\|w\|_{\Theta,\Omega}+\eps\|\nabla w\|_{\Theta,\Omega}$ and using the $\Delta_2$ condition for $\Theta$ gives the
norm bound.
\end{proof}

\section{Microscopic problem, coercivity, and uniform estimates}
\label{sec:microscopic-estimates}

Throughout this section, Assumptions~\ref{ass:geometry} and~\ref{ass:orlicz} hold.  We convert bulk and
surface coercivity into $\eps$-uniform velocity control: first on one cell, then on the complete-cell region
and the boundary remainder, and finally for the microscopic problem.

\subsection{Cell Orlicz--Korn coercivity with Navier friction}
\label{subsec:cellkorn}

The cell inequality controls the velocity modulo rigid motions, which impermeability and friction then
remove.  Denote the finite-dimensional space of infinitesimal rigid motions by
\begin{equation}
\label{rigidspace}
\mathcal R:=\left\{r:\R^N\to\R^N:\ r(y)=a+Wy,
\quad a\in\R^N,\quad W\in\R^{N\times N},\quad W^T=-W\right\}.
\end{equation}
Every $r\in\mathcal R$ satisfies $\D r=0$.  Conversely, if $Q\subset\R^N$ is connected
and $v\in W^{1,1}(Q)^N$ satisfies $\D v=0$ in distributions, then $v$ agrees almost
everywhere with an element of $\mathcal R$.

\begin{lemma}
\label{lem:kornmod}
Let $Q\subset\R^N$ be bounded, connected and Lipschitz.  There
exist constants $C_{\rm K},c_{\rm K}>0$, depending only on $Q$, $N$ and the $\Delta_2$ data of
$\Phi,\widetilde\Phi$, such that for every $v\in W^1L^{\Phi}(Q)^N$ one can choose
$r_v(y)=a_v+W_vy\in\mathcal R$ satisfying
\begin{equation}
\label{kornmod}
\int_Q\!\left[
\Phi(\abs{v-r_v})+\Phi(\abs{\nabla v-W_v})
\right]\dd y
\le
C_{\rm K}\int_Q\Phi(c_{\rm K}\abs{\D v})\,\dd y.
\end{equation}
\end{lemma}

\begin{proof}
Since $\Phi,\widetilde\Phi\in\Delta_2$, the Young function $\Phi$ satisfies the $\Delta_2$ and $\nabla_2$ conditions
globally.  Hence \cite[Corollaries~3.4 and~3.6]{Cianchi2014}, applied with both Young functions equal to $\Phi$ on the Lipschitz domain
$Q$, give inequality (3.8) of that paper with vanishing additive constant: denoting by $\langle\cdot\rangle_Q$
the mean over $Q$,
\[
\int_Q\Phi\big(\abs{\nabla v-\langle\nabla v\rangle_Q}\big)\,\dd y
\le \int_Q\Phi\big(C\abs{\D v-\langle\D v\rangle_Q}\big)\,\dd y.
\]
The matrix
\[
W_v:=\langle\nabla v\rangle_Q-\langle\D v\rangle_Q=\frac1{|Q|}\int_Q\frac12\big(\nabla v-(\nabla v)^T\big)\,\dd y
\]
is skew-symmetric, and $\abs{\nabla v-W_v}\le\abs{\nabla v-\langle\nabla v\rangle_Q}+\abs{\langle\D v\rangle_Q}$.  Convexity,
Jensen's inequality for the means and \eqref{scalingPhi} therefore give
\[
\int_Q\Phi(\abs{\nabla v-W_v})\,\dd y
\le C\int_Q\Phi(c\abs{\D v})\,\dd y.
\]
Choose the translation $a_v$ so that
\[
\int_Q\big(v(y)-a_v-W_vy\big)\,\dd y=0.
\]
The modular Poincar\'e inequality \eqref{modularPoincare} on $Q$ then yields
\[
\int_Q\Phi(\abs{v-r_v})\,\dd y
\le C\int_Q\Phi(c\abs{\nabla v-W_v})\,\dd y.
\]
Combining the last two estimates with \eqref{scalingPhi} proves \eqref{kornmod}.
\end{proof}

A rigid motion is controlled by its trace on the pore wall.

\begin{lemma}
\label{lem:rigidboundary}
There exist constants $C_R,c_R>0$, depending only on $Y^*$, $T$ and the $\Delta_2$ data of $\Phi$, such that every
$r(y)=a+Wy\in\mathcal R$ satisfies
\begin{equation}
\label{rigidboundary}
\int_{Y^*}\left[\Phi(\abs{r})+\Phi(\abs{W})\right]\dd y
\le
C_R\int_{\partial T}\Phi(c_R\abs{r})\,\dd\sigma_y.
\end{equation}
In particular, if $r=0$ for $\sigma_y$-almost every point of $\partial T$, then
$r\equiv0$ in $\R^N$.
\end{lemma}

\begin{proof}
We first prove injectivity of the restriction map $r\mapsto r|_{\partial T}$ on $\mathcal R$.
Suppose $r(y)=a+Wy$ vanishes for $\sigma_y$-almost every $y\in\partial T$.
By continuity of $r$ and of the $C^{1,1}$ hypersurface, $r=0$ on all of $\partial T$.
Fix $y_0\in\partial T$.  Then
\[
W(y-y_0)=0\qquad\text{for every }y\in\partial T.
\]
If $W\ne0$, a skew-symmetric matrix has rank at least two, hence $\dim\ker W\le N-2$.
Consequently $\partial T\subset y_0+\ker W$, an affine subspace of dimension at most
$N-2$.  This is impossible because the $C^{1,1}$ hypersurface $\partial T$ has positive
surface measure $\sigma_y(\partial T)$.  Therefore $W=0$, and then $a=0$.

Since $\mathcal R$ is finite-dimensional and the restriction map is injective,
\[
\abs{a}+\abs{W}
\le C\int_{\partial T}\abs{a+Wy}\,\dd\sigma_y.
\]
Indeed the right-hand side is a norm on $\mathcal R$, and all norms on a finite-dimensional
space are equivalent.  Jensen's inequality for the probability measure $|\partial T|^{-1}\dd\sigma_y$ and
\eqref{scalingPhi} therefore give
\begin{equation}
\label{rigid-jensen}
\Phi\big(c(\abs{a}+\abs{W})\big)
\le
C\int_{\partial T}\Phi(c_R\abs{r})\,\dd\sigma_y.
\end{equation}
On the bounded set $Y^*$,
\[
\abs{r(y)}+\abs{W}\le C(\abs{a}+\abs{W}),
\]
and hence \eqref{rigidboundary} follows from \eqref{rigid-jensen} and \eqref{scalingPhi}.
\end{proof}

Define the local impermeable space
\begin{equation}
\label{Wcell}
\mathcal W^{\rm imp}_\Phi(Y^*):=
\left\{v\in W^1L^{\Phi}(Y^*)^N:\
\gamma_0 v\cdot n=0\quad\text{on }\partial T\right\}.
\end{equation}
No condition is imposed on the outer faces $\partial Y$.

\begin{proposition}
\label{prop:cellkorn}
There exist constants $C,c>0$, depending only on
$N$, $Y^*$, $T$ and the Orlicz data, such that every
$v\in\mathcal W^{\rm imp}_\Phi(Y^*)$ satisfies
\begin{equation}
\label{cellkorn}
\int_{Y^*}\left[
\Phi(\abs{v})+\Phi(\abs{\nabla v})
\right]\dd y
\le
C\left\{
\int_{Y^*}\Phi(c\abs{\D v})\,\dd y
+
\int_{\partial T}\Phi(c\abs{(\gamma_0 v)_\tau})\,\dd\sigma_y
\right\}.
\end{equation}
\end{proposition}

\begin{proof}
Fix $v\in\mathcal W^{\rm imp}_\Phi(Y^*)$ and choose, by Lemma~\ref{lem:kornmod}, a rigid motion
\[
r(y)=a+Wy\in\mathcal R
\]
such that, with $w:=v-r$,
\begin{equation}
\label{wbound}
\int_{Y^*}\left[
\Phi(\abs{w})+\Phi(\abs{\nabla w})
\right]\dd y
\le
C\int_{Y^*}\Phi(c\abs{\D v})\,\dd y.
\end{equation}
Here we used $\nabla w=\nabla v-W$.

We next estimate the rigid part from its trace on $\partial T$.  Since
$\gamma_0 v\cdot n=0$,
\[
r\cdot n=(v-w)\cdot n=-w\cdot n
\qquad\text{on }\partial T.
\]
For the tangential component, $r_\tau=v_\tau-w_\tau$.
Consequently
\begin{equation}
\label{rtracepoint}
\abs{r}
\le \abs{r_\tau}+\abs{r\cdot n}
\le \abs{v_\tau}+2\abs{w}
\qquad\text{on }\partial T.
\end{equation}
By convexity and $\Delta_2$,
\begin{equation}
\label{rtracephi}
\int_{\partial T}\Phi(c_R\abs{r})\,\dd\sigma_y
\le C\left[
\int_{\partial T}\Phi(c\abs{v_\tau})\,\dd\sigma_y
+
\int_{\partial T}\Phi(c\abs{\gamma_0 w})\,\dd\sigma_y
\right].
\end{equation}
The fixed-cell trace estimate \eqref{fixedtrace}, combined with \eqref{scalingPhi} and \eqref{wbound},
gives
\begin{equation}
\label{wtracebound}
\int_{\partial T}\Phi(c\abs{\gamma_0 w})\,\dd\sigma_y
\le
C\int_{Y^*}\Phi(c\abs{\D v})\,\dd y.
\end{equation}
Combining Lemma~\ref{lem:rigidboundary}, \eqref{rtracephi} and
\eqref{wtracebound}, we obtain
\begin{equation}
\label{rigidfinal}
\int_{Y^*}\left[
\Phi(\abs{r})+\Phi(\abs{W})
\right]\dd y
\le C\left[
\int_{Y^*}\Phi(c\abs{\D v})\,\dd y
+
\int_{\partial T}\Phi(c\abs{v_\tau})\,\dd\sigma_y
\right].
\end{equation}
Finally $v=w+r$ and $\nabla v=\nabla w+W$.  A last application of convexity and the
$\Delta_2$ property, together with \eqref{wbound} and \eqref{rigidfinal}, yields
\eqref{cellkorn}.
\end{proof}

\begin{remark}
If the right-hand side of \eqref{cellkorn} vanishes, then $\D v=0$, so $v$ is a rigid
motion.  The friction term gives $v_\tau=0$ on $\partial T$, whereas membership in
$\mathcal W^{\rm imp}_\Phi(Y^*)$ gives $v\cdot n=0$.  Hence $v=0$ on $\partial T$, and
Lemma~\ref{lem:rigidboundary} implies $v\equiv0$.  Thus the boundary friction and the
impermeability condition jointly remove every translation and infinitesimal rotation.  No periodicity
is used, which is what the rescaling below requires: the traces of a global velocity on opposite faces of
a single cell $\eps(k+Y^*)$ need not coincide, and periodicity enters only in the limiting cell problem.
\end{remark}

In the limiting cell problem, periodicity alone removes the rigid motions, without friction and for every
hole allowed by Assumption~\ref{ass:geometry}.

\begin{lemma}
\label{lem:periodicKorn}
There is $C_P>0$ such that
\begin{equation}
\label{periodicKorn}
 \|v\|_{W^1L^\Phi(Y^*)}\le C_P\|\D_yv\|_{\Phi,Y^*}
 \qquad\forall\,v\in W^1_{\rm per}L^\Phi(Y^*)^N\text{ with }\gamma_0v\cdot n=0\text{ on }\partial T.
\end{equation}
\end{lemma}

\begin{proof}
Applying \eqref{kornmod} on $Q=Y^*$ to $v/\lambda$, $\lambda=\|\D_yv\|_{\Phi,Y^*}$, and using the linearity of
$v\mapsto r_v$ gives $\|v-r_v\|_{W^1L^\Phi(Y^*)}\le C\|\D_yv\|_{\Phi,Y^*}$.  If \eqref{periodicKorn} failed, there would be
admissible $v_n$ with $\|v_n\|_{W^1L^\Phi(Y^*)}=1$ and $\|\D_yv_n\|_{\Phi,Y^*}\to0$; then $v_n-r_{v_n}\to0$, the rigid
motions $r_{v_n}$ are bounded in the finite-dimensional space $\mathcal R$, and a subsequence converges in
$W^1L^\Phi(Y^*)$ to some $r(y)=a+Wy$ with $\|r\|_{W^1L^\Phi(Y^*)}=1$, which is admissible because the admissible
fields form a closed subspace.  Since $\overline T\subset Y$, the faces of $Y$ lie in $\overline{Y^*}$, and periodicity
of the traces gives $r(y+e_i)=r(y)$ for $y$ on the face $\{y_i=0\}$, that is, $We_i=0$ for every $i$; hence $W=0$.
Then $a\cdot n=0$ on $\partial T$, and
$|a|^2|T|=\int_T\diver\big((a\cdot y)a\big)\,\dd y=\int_{\partial T}(a\cdot y)(a\cdot n)\,\dd\sigma_y=0$ gives $a=0$, a
contradiction.
\end{proof}

\subsection{Rescaling to a microscopic pore cell}
\label{subsec:rescaling}

Fix $\eps\in(0,1)$ and denote by
\[
\Omega_\eps^{\rm cell,*}:=\bigcup_{k\in K_\eps}\eps(k+Y^*)
\]
the union of the complete fluid cells.  For $k\in K_\eps$ set
\[
Y^*_{\eps,k}:=\eps(k+Y^*),\qquad
\Gamma_{\eps,k}:=\eps(k+\partial T).
\]
Let $v\in W^1L^{\Phi}(Y^*_{\eps,k})^N$ satisfy
\[
\gamma_0 v\cdot n_\eps=0\qquad\text{on }\Gamma_{\eps,k}.
\]
Define
\[
V(y):=v(\eps(k+y)),\qquad y\in Y^*.
\]
Then
\begin{equation}
\label{scalegrad}
\nabla_yV(y)=\eps\nabla_xv(\eps(k+y)),\qquad
\D_yV(y)=\eps\D_xv(\eps(k+y)).
\end{equation}
The tangential projection is preserved by the dilation because the normal of
$\Gamma_{\eps,k}$ at $x=\eps(k+y)$ is $n(y)$.

\begin{proposition}
\label{prop:rescaled}
There exist constants $C,c>0$, independent of $\eps$ and $k$, such that
{\small
\begin{equation}
\label{rescaledcell}
\int_{Y^*_{\eps,k}}[\Phi(\abs v)+\Phi(\eps\abs{\nabla v})]\dd x\le C\left\{\int_{Y^*_{\eps,k}}\Phi(c\eps\abs{\D v})\dd x+\eps\int_{\Gamma_{\eps,k}}\Phi(c\abs{v_\tau})\dd\sigma_x\right\}.
\end{equation}
}
\end{proposition}

\begin{proof}
Apply Proposition~\ref{prop:cellkorn} to $V$.  By \eqref{scalegrad},
\[
\abs{\nabla_yV}=\eps\abs{\nabla_xv},\qquad
\abs{\D_yV}=\eps\abs{\D_xv}.
\]
The changes of variables $x=\eps(k+y)$, $\dd x=\eps^N\dd y$ and, on the pore wall,
$\dd\sigma_x=\eps^{N-1}\dd\sigma_y$ give
\[
\eps^N\int_{\partial T}\Phi(c\abs{V_\tau})\,\dd\sigma_y
=
\eps\int_{\Gamma_{\eps,k}}\Phi(c\abs{v_\tau})\,\dd\sigma_x.
\]
Multiplying the fixed-cell inequality by $\eps^N$ therefore gives exactly
\eqref{rescaledcell}.  All constants come from the reference cell and are independent of
$\eps$ and $k$; the factor $\eps$ in front of the surface term is the quotient $\eps^N/\eps^{N-1}$ of the
volume and surface scaling factors, that is, the critical surface measure $\eps\,\dd\sigma_x$.
\end{proof}

\begin{corollary}
\label{cor:completecells}
Let $v\in W^1L^{\Phi}(\Omega_\eps)^N$ satisfy $v\cdot n_\eps=0$ on $\Gamma_\eps$.
Then
{\small
\begin{equation}
\label{completecells}
\int_{\Omega_\eps^{\rm cell,*}}[\Phi(\abs v)+\Phi(\eps\abs{\nabla v})]\dd x\le C\left\{\int_{\Omega_\eps^{\rm cell,*}}\Phi(c\eps\abs{\D v})\dd x+\eps\int_{\Gamma_\eps}\Phi(c\abs{v_\tau})\dd\sigma_x\right\}.
\end{equation}
}
where $C,c$ are independent of $\eps$.
\end{corollary}

\begin{proof}
Since $\Omega_\eps^{\rm cell,*}=\bigcup_{k\in K_\eps}Y^*_{\eps,k}$ and
$\Gamma_\eps=\bigcup_{k\in K_\eps}\Gamma_{\eps,k}$, with pairwise disjoint pieces up to null sets, the sum of
\eqref{rescaledcell} over $k\in K_\eps$ is exactly \eqref{completecells}.
\end{proof}

\subsection{Uniform control of the outer boundary layer}
\label{subsec:boundarylayer}

The complete-cell estimate \eqref{completecells} does not control the unperforated remainder $R_\eps$
defined in Subsection~\ref{subsec:cell-estimates}.  Instead of proving a Korn inequality on the possibly
jagged set $R_\eps$, we fill only the periodic holes and work on the fixed domain $\Omega$, so that no
constant depends on sliver-shaped intersections of boundary cells with $\Omega$.  The argument combines
\eqref{completecells}, a thin-layer Orlicz--Poincar\'e inequality at $\partial\Omega$ and the extension
of Proposition~\ref{prop:uniform-extension-prelim}.

For $\delta>0$ set
\begin{equation}
\label{Sdelta}
\Sigma_\delta:=\{x\in\Omega:\ \operatorname{dist}(x,\partial\Omega)<\delta\}.
\end{equation}

\begin{lemma}
\label{lem:Rstrip}
Up to a set of Lebesgue measure zero,
\begin{equation}
\label{Rstrip}
R_\eps\subset \Sigma_{\sqrt N\,\eps}.
\end{equation}
\end{lemma}

\begin{proof}
Take $x\in\Omega$ which does not belong to a grid hyperplane and assume
\[
\operatorname{dist}(x,\partial\Omega)>\sqrt N\,\eps.
\]
There is a unique $k\in\Z^N$ such that $x\in\eps(k+Y)$.  Every
$z\in\eps(k+\overline Y)$ satisfies
\[
\abs{z-x}\le \operatorname{diam}(\eps Y)=\sqrt N\,\eps
<\operatorname{dist}(x,\partial\Omega).
\]
Hence $z\in\Omega$ for every $z\in\eps(k+\overline Y)$, so
$\eps(k+\overline Y)\subset\Omega$ and therefore $k\in K_\eps$.  Thus
$x\notin R_\eps$.  The grid hyperplanes have zero Lebesgue measure, which proves
\eqref{Rstrip}.
\end{proof}

The next estimate is the only place in this boundary-layer argument where the $C^{1,1}$ regularity of the
macroscopic boundary is used; a bi-Lipschitz collar of $\partial\Omega$ would suffice.

\begin{lemma}
\label{lem:stripPoincare}
There exist $\delta_0>0$ and constants $C,c>0$, depending only on $\Omega$ and the
$\Delta_2$ data of $\Phi$, such that for every $0<\delta<\delta_0$ and every
$w\in W^1_0L^{\Phi}(\Omega)^N$,
\begin{equation}
\label{stripPoincare}
\int_{\Sigma_\delta}\Phi(\abs{w})\,\dd x
\le
C\int_{\Sigma_\delta}\Phi(c\delta\abs{\nabla w})\,\dd x.
\end{equation}
\end{lemma}

\begin{proof}
Because $\partial\Omega$ is $C^{1,1}$, there is $\delta_0>0$ such that the inner normal map
\[
\Theta:\partial\Omega\times(0,\delta_0)\longrightarrow \Sigma_{\delta_0},
\qquad
\Theta(z,s)=z-s\nu(z),
\]
is a bi-Lipschitz $C^{0,1}$ parametrization.
Its Jacobian $J_\Theta(z,s)$ satisfies
\begin{equation}
\label{Jacobianstrip}
0<c_\Theta\le J_\Theta(z,s)\le C_\Theta<\infty
\qquad (z\in\partial\Omega,\ 0<s<\delta_0).
\end{equation}
It is enough first to argue for smooth fields with zero boundary trace and then use density.
For almost every $z\in\partial\Omega$ and $0<s<\delta$,
\[
w(\Theta(z,s))
=-\int_0^s \nabla w(\Theta(z,r))\nu(z)\,\dd r,
\]
and hence
\[
\abs{w(\Theta(z,s))}
\le \int_0^\delta \abs{\nabla w(\Theta(z,r))}\,\dd r.
\]
Writing the last integral as $\delta$ times an average and applying Jensen's inequality,
\begin{equation}
\label{stripJensen}
\Phi\big(\abs{w(\Theta(z,s))}\big)
\le
\frac1\delta\int_0^\delta
\Phi\big(\delta\abs{\nabla w(\Theta(z,r))}\big)\,\dd r.
\end{equation}
Integrate \eqref{stripJensen} with respect to $s\in(0,\delta)$ and
$z\in\partial\Omega$.  The uniform Jacobian bounds \eqref{Jacobianstrip} give
\eqref{stripPoincare}, with \eqref{scalingPhi}.  Density of
$C^\infty(\overline\Omega)^N\cap W^1_0L^{\Phi}(\Omega)^N$ in the reflexive
Orlicz--Sobolev space completes the proof.
\end{proof}

We can now close the outer-layer estimate; its right-hand side is the global Darcy energy of the
microscopic problem.

\begin{theorem}
\label{thm:boundarylayer}
There are
$\eps_0\in(0,1]$ and constants $C,c>0$, independent of $0<\eps<\eps_0$, such that every
$v\in W^1L^{\Phi}(\Omega_\eps)^N$ satisfying
\begin{equation}
\label{boundaryadmissible}
\gamma_0 v=0\quad\text{on }\partial\Omega,
\qquad
\gamma_0 v\cdot n_\eps=0\quad\text{on }\Gamma_\eps,
\end{equation}
satisfies
{\small
\begin{equation}
\label{boundarylayerproved}
\int_{R_\eps}[\Phi(\abs v)+\Phi(\eps\abs{\nabla v})]\dd x\le C\left\{\int_{\Omega_\eps}\Phi(c\eps\abs{\D v})\dd x+\eps\int_{\Gamma_\eps}\Phi(c\abs{v_\tau})\dd\sigma_x\right\}.
\end{equation}
}
\end{theorem}

\begin{proof}
Let $\widehat v_\eps:=P_\eps v$.  By Proposition~\ref{prop:uniform-extension-prelim}, applied componentwise,
$\widehat v_\eps\in W^1_0L^{\Phi}(\Omega)^N$ and $\widehat v_\eps=v$ in $\Omega_\eps$, hence in particular
on $R_\eps$.

The zero-trace Orlicz--Korn inequality on the fixed Lipschitz domain $\Omega$, namely inequality (3.4) of
\cite[Theorem~3.1 and Corollary~3.6]{Cianchi2014} with both Young functions equal to $\Phi$ and vanishing
additive constant, applied to $\eps w$, gives constants $C,c>0$, independent of $\eps$ and $w$, with
\begin{equation}
\label{scaledKornOmega}
\int_\Omega\Phi(\eps\abs{\nabla w})\,\dd x
\le
C\int_\Omega\Phi(c\eps\abs{\D w})\,\dd x
\qquad\forall\,w\in W^1_0L^{\Phi}(\Omega)^N.
\end{equation}
We first control the scaled gradient on $R_\eps$.  By \eqref{scaledKornOmega},
\begin{equation}
\label{BLstep1a}
\int_{R_\eps}\Phi(\eps\abs{\nabla v})\,\dd x
\le
\int_\Omega\Phi(\eps\abs{\nabla \widehat v_\eps})\,\dd x
\le
C\int_\Omega\Phi(c\eps\abs{\D \widehat v_\eps})\,\dd x.
\end{equation}
Since $\widehat v_\eps=v$ on $\Omega_\eps$ and
$\abs{\D \widehat v_\eps}\le\abs{\nabla \widehat v_\eps}$ inside the holes,
{\small
\begin{equation}
\label{BLdecompose}
\int_\Omega\Phi(c\eps\abs{\D \widehat v_\eps})\dd x\le\int_{\Omega_\eps}\Phi(c\eps\abs{\D v})\dd x+\sum_{k\in K_\eps}\int_{\eps(k+T)}\Phi(c\eps\abs{\nabla \widehat v_\eps})\dd x.
\end{equation}
}
Estimate \eqref{holeextensionlocal}, summed over $k\in K_\eps$, gives
\begin{equation}
\label{BLholecost}
\sum_{k\in K_\eps}
\int_{\eps(k+T)}
\Phi(c\eps\abs{\nabla \widehat v_\eps})\,\dd x
\le
C\int_{\Omega_\eps^{\rm cell,*}}
\Phi(C\eps\abs{\nabla v})\,\dd x.
\end{equation}
By the $\Delta_2$ condition and Corollary~\ref{cor:completecells}, the last term is bounded by
\begin{equation}
\label{BLholecost2}
C\left\{
\int_{\Omega_\eps^{\rm cell,*}}
\Phi(C\eps\abs{\D v})\,\dd x
+
\eps\int_{\Gamma_\eps}\Phi(C\abs{v_\tau})\,\dd\sigma_x
\right\}.
\end{equation}
Combining \eqref{BLstep1a}--\eqref{BLholecost2} with \eqref{scalingPhi} yields
\begin{equation}
\label{BLgradientfinal}
\int_{R_\eps}\Phi(\eps\abs{\nabla v})\,\dd x
\le
C\left\{
\int_{\Omega_\eps}\Phi(c\eps\abs{\D v})\,\dd x
+
\eps\int_{\Gamma_\eps}\Phi(c\abs{v_\tau})\,\dd\sigma_x
\right\}.
\end{equation}

We next control $v$ itself on $R_\eps$.  Choose $\eps_0>0$ so small that $\sqrt N\eps_0<\delta_0$, where $\delta_0$ is given by
Lemma~\ref{lem:stripPoincare}.  By Lemma~\ref{lem:Rstrip} and the identity
$\widehat v_\eps=v$ on $R_\eps$,
\begin{align}
\int_{R_\eps}\Phi(\abs v)\,\dd x
&\le
\int_{\Sigma_{\sqrt N\eps}}\Phi(\abs{\widehat v_\eps})\,\dd x
\le
C\int_{\Sigma_{\sqrt N\eps}}
\Phi(C\eps\abs{\nabla \widehat v_\eps})\,\dd x
\label{BLvalue1}\\
&\le
C\int_\Omega\Phi(C\eps\abs{\nabla \widehat v_\eps})\,\dd x.
\notag
\end{align}
The second line follows from Lemma~\ref{lem:stripPoincare}.  Apply
\eqref{scaledKornOmega}, with $w$ replaced by $C\widehat v_\eps$, to the last integral and repeat the decomposition
\eqref{BLdecompose}--\eqref{BLholecost2}.  We obtain
\begin{equation}
\label{BLvaluefinal}
\int_{R_\eps}\Phi(\abs v)\,\dd x
\le
C\left\{
\int_{\Omega_\eps}\Phi(c\eps\abs{\D v})\,\dd x
+
\eps\int_{\Gamma_\eps}\Phi(c\abs{v_\tau})\,\dd\sigma_x
\right\}.
\end{equation}
Adding \eqref{BLgradientfinal} and \eqref{BLvaluefinal} proves
\eqref{boundarylayerproved}.
\end{proof}

Combining the complete-cell estimate with the boundary-layer control yields the global inequality.

\begin{corollary}
\label{cor:globalkorn}
Under the assumptions of Theorem~\ref{thm:boundarylayer}, there exist $C,c>0$, independent
of $\eps\in(0,\eps_0)$, such that every admissible $v$ satisfying
\eqref{boundaryadmissible} obeys
{\small
\begin{equation}
\label{globalgoal}
\int_{\Omega_\eps}[\Phi(\abs v)+\Phi(\eps\abs{\nabla v})]\dd x\le C\left\{\int_{\Omega_\eps}\Phi(c\eps\abs{\D v})\dd x+\eps\int_{\Gamma_\eps}\Phi(c\abs{v_\tau})\dd\sigma_x\right\}.
\end{equation}
}
\end{corollary}

\begin{proof}
Up to null sets, $\Omega_\eps$ is the disjoint union of
$\Omega_\eps^{\rm cell,*}$ and $R_\eps$.  Apply Corollary~\ref{cor:completecells} on the
first set and Theorem~\ref{thm:boundarylayer} on the second one, then add the two
inequalities and use $\Phi\in\Delta_2$ to unify the fixed constants inside the modulars.
\end{proof}

\subsection{Weak formulation and existence}
\label{sec:microscopic}

In this subsection $\eps\in(0,\eps_0)$ is fixed and constants may depend on $\eps$.  Define the
solenoidal velocity space
{\small
\begin{equation}
\label{Vepsdef}
\mathbb V_\eps:=\{v\in W^1L^{\Phi}(\Omega_\eps)^N:\ \diver v=0\text{ in }\Omega_\eps,\ \gamma_0 v=0\text{ on }\partial\Omega,\ \gamma_0 v\cdot n_\eps=0\text{ on }\Gamma_\eps\}.
\end{equation}
}
Since the trace and divergence operators are continuous, $\mathbb V_\eps$ is a closed subspace of
the reflexive space $W^1L^{\Phi}(\Omega_\eps)^N$.  It is therefore reflexive.  It is convenient
to equip it with the Darcy-scaled norm
\begin{equation}
\label{Vepsnorm}
\norm{v}_{\mathbb V_\eps}:=\norm{v}_{\Phi,\Omega_\eps}
+\norm{\eps\nabla v}_{\Phi,\Omega_\eps}.
\end{equation}
For each fixed $\eps$, this norm is equivalent to the usual norm of $W^1L^{\Phi}(\Omega_\eps)^N$.

The convection is treated through a skew-symmetric form.  Set $r:=2\underline p'=2\underline p/(\underline p-1)$, with $\underline p$ as in \eqref{underlinep}.  By
\eqref{underlinep}, $t^{\underline p}\le C(1+\Phi(t))$ for $t\ge0$.  For every
measurable set $D\subset\mathbb R^N$ of finite measure, applying this to $v/\norm v_{\Phi,D}$ gives
\begin{equation}
\label{OrliczLpembed}
L^{\Phi}(D)\hookrightarrow L^{\underline p}(D),
\qquad\text{hence}\qquad
W^1L^{\Phi}(\Omega_\eps)^N\hookrightarrow W^{1,\underline p}(\Omega_\eps)^N.
\end{equation}
If $\underline p<N$, the condition $\underline p>3N/(N+2)$ is equivalent to $r<\underline p^{\,*}$.
If $\underline p\ge N$, the same conclusion below follows from the usual Sobolev embeddings into every
finite power space.  Thus, for every fixed $\eps$,
\begin{equation}
\label{compactconvective}
\mathbb V_\eps\hookrightarrow\hookrightarrow L^r(\Omega_\eps)^N.
\end{equation}
The embedding constant is allowed to depend on $\eps$ in this existence section.

For $a,b,v\in W^1L^{\Phi}(\Omega_\eps)^N$, define the skew-symmetrized trilinear form
\begin{equation}
\label{convectiveform}
c_\eps(a,b,v):=\frac{\eps^\kappa}{2}\int_{\Omega_\eps}[(a\cdot\nabla)b\cdot v-(a\cdot\nabla)v\cdot b]\,\dd x.
\end{equation}
Since $\frac2r+\frac1{\underline p}=1$,
by H\"older's inequality together with \eqref{OrliczLpembed} and \eqref{compactconvective}, the
form $c_\eps$ is continuous and trilinear on $\big(W^1L^{\Phi}(\Omega_\eps)^N\big)^3$, hence on
$\mathbb V_\eps^3$.  Moreover,
\begin{equation}
\label{convectivezero}
c_\eps(a,b,b)=0\qquad\forall\,a,b\in \mathbb V_\eps,
\end{equation}
in particular $c_\eps(v,v,v)=0$.  For $u,v\in\mathbb V_\eps$, integration by parts gives
\begin{equation}
\label{convectiveIBP}
c_\eps(u,u,v)
=-\eps^\kappa\int_{\Omega_\eps}(u\otimes u):\nabla v\,\dd x,
\end{equation}
so \eqref{convectiveform} is the weak form of the convective term in \eqref{Peps}.  The
skew-symmetric expression is preferable at the existence level because
\eqref{convectivezero} is then algebraic and requires no additional approximation argument.

\begin{lemma}
\label{lem:convweak}
Define $C_\eps:\mathbb V_\eps\to \mathbb V_\eps^*$ by
\[
\langle C_\eps(u),v\rangle:=c_\eps(u,u,v).
\]
If $u_n\rightharpoonup u$ in $\mathbb V_\eps$, then
\begin{equation}
\label{convweakstar}
C_\eps(u_n)\rightharpoonup C_\eps(u)
\qquad\text{in }\mathbb V_\eps^*.
\end{equation}
Furthermore,
\begin{equation}
\label{convdiagzero}
\langle C_\eps(u),u\rangle=0\qquad\forall\,u\in \mathbb V_\eps.
\end{equation}
\end{lemma}

\begin{proof}
By \eqref{compactconvective},
\[
u_n\to u\qquad\text{strongly in }L^r(\Omega_\eps)^N.
\]
Also $\nabla u_n\rightharpoonup\nabla u$ in $L^{\underline p}(\Omega_\eps)^{N\times N}$, because
$W^1L^{\Phi}(\Omega_\eps)^N\hookrightarrow W^{1,\underline p}(\Omega_\eps)^N$.  Fix $v\in \mathbb V_\eps$.  Since
$u_n\to u$ strongly in $L^r(\Omega_\eps)^N$ and $v\in L^r(\Omega_\eps)^N$,
\[
v\otimes u_n\to v\otimes u
\qquad\text{strongly in }L^{\underline p'}(\Omega_\eps)^{N\times N}.
\]
Therefore
\[
\int_{\Omega_\eps}(u_n\cdot\nabla)u_n\cdot v\,\dd x
\longrightarrow
\int_{\Omega_\eps}(u\cdot\nabla)u\cdot v\,\dd x.
\]
Likewise $u_n\otimes u_n\to u\otimes u$ strongly in $L^{\underline p'}(\Omega_\eps)^{N\times N}$, so
\[
\int_{\Omega_\eps}(u_n\cdot\nabla)v\cdot u_n\,\dd x
\longrightarrow
\int_{\Omega_\eps}(u\cdot\nabla)v\cdot u\,\dd x.
\]
This proves \eqref{convweakstar}.  Identity \eqref{convdiagzero} follows immediately from
the definition \eqref{convectiveform}.
\end{proof}

The convective term is used in this skew-symmetric form in the weak formulation.

\begin{definition}
\label{def:weakvelocity}
For fixed $\eps\in(0,\eps_0)$, a field $u_\eps\in \mathbb V_\eps$ is called a weak velocity solution of
\eqref{Peps} if
\begin{equation}
\label{weakPeps}
\begin{aligned}
&\int_{\Omega_\eps}
\mathbb S\!\left(\frac{x}{\eps},\eps\D u_\eps\right):\eps\D v\,\dd x
+\eps\int_{\Gamma_\eps}
\mathcal B\!\left(\frac{x}{\eps},(u_\eps)_\tau\right)\cdot v_\tau\,\dd\sigma_x
\\
&\qquad +c_\eps(u_\eps,u_\eps,v)
=\int_{\Omega_\eps}f\cdot v\,\dd x
\qquad\forall\,v\in \mathbb V_\eps.
\end{aligned}
\end{equation}
\end{definition}

Define $A_\eps:\mathbb V_\eps\to \mathbb V_\eps^*$ by
{\small
\begin{equation}
\label{Aepsdef}
\langle A_\eps(u),v\rangle:=\int_{\Omega_\eps}\mathbb S\!\left(\frac{x}{\eps},\eps\D u\right):\eps\D v\,\dd x+\eps\int_{\Gamma_\eps}\mathcal B\!\left(\frac{x}{\eps},u_\tau\right)\cdot v_\tau\,\dd\sigma_x.
\end{equation}
}
By Proposition~\ref{prop:scaled-trace-prelim}, the trace operator
$W^1L^{\Phi}(\Omega_\eps)^N\to L^{\Phi}(\Gamma_\eps)^N$ is continuous for fixed $\eps$.  The growth assumptions \eqref{Sgrowth} and \eqref{Bgrowth}, generalized
H\"older's inequality, and the $\Delta_2$ conditions therefore show that $A_\eps$ maps bounded
sets into bounded sets of $\mathbb V_\eps^*$.  Continuity of the Nemytskii maps also gives
hemicontinuity.  Indeed, if $t_n\to t$, continuity in the constitutive variables gives pointwise
convergence, while the growth estimates and the bound
$\Phi(\abs{E+t_nG})\le C\big[\Phi(C\abs E)+\Phi(C\abs G)\big]$, valid for $E,G\in\R^{N\times N}_{\rm sym}$ and
$|t_n|\le1$, provide uniform integrability; Vitali's theorem yields the required convergence in the dual
Orlicz spaces.

Monotonicity follows directly from \eqref{Smonotone} and \eqref{Bmonotone}:
\begin{equation}
\label{Amonotone}
\begin{aligned}
\langle A_\eps(u)-A_\eps(w),u-w\rangle
={}&\int_{\Omega_\eps}
\Big[\mathbb S\!\left(\frac{x}{\eps},\eps\D u\right)
-\mathbb S\!\left(\frac{x}{\eps},\eps\D w\right)\Big]
:\eps\D(u-w)\,\dd x\\
&+\eps\int_{\Gamma_\eps}
\Big[\mathcal B\!\left(\frac{x}{\eps},u_\tau\right)
-\mathcal B\!\left(\frac{x}{\eps},w_\tau\right)\Big]
\cdot(u_\tau-w_\tau)\,\dd\sigma_x\\
\ge{}&0.
\end{aligned}
\end{equation}
Thus $A_\eps$ is a bounded monotone hemicontinuous operator and, in particular,
pseudomonotone.

We next verify coercivity in the scaled norm.  For $v\in\mathbb V_\eps$, set
\begin{equation}
\label{Darcyenergydef}
\mathcal E_\eps(v):=
\int_{\Omega_\eps}\Phi(\eps\abs{\D v})\,\dd x
+\eps\int_{\Gamma_\eps}\Phi(\abs{v_\tau})\,\dd\sigma_x.
\end{equation}
By the bulk and surface coercivity estimates \eqref{Scoercive} and \eqref{Bcoercive},
\begin{equation}
\label{Acoerciveenergy}
\langle A_\eps(v),v\rangle
\ge a_0\mathcal E_\eps(v),
\qquad a_0:=\min\{a_S,a_B\}>0.
\end{equation}
On the other hand, Corollary~\ref{cor:globalkorn} and $\Phi\in\Delta_2$ imply
\begin{equation}
\label{energycontrolsmodular}
\int_{\Omega_\eps}
\Big[\Phi(\abs v)+\Phi(\eps\abs{\nabla v})\Big]\,\dd x
\le C_0\mathcal E_\eps(v),
\end{equation}
where $C_0$ is independent of $\eps\in(0,\eps_0)$. For fixed $\eps\in(0,\eps_0)$, the
modular--norm inequalities applied to the left-hand side show that, when
$\norm v_{\mathbb V_\eps}\ge2$,
\begin{equation}
\label{modularcoercivitynorm}
\mathcal E_\eps(v)\ge c\norm v_{\mathbb V_\eps}^{p^-}.
\end{equation}
Consequently
\begin{equation}
\label{operatorcoercivity}
\frac{\langle A_\eps(v),v\rangle}{\norm v_{\mathbb V_\eps}}
\longrightarrow+\infty
\qquad\text{as }\norm v_{\mathbb V_\eps}\to\infty.
\end{equation}

Set
\[
\mathcal T_\eps:=A_\eps+C_\eps:\mathbb V_\eps\longrightarrow \mathbb V_\eps^*.
\]
The skew-symmetric convection preserves the pseudomonotone structure.

\begin{lemma}
\label{lem:Tpseudo}
Under Assumptions~\ref{ass:stress}--\ref{ass:friction} and \eqref{NSexponent},
$\mathcal T_\eps$ is bounded, coercive and pseudomonotone on $\mathbb V_\eps$.
\end{lemma}

\begin{proof}
Boundedness follows from the boundedness of $A_\eps$, the continuity of the trilinear form,
and the fixed-$\eps$ Sobolev embedding.  By \eqref{convdiagzero},
\[
\langle \mathcal T_\eps(v),v\rangle=\langle A_\eps(v),v\rangle,
\]
so coercivity is exactly \eqref{operatorcoercivity}.

We finally prove pseudomonotonicity.  Let
$v_n\rightharpoonup v$ in $\mathbb V_\eps$ and suppose
\begin{equation}
\label{pseudolimassumption}
\limsup_{n\to\infty}
\langle \mathcal T_\eps(v_n),v_n-v\rangle\le0.
\end{equation}
By Lemma~\ref{lem:convweak} and \eqref{convdiagzero},
\[
\langle C_\eps(v_n),v_n-v\rangle
=-\langle C_\eps(v_n),v\rangle\longrightarrow
-\langle C_\eps(v),v\rangle=0.
\]
Hence
\[
\limsup_{n\to\infty}
\langle A_\eps(v_n),v_n-v\rangle\le0.
\]
Since the bounded, monotone and hemicontinuous operator $A_\eps$ is pseudomonotone
\cite[Lemma~2.9]{Roubicek2013}, for every
$w\in \mathbb V_\eps$,
\begin{equation}
\label{Apseudoineq}
\langle A_\eps(v),v-w\rangle
\le\liminf_{n\to\infty}
\langle A_\eps(v_n),v_n-w\rangle.
\end{equation}
Furthermore,
\[
\langle C_\eps(v_n),v_n-w\rangle
=-\langle C_\eps(v_n),w\rangle
\longrightarrow
-\langle C_\eps(v),w\rangle
=\langle C_\eps(v),v-w\rangle.
\]
Adding this convergence to \eqref{Apseudoineq} yields
\[
\langle \mathcal T_\eps(v),v-w\rangle
\le\liminf_{n\to\infty}
\langle \mathcal T_\eps(v_n),v_n-w\rangle,
\]
which is the definition of pseudomonotonicity.
\end{proof}

These properties give existence at fixed $\eps$.

\begin{theorem}
\label{thm:microexistence}
Assume Assumptions~\ref{ass:geometry}--\ref{ass:force-exponent}, and let $\kappa\in\R$.
Then, for every $\eps\in(0,\eps_0)$, there exists
$u_\eps\in \mathbb V_\eps$ which is a weak velocity solution in the sense of Definition~\ref{def:weakvelocity}, equivalently satisfying \eqref{weakPeps}.
\end{theorem}

\begin{proof}
The linear functional $F_\eps(v):=\int_{\Omega_\eps}f\cdot v\,\dd x$, $v\in\mathbb V_\eps$, belongs to $\mathbb V_\eps^*$ by generalized H\"older's inequality.  The space $\mathbb V_\eps$ is separable and reflexive,
and Lemma~\ref{lem:Tpseudo} shows that $\mathcal T_\eps$ is bounded, coercive and
pseudomonotone.  The surjectivity theorem for coercive pseudomonotone operators
\cite[Theorem~2.6]{Roubicek2013} therefore gives $u_\eps\in \mathbb V_\eps$ such that
\[
\mathcal T_\eps(u_\eps)=F_\eps\quad\text{in }\mathbb V_\eps^*.
\]
This is exactly \eqref{weakPeps}.
\end{proof}

\begin{remark}
\label{rem:pthreshold}
For fixed $\eps$, the number $\eps^\kappa$ is a finite positive coefficient of the skew-symmetric
convection, whose diagonal contribution vanishes; Theorem~\ref{thm:microexistence} therefore holds for
every $\kappa\in\R$, and the restriction on $\kappa$ enters only when $\eps\to0$ and microscopic test
gradients have size $O(\eps^{-1})$.  The threshold \eqref{NSexponent} is used only to treat convection by
the compact Sobolev embedding \eqref{compactconvective} (see also Remark~\ref{rem:energyinequality}).
\end{remark}

\begin{remark}
\label{rem:fixedpressure}
For fixed $\eps$, $\Omega_\eps$ is a bounded John domain, hence satisfies the emanating chain condition
\cite[Theorem~3.8]{DieningRuzickaSchumacher2010}, and the Orlicz Bogovskii operator
\cite[Theorem~6.6]{DieningRuzickaSchumacher2010} yields, by the standard de Rham argument, a unique
$p_\eps\in L^{\widetilde\Phi}(\Omega_\eps)$ with $\int_{\Omega_\eps}p_\eps\,\dd x=0$ such that $(u_\eps,p_\eps)$
satisfies the weak form of \eqref{Peps} for every impermeable test field vanishing on $\partial\Omega$.
The Bogovskii constant depends on the perforated geometry, hence on $\eps$.  The homogenization argument
below is therefore carried out in pressure-free solenoidal form; both limit pressures are reconstructed from
the limit system in Section~\ref{sec:effectivedarcy}, where the extended microscopic pressure is also shown to
converge (Theorem~\ref{thm:pressureconv}).
\end{remark}

\subsection{Uniform Darcy-scale estimates}
\label{subsec:uniform}

The $\eps$-uniform Darcy-scale energy estimate provides the compactness required for homogenization.

\begin{theorem}
\label{thm:uniformenergy}
Under the assumptions of Theorem~\ref{thm:microexistence}, let $u_\eps$ be any weak
velocity solution of \eqref{Peps} in the sense of Definition~\ref{def:weakvelocity}, $0<\eps<\eps_0$.  There exists a constant $C>0$, depending
only on the structural data, $\Omega$, $Y^*$, $T$, and $f$, but not on
$\eps\in(0,\eps_0)$, such that
{\small
\begin{equation}
\label{uniformenergy}
\int_{\Omega_\eps}\!\Phi(\abs{u_\eps})\dd x+\int_{\Omega_\eps}\!\Phi(\eps\abs{\nabla u_\eps})\dd x+\int_{\Omega_\eps}\!\Phi(\eps\abs{\D u_\eps})\dd x+\eps\int_{\Gamma_\eps}\!\Phi(\abs{(u_\eps)_\tau})\dd\sigma_x\le C.
\end{equation}
}
Moreover, the nonlinear bulk and surface fluxes satisfy the dual modular bounds
\begin{align}
\label{stressdualbound}
\int_{\Omega_\eps}
\widetilde\Phi\!\left(
\frac{\abs{\mathbb S(x/\eps,\eps\D u_\eps)}}{c_S}
\right)\dd x&\le C,\\
\label{frictiondualbound}
\eps\int_{\Gamma_\eps}
\widetilde\Phi\!\left(
\frac{\abs{\mathcal B(x/\eps,(u_\eps)_\tau)}}{c_B}
\right)\dd\sigma_x&\le C.
\end{align}
\end{theorem}

\begin{proof}
Take $v=u_\eps$ in \eqref{weakPeps}.  This is legitimate because
$u_\eps\in \mathbb V_\eps$.  By \eqref{convdiagzero}, the convective contribution is exactly zero.
Using \eqref{Scoercive} and \eqref{Bcoercive} gives
\begin{equation}
\label{energyfirst}
a_0\mathcal E_\eps(u_\eps)
\le
\int_{\Omega_\eps}f\cdot u_\eps\,\dd x,
\end{equation}

with $a_0$ as in \eqref{Acoerciveenergy}.  Fix $0<\eta<1$.  By \eqref{youngPhi} and \eqref{scalingPhi},
\[
\abs{f\cdot u_\eps}
\le
\Phi(\eta\abs{u_\eps})
+\widetilde\Phi\!\left(\frac{\abs f}{\eta}\right)
\le\eta^{p^-}\Phi(\abs{u_\eps})+\widetilde\Phi\!\left(\frac{\abs f}{\eta}\right).
\]
Consequently, by \eqref{energycontrolsmodular} applied to $u_\eps$,
\begin{equation}
\label{forcingabsorb}
\begin{aligned}
\int_{\Omega_\eps}f\cdot u_\eps\,\dd x
&\le
\eta^{p^-}\int_{\Omega_\eps}\Phi(\abs{u_\eps})\,\dd x
+\int_\Omega\widetilde\Phi\!\left(\frac{\abs f}{\eta}\right)\dd x\\
&\le
C_0\eta^{p^-}\mathcal E_\eps(u_\eps)
+\int_\Omega\widetilde\Phi\!\left(\frac{\abs f}{\eta}\right)\dd x.
\end{aligned}
\end{equation}
Choose $\eta\in(0,1)$, independently of $\eps$, so small that $C_0\eta^{p^-}\le a_0/2$.
Combining \eqref{energyfirst} and \eqref{forcingabsorb} gives
\begin{equation}
\label{energyEbound}
\mathcal E_\eps(u_\eps)
\le
\frac{2}{a_0}
\int_\Omega\widetilde\Phi\!\left(\frac{\abs f}{\eta}\right)\dd x
=:C_f.
\end{equation}
The number $C_f$ is finite because $f\in L^{\widetilde\Phi}(\Omega)^N$ and
$\widetilde\Phi\in\Delta_2$, and it is independent of $\eps$.

The bounds for
$\int\Phi(\eps\abs{\D u_\eps})$ and the surface modular follow directly from
\eqref{energyEbound}; the first two terms of \eqref{uniformenergy} then follow from
\eqref{energycontrolsmodular}, and the remaining term is controlled because
$\abs{\D u_\eps}\le\abs{\nabla u_\eps}$.

Finally, \eqref{Sgrowth} and \eqref{energyEbound} give
\[
\int_{\Omega_\eps}
\widetilde\Phi\!\left(
\frac{\abs{\mathbb S(x/\eps,\eps\D u_\eps)}}{c_S}
\right)\dd x
\le
C_S\int_{\Omega_\eps}\Phi(\eps\abs{\D u_\eps})\dd x\le C.
\]
Likewise, \eqref{Bgrowth} yields
\[
\eps\int_{\Gamma_\eps}
\widetilde\Phi\!\left(
\frac{\abs{\mathcal B(x/\eps,(u_\eps)_\tau)}}{c_B}
\right)\dd\sigma_x
\le
C_B\eps\int_{\Gamma_\eps}\Phi(\abs{(u_\eps)_\tau})\dd\sigma_x\le C.
\]
This proves \eqref{stressdualbound}--\eqref{frictiondualbound}.
\end{proof}

These uniform estimates are those required by the volume and surface two-scale framework of
Section~\ref{sec:preliminaries}; we now pass to the limit.


\section{Homogenization: compactness, nonlinear identification, and limit equation}
\label{sec:homogenization}

Throughout this section, Assumptions~\ref{ass:geometry}--\ref{ass:force-exponent} hold,
$\mathscr E\subset(0,\eps_0)$, and $(u_\eps)_{\eps\in\mathscr E}$ is any family of weak solutions of
\eqref{Peps} in the sense of Definition~\ref{def:weakvelocity}; from Subsection~\ref{subsec:inertia} on,
Assumption~\ref{ass:kappa} holds as well.  We first establish the compactness of the velocities and of their
traces, then show that inertia vanishes, construct exactly solenoidal recovery tests, identify the nonlinear
limits by a coupled Minty argument, and finally derive the corresponding strong corrector result.

\subsection{Darcy-scaled two-scale compactness}
\label{subsec:compactness}

The energy estimate of Theorem~\ref{thm:uniformenergy} controls $u_\eps$ in $L^{\Phi}(\Omega_\eps)^N$ and the
scaled gradient $\eps\nabla u_\eps$ in $L^{\Phi}(\Omega_\eps)^{N\times N}$, but it does not control the unscaled gradient.
The corresponding compactness statement is therefore different from the usual
unscaled two-scale gradient theorem in $W^1L^{\Phi}(\Omega)$.  At the Darcy scale the whole
leading-order gradient is microscopic: the limit of $\eps\nabla u_\eps$ is a $y$-gradient,
with no macroscopic $x$-gradient term.

We first isolate the Darcy-scaled compactness mechanism on the fixed domain $\Omega$, before restricting
the limit to the fluid phase.  The space $\mathbb W^\Phi_{\rm per}(\Omega\times Y)$ is defined as in
\eqref{productOrliczSpace}, with $Y$ in place of $Y^*$.

\begin{theorem}
\label{thm:fixedscaledcompactness}
Let $(v_\eps)_{\eps\in\mathscr E}$ be a family in $W^1L^{\Phi}(\Omega)^N$ satisfying
\begin{equation}
\label{fixedscaledbound}
\sup_{\eps\in\mathscr E}\int_\Omega
\left[\Phi(\abs{v_\eps})+\Phi(\eps\abs{\nabla v_\eps})\right]\dd x<\infty.
\end{equation}
Then there exist a subsequence, not relabelled, and
$v_0\in\mathbb W^\Phi_{\rm per}(\Omega\times Y)^N$ such that
\begin{align}
\label{fixedUtwoscale}
v_\eps&\overset{2s}{\rightharpoonup}v_0
&&\text{weakly two-scale in }L^{\Phi}(\Omega)^N,\\
\label{fixedgradtwoscale}
\eps\nabla v_\eps&\overset{2s}{\rightharpoonup}\nabla_yv_0
&&\text{weakly two-scale in }L^{\Phi}(\Omega)^{N\times N},\\
\label{fixedDtwoscale}
\eps\D v_\eps&\overset{2s}{\rightharpoonup}\D_yv_0
&&\text{weakly two-scale in }L^{\Phi}(\Omega)_{\rm sym}^{N\times N},
\end{align}
where
\[
\D_yv_0:=\frac12\big(\nabla_yv_0+(\nabla_yv_0)^T\big).
\]
\end{theorem}

\begin{proof}
By the modular--norm equivalence under $\Delta_2$, \eqref{fixedscaledbound} implies that
$(v_\eps)$ and each component of $(\eps\nabla v_\eps)$ are bounded in $L^{\Phi}(\Omega)$.
Theorem~\ref{thm:volume-two-scale-compactness}, applied componentwise and followed by
a diagonal extraction, yields
\begin{equation}
\label{preliminaryscaledlimits}
v_\eps\overset{2s}{\rightharpoonup}v_0,
\qquad
\eps D_{x_j}v_{\eps,i}\overset{2s}{\rightharpoonup}G_{ij},
\quad 1\le i,j\le N,
\end{equation}
for some
\[
v_0\in L^{\Phi}(\Omega\times Y)^N,
\qquad
G=(G_{ij})\in L^{\Phi}(\Omega\times Y)^{N\times N}.
\]
We identify $G$ as follows.

Fix $1\le i,j\le N$ and $\psi\in C_c^\infty\big(\Omega;C_{\rm per}^\infty(Y)\big)$.
Since the $x$-support of $\psi$ is compactly contained in $\Omega$, integration by parts on the
fixed domain gives
{\small
\begin{equation}
\label{scaledIBP}
\int_\Omega \eps D_{x_j}v_{\eps,i}\,\psi\!\left(x,\frac{x}{\eps}\right)\dd x=-\eps\int_\Omega v_{\eps,i}D_{x_j}\psi\!\left(x,\frac{x}{\eps}\right)\dd x-\int_\Omega v_{\eps,i}D_{y_j}\psi\!\left(x,\frac{x}{\eps}\right)\dd x.
\end{equation}
}
The first term on the right tends to zero.  Indeed, generalized H\"older's inequality, the uniform
$L^{\Phi}(\Omega)$ bound of $v_\eps$, and periodic sampling in $L^{\widetilde\Phi}(\Omega)$ give
\[
\left|\eps\int_\Omega v_{\eps,i}
D_{x_j}\psi\!\left(x,\frac{x}{\eps}\right)\dd x\right|
\le C\eps\norm{v_{\eps,i}}_{\Phi,\Omega}
\norm{D_{x_j}\psi(\cdot,\cdot/\eps)}_{\widetilde\Phi,\Omega}
\longrightarrow0.
\]
Passing to the two-scale limit in the remaining terms of \eqref{scaledIBP} yields
\begin{equation}
\label{ydistributionidentity}
\int_\Omega\int_Y G_{ij}(x,y)\psi(x,y)\dd y\dd x
=-\int_\Omega\int_Y v_{0,i}(x,y)D_{y_j}\psi(x,y)\dd y\dd x.
\end{equation}
Thus $G_{ij}=D_{y_j}v_{0,i}$ in the sense of distributions on $\Omega\times\R^N$, for the $Y$-periodic extensions in $y$.  Taking
$\psi(x,y)=\eta(x)\varpi(y)$, with $\eta$ and $\varpi$ in countable dense subsets of $C_c^\infty(\Omega)$ and
$C^\infty_{\rm per}(Y)$, shows that, outside a common null set of $x\in\Omega$,
$\int_YG_{ij}(x,\cdot)\varpi\,\dd y=-\int_Yv_{0,i}(x,\cdot)D_{y_j}\varpi\,\dd y$ for all such $\varpi$, hence
for all $\varpi\in C^\infty_{\rm per}(Y)$ by density; thus $v_0(x,\cdot)\in W^1_{\rm per}L^\Phi(Y)^N$ with
$\nabla_yv_0(x,\cdot)=G(x,\cdot)$ for almost every $x\in\Omega$.  Since $G_{ij}\in L^{\Phi}(\Omega\times Y)$, this
proves that $v_0\in\mathbb W^\Phi_{\rm per}(\Omega\times Y)^N$ and establishes
\eqref{fixedgradtwoscale}.  Taking symmetric parts, which is a continuous linear operation,
gives \eqref{fixedDtwoscale}.
\end{proof}

We next transfer the fixed-domain convergence to the fluid region.  Because only complete cells are
perforated, the difference between $\chi_\eps$ and
$\chi_{Y^*}(x/\eps)$ is supported in the boundary remainder $R_\eps$, whose measure tends to
zero.

\begin{lemma}
\label{lem:phaseindicator}
If $(z_\eps)$ is bounded in $L^{\Phi}(\Omega)^m$ and $z_\eps\overset{2s}{\rightharpoonup}z_0$ in
$L^{\Phi}(\Omega)^m$, then
\begin{equation}
\label{indicatorproduct}
\chi_\eps z_\eps\overset{2s}{\rightharpoonup}
\chi_{Y^*}(y)z_0(x,y)
\quad\text{in }L^{\Phi}(\Omega)^m.
\end{equation}
\end{lemma}

\begin{proof}
Let $\psi\in C_c(\Omega;C_{\rm per}(Y))^m$.  For small $\eps$, the $x$-support of $\psi$ meets only complete
cells, on which $\chi_\eps=\chi_{Y^*}(x/\eps)$; hence $\chi_\eps\psi^\eps=(\chi_{Y^*}\psi)^\eps$.  Since
$y\mapsto\psi(\cdot,y)$ is continuous into $C(\overline\Omega)^m$, the field $\chi_{Y^*}\psi$ belongs to
$L^\infty_{\rm per}(Y;C(\overline\Omega))^m$ and vanishes outside the $x$-support of $\psi$, and
Lemma~\ref{lem:Lqadmissible} shows
that $\chi_\eps\psi^\eps\to\chi_{Y^*}\psi$ strongly two-scale in $L^{\widetilde\Phi}(\Omega)^m$, and
Lemma~\ref{lem:weak-strong-two-scale-pairing} gives
$\int_\Omega\chi_\eps z_\eps\cdot\psi^\eps\,\dd x\to\int_\Omega\int_Y\chi_{Y^*}z_0\cdot\psi\,\dd y\dd x$.  For a
general $\psi\in L^{\widetilde\Phi}(\Omega;C_{\rm per}(Y))^m$, approximate $\psi$ by $\psi_j\in C_c(\Omega;C_{\rm per}(Y))^m$
in that space; since $|(\psi-\psi_j)^\eps|\le\sup_{y\in Y}|\psi(\cdot,y)-\psi_j(\cdot,y)|$, \eqref{OrliczHolder} and the
bound on $(z_\eps)$ make the replacement error small uniformly in $\eps$.
\end{proof}

For the weak solutions $u_\eps$, let
\[
\bar u_\eps:=P_\eps u_\eps\in W^1_0L^{\Phi}(\Omega)^N
\]
be the componentwise hole-filling extension of Proposition~\ref{prop:uniform-extension-prelim}.  We also define
the zero extensions to $\Omega$ by
\begin{equation}
\label{zeroextensions}
\widetilde u_\eps:=\chi_\eps \bar u_\eps,
\qquad
\widetilde{\nabla u_\eps}:=\chi_\eps\nabla \bar u_\eps,
\qquad
\widetilde{\D u_\eps}:=\chi_\eps\D \bar u_\eps.
\end{equation}
These coincide with $u_\eps$, $\nabla u_\eps$ and $\D u_\eps$ on $\Omega_\eps$ and vanish
inside the holes.

\begin{theorem}
\label{thm:darcyscaledcompactness}
The extensions satisfy
\begin{equation}
\label{extensionfixedbound}
\sup_{\eps\in\mathscr E}\int_\Omega
\left[\Phi(\abs{\bar u_\eps})+\Phi(\eps\abs{\nabla \bar u_\eps})\right]\dd x<\infty,
\end{equation}
and, after extraction, there exists
\begin{equation}
\label{u0space}
u_0\in \mathbb W^\Phi_{\rm per}(\Omega\times Y^*)^N,
\end{equation}
that is,
$u_0\in L^{\Phi}(\Omega\times Y^*)^N$ and
$\nabla_yu_0\in L^{\Phi}(\Omega\times Y^*)^{N\times N}$ with periodic outer-face traces, such that
\begin{align}
\label{u0twoscale}
\widetilde u_\eps
&\overset{2s}{\rightharpoonup}
\chi_{Y^*}(y)u_0(x,y),\\
\label{grad0twoscale}
\eps\widetilde{\nabla u_\eps}
&\overset{2s}{\rightharpoonup}
\chi_{Y^*}(y)\nabla_yu_0(x,y),\\
\label{D0twoscale}
\eps\widetilde{\D u_\eps}
&\overset{2s}{\rightharpoonup}
\chi_{Y^*}(y)\D_yu_0(x,y).
\end{align}
All convergences are weak two-scale in $L^{\Phi}(\Omega)^N$ or $L^{\Phi}(\Omega)^{N\times N}$.
\end{theorem}

\begin{proof}
We first prove the uniform fixed-domain bounds for $\bar u_\eps=P_\eps u_\eps$.  The value estimate in Proposition~\ref{prop:uniform-extension-prelim}, applied componentwise, gives
\begin{equation}
\label{extensionvaluebound}
\int_\Omega\Phi(\abs{\bar u_\eps})\dd x
\le C\int_{\Omega_\eps}\Phi(C\abs{u_\eps})\dd x\le C,
\end{equation}
where the last inequality follows from \eqref{uniformenergy} and $\Phi\in\Delta_2$.
For the scaled gradient, split the fixed domain into the fluid region and the holes.  Since
$\bar u_\eps=u_\eps$ on $\Omega_\eps$,
\begin{align}
\int_\Omega\Phi(\eps\abs{\nabla \bar u_\eps})\dd x
={}&\int_{\Omega_\eps}\Phi(\eps\abs{\nabla u_\eps})\dd x
+\sum_{k\in K_\eps}
\int_{\eps(k+T)}\Phi(\eps\abs{\nabla \bar u_\eps})\dd x.
\label{extensiongraddecomp}
\end{align}
The first term is uniformly bounded by \eqref{uniformenergy}; the second is controlled by the
local extension estimate \eqref{holeextensionlocal} and then again by
\eqref{uniformenergy}, which proves \eqref{extensionfixedbound}.

Apply Theorem~\ref{thm:fixedscaledcompactness}.  There exist
$\bar u_0\in\mathbb W^\Phi_{\rm per}(\Omega\times Y)^N$ and a subsequence such that
\begin{equation}
\label{extendedcompactness}
\bar u_\eps\overset{2s}{\rightharpoonup}\bar u_0,
\qquad
\eps\nabla \bar u_\eps\overset{2s}{\rightharpoonup}\nabla_y\bar u_0,
\qquad
\eps\D \bar u_\eps\overset{2s}{\rightharpoonup}\D_y\bar u_0.
\end{equation}
Define
\[
u_0(x,y):=\bar u_0(x,y)\big|_{Y^*}.
\]
Lemma~\ref{lem:phaseindicator}, applied successively to the three bounded sequences in
\eqref{extendedcompactness}, yields
\eqref{u0twoscale}--\eqref{D0twoscale}.  The function-space properties in
\eqref{u0space} follow by restriction from $Y$ to $Y^*$.  The limit $\bar u_0$ inside $T$ may depend on the
extension, but $u_0$ does not, since $\widetilde u_\eps$ and $\eps\widetilde{\nabla u_\eps}$ depend only on $u_\eps$.
\end{proof}

\begin{corollary}
\label{cor:limitincompressibility}
Let $u_0$ be obtained in Theorem~\ref{thm:darcyscaledcompactness}.  Then for almost every
$x\in\Omega$,
\begin{equation}
\label{fastsolenoidalweak}
\int_{Y^*}u_0(x,y)\cdot\nabla_y\psi(y)\dd y=0
\qquad\forall\,\psi\in C_{\rm per}^\infty(Y).
\end{equation}
Consequently
\begin{equation}
\label{fastdivnormal}
\operatorname{div}_yu_0=0\quad\text{in }Y^*,
\qquad
\gamma_0^y u_0\cdot n=0\quad\text{on }\partial T
\end{equation}
in the weak periodic sense.  Moreover, the filtration velocity of Theorem~A,
\begin{equation}
\label{filtrationvelocity}
U(x)=\int_{Y^*}u_0(x,y)\dd y,
\end{equation}
satisfies
\begin{equation}
\label{macroincompressibility}
\diver_x U=0\quad\text{in }\Omega,
\qquad
U\cdot\nu=0\quad\text{on }\partial\Omega
\end{equation}
in the weak sense.
\end{corollary}

\begin{proof}
Fix $\eta\in C_c^\infty(\Omega)$ and
$\psi\in C_{\rm per}^\infty(Y)$.  Use the scalar test
\[
\omega_\eps(x):=\eps\eta(x)\psi\!\left(\frac{x}{\eps}\right)
\]
in the distributional identity $\diver u_\eps=0$ on $\Omega_\eps$.  Integration by parts
produces no boundary contribution: $u_\eps=0$ on $\partial\Omega$ and
$u_\eps\cdot n_\eps=0$ on $\Gamma_\eps$.  Hence
\begin{equation}
\label{microdivtest}
0=\eps\int_{\Omega_\eps}u_\eps\cdot\nabla\eta\,
\psi\!\left(\frac{x}{\eps}\right)\dd x
+\int_{\Omega_\eps}\eta(x)u_\eps(x)\cdot
\nabla_y\psi\!\left(\frac{x}{\eps}\right)\dd x.
\end{equation}
The first term tends to zero by the uniform bound of $u_\eps$ in $L^{\Phi}(\Omega_\eps)^N$.  Passing to the limit in the second
term by \eqref{u0twoscale} gives
\[
\int_\Omega\eta(x)\int_{Y^*}u_0(x,y)\cdot\nabla_y\psi(y)\dd y\dd x=0.
\]
The fundamental lemma in $x$ proves \eqref{fastsolenoidalweak}.  Since
$u_0(x,\cdot)\in W^1L^{\Phi}(Y^*)^N$ and is periodic across opposite outer faces, the weak Green
formula on $Y^*$ shows that \eqref{fastsolenoidalweak} is equivalent to
\eqref{fastdivnormal}: the outer-face fluxes cancel by periodicity, while the remaining boundary
is $\partial T$.

For the macroscopic statement, let $\eta\in C^\infty(\overline\Omega)$.  The same divergence
identity, now tested directly by $\eta(x)$, gives
\[
\int_{\Omega_\eps}u_\eps\cdot\nabla\eta\dd x=0.
\]
Passing to the two-scale limit yields
\[
\int_\Omega U(x)\cdot\nabla\eta(x)\dd x=0
\qquad\forall\,\eta\in C^\infty(\overline\Omega).
\]
Because $\widetilde\Phi\in\Delta_2$ and $\Omega$ is $C^{1,1}$, smooth functions on
$\overline\Omega$ are dense in $W^1L^{\widetilde\Phi}(\Omega)$; see the Orlicz--Sobolev density
results of \cite[Theorem~2.3]{DonaldsonTrudinger1971}.  Since
$U\in L^\Phi(\Omega)^N$, generalized Orlicz H\"older inequality therefore extends the
preceding identity to every $\eta\in W^1L^{\widetilde\Phi}(\Omega)$.  This is precisely
\eqref{macroincompressibility} in weak form and, in particular, encodes both
$\diver_x U=0$ in $\Omega$ and $U\cdot\nu=0$ on $\partial\Omega$.
\end{proof}

\subsection{Darcy-scaled two-scale trace compactness on periodic surfaces}

Corollary~\ref{cor:limitincompressibility} identifies the normal component of the limit velocity on the
pore walls; the full trace is recovered from the pair formed by the volume limit and the scaled-gradient
limit, through a Green formula which identifies the surface limit with the $y$-trace of the limit.  Since $\eps\nabla v_\eps$
survives at leading order, it reconstructs the $y$-trace of the limit directly.

For a matrix field $\mathbf\Psi=(\Psi_{\alpha i})_{1\le\alpha\le m,\,1\le i\le N}$ we use the
row-wise divergences
\[
(\diver_x\mathbf\Psi)_\alpha:=\sum_{i=1}^N D_{x_i}\Psi_{\alpha i},
\qquad
(\diver_y\mathbf\Psi)_\alpha:=\sum_{i=1}^N D_{y_i}\Psi_{\alpha i}.
\]
On $\partial T$, $\mathbf\Psi n$ denotes the vector in $\R^m$ whose $\alpha$-th component is
$\sum_{i=1}^N\Psi_{\alpha i}n_i$.

\begin{lemma}
\label{lem:normalfluxdense}
Let $m\in\mathbb N$.  The set
\begin{equation}
\label{normalfluxclass}
\mathcal N_m:=\left\{
(\mathbf\Psi n)|_{\partial T}:\
\mathbf\Psi\in C_{\rm per}^{\infty}(Y)^{m\times N}
\right\}
\end{equation}
is dense in $C(\partial T)^m$ for the uniform norm.  Consequently finite sums of functions
$\eta(x)a(y)$, with $\eta\in C_c^\infty(\Omega)$ and $a\in\mathcal N_m$, are dense in
$C_c(\Omega;C(\partial T))^m$ for the natural uniform topology on compact
$x$-supports.
\end{lemma}

\begin{proof}
Let $\psi\in C(\partial T)^m$.  Since $\partial T$ is of class $C^{1,1}$, its unit normal
$n$ is continuous.  Define on $\partial T$ the continuous matrix field $A(y):=\psi(y)\otimes n(y)$, so that $A(y)n(y)=\psi(y)$ because $\abs{n(y)}=1$.  Since
$\overline T\subset Y$, a continuous extension of $A$ can be chosen in a neighbourhood of
$\partial T$ whose closure is compactly contained in $Y$.  Multiplication by a cut-off which is
identically one near $\partial T$, followed by standard mollification, gives matrices
$A_j\in C_c^\infty(Y)^{m\times N}$ such that
\[
\sup_{y\in\partial T}\abs{A_j(y)-A(y)}\longrightarrow0.
\]
Regarded as periodic functions (they vanish near $\partial Y$), the fields $A_j$ belong to
$C_{\rm per}^\infty(Y)^{m\times N}$ and
\[
\sup_{y\in\partial T}\abs{A_j(y)n(y)-\psi(y)}
\le \sup_{y\in\partial T}\abs{A_j(y)-A(y)}\longrightarrow0.
\]
This proves density in $C(\partial T)^m$.  The $x$-dependent assertion follows from the density
of finite sums of separated variables in continuous Bochner-valued functions with compact
$x$-support.
\end{proof}

For a family $(v_\eps)$ we write
\[
\widetilde v_\eps:=\chi_\eps v_\eps,
\qquad
\widetilde{\nabla v_\eps}:=\chi_\eps\nabla v_\eps
\]
for the zero extensions to $\Omega$.

\begin{theorem}
\label{thm:darcyscaledtrace}
Let Assumption~\ref{ass:orlicz} hold, let $m\in\mathbb N$, and let
$v_\eps\in W^1L^{\Phi}(\Omega_\eps)^m$ satisfy
\begin{equation}
\label{darcyscaledtracebound}
\sup_{\eps\in\mathscr E}
\int_{\Omega_\eps}
\left[
\Phi(\abs{v_\eps})+\Phi(\eps\abs{\nabla v_\eps})
\right]\dd x<\infty.
\end{equation}
Suppose that, along a subsequence,
\begin{align}
\label{tracevolumeconv1}
\widetilde v_\eps
&\overset{2s}{\rightharpoonup}
\chi_{Y^*}(y)v_0(x,y),\\
\label{tracevolumeconv2}
\eps\widetilde{\nabla v_\eps}
&\overset{2s}{\rightharpoonup}
\chi_{Y^*}(y)\nabla_yv_0(x,y),
\end{align}
where
\begin{equation}
\label{tracev0space}
v_0\in L^{\Phi}(\Omega\times Y^*)^m,
\qquad
\nabla_yv_0\in L^{\Phi}(\Omega\times Y^*)^{m\times N},
\end{equation}
and $v_0$ has periodic traces across opposite faces of $Y$.

Then the $y$-trace $\gamma_0^y v_0$ belongs to
$L^{\Phi}(\Omega\times\partial T)^m$, the microscopic traces satisfy
\begin{equation}
\label{tracesurfacebound}
\sup_{\eps\in\mathscr E}
\eps\int_{\Gamma_\eps}
\Phi(\abs{v_\eps})\dd\sigma_x<\infty,
\end{equation}
and, along the same subsequence,
\begin{equation}
\label{darcyscaledtraceconv}
v_\eps
\overset{2s,\Gamma}{\rightharpoonup}
\gamma_0^y v_0
\quad\text{weakly in }L^{\Phi}(\Omega\times\partial T)^m.
\end{equation}
Equivalently,
\begin{equation}
\label{darcyscaledtraceintegral}
\lim_{\eps\to0}
\eps\int_{\Gamma_\eps}
v_\eps(x)\cdot
\psi\!\left(x,\frac{x}{\eps}\right)\dd\sigma_x
=
\int_\Omega\int_{\partial T}
\gamma_0^y v_0(x,y)\cdot\psi(x,y)\dd\sigma_y\dd x
\qquad\forall\,\psi\in C_c(\Omega;C(\partial T))^m.
\end{equation}
\end{theorem}

\begin{proof}
We first prove the surface bound and extract a surface limit.  Proposition~\ref{prop:scaled-trace-prelim}, applied componentwise and combined with
\eqref{darcyscaledtracebound}, gives
\[
\eps\int_{\Gamma_\eps}\Phi(\abs{v_\eps})\dd\sigma_x
\le C\int_{\Omega_\eps}
\left[
\Phi(C\abs{v_\eps})+\Phi(C\eps\abs{\nabla v_\eps})
\right]\dd x
\le C.
\]
Thus \eqref{tracesurfacebound} holds.  Theorem~\ref{thm:surface-two-scale-compactness}
yields, after a further extraction if necessary, a field $z_0\in L^{\Phi}(\Omega\times\partial T)^m$ such that
\begin{equation}
\label{preliminarysurfacelimit}
v_\eps\overset{2s,\Gamma}{\rightharpoonup}z_0.
\end{equation}
We identify $z_0$ below.

The fixed-cell trace inequality, integrated with respect to $x$, gives from
\eqref{tracev0space}
\begin{equation}
\label{limittraceintegrability}
\int_\Omega\int_{\partial T}
\Phi(\abs{\gamma_0^y v_0})\dd\sigma_y\dd x
\le C\int_\Omega\int_{Y^*}
\left[
\Phi(C\abs{v_0})+\Phi(C\abs{\nabla_yv_0})
\right]\dd y\dd x<\infty.
\end{equation}
Hence $\gamma_0^y v_0\in L^{\Phi}(\Omega\times\partial T)^m$.

To identify $z_0$, we use normal-flux tests.  Let $\mathbf\Psi\in C_c^\infty\big(\Omega;C_{\rm per}^\infty(Y)^{m\times N}\big)$ and set
\[
\psi(x,y):=\mathbf\Psi(x,y)n(y)
\qquad (y\in\partial T).
\]
For shortness write
$\mathbf\Psi^\eps(x)=\mathbf\Psi(x,x/\eps)$.  Applying the divergence theorem row by row to
$v_\eps\mathbf\Psi^\eps$ on $\Omega_\eps$ gives no contribution on $\partial\Omega$, because
the $x$-support of $\mathbf\Psi$ is compactly contained in $\Omega$.  On $\Gamma_\eps$, the
outward normal of the fluid is $n_\eps(x)=n(x/\eps)$.  Therefore
\begin{align}
\label{surfacegreenidentity}
&\eps\int_{\Gamma_\eps}
v_\eps(x)\cdot
\left(\mathbf\Psi\!\left(x,\frac{x}{\eps}\right)n_\eps(x)\right)
\dd\sigma_x\notag\\
&\quad=
\int_{\Omega_\eps}
\eps\nabla v_\eps:\mathbf\Psi^\eps\dd x
+
\eps\int_{\Omega_\eps}
v_\eps\cdot(\diver_x\mathbf\Psi)^\eps\dd x
+
\int_{\Omega_\eps}
v_\eps\cdot(\diver_y\mathbf\Psi)^\eps\dd x.
\end{align}
Indeed
\[
\diver_x\mathbf\Psi^\eps
=(\diver_x\mathbf\Psi)^\eps
+\frac1\eps(\diver_y\mathbf\Psi)^\eps
\]
row by row, which is exactly the origin of the three terms on the right-hand side.

We now pass to the volume two-scale limit.  The middle term in \eqref{surfacegreenidentity} tends to zero.  Generalized H\"older's inequality,
the uniform bound of $v_\eps$ in $L^{\Phi}(\Omega_\eps)^m$ and periodic sampling give
\[
\left|
\eps\int_{\Omega_\eps}
v_\eps\cdot(\diver_x\mathbf\Psi)^\eps\dd x
\right|
\le C\eps\norm{v_\eps}_{\Phi,\Omega_\eps}
\norm{(\diver_x\mathbf\Psi)^\eps}_{\widetilde\Phi,\Omega}
\longrightarrow0.
\]
Using \eqref{tracevolumeconv1}--\eqref{tracevolumeconv2} in the other two terms yields
{\small
\begin{equation}
\label{surfacegreenlimit}
\lim_{\eps\to0}\eps\int_{\Gamma_\eps}v_\eps(x)\cdot\left(\mathbf\Psi\!\left(x,\frac{x}{\eps}\right)n_\eps(x)\right)\dd\sigma_x=\int_\Omega\int_{Y^*}[\nabla_yv_0:\mathbf\Psi+v_0\cdot\diver_y\mathbf\Psi]\dd y\dd x.
\end{equation}
}
For almost every $x$, the Orlicz Green formula on $Y^*$ gives
{\small
\begin{equation}
\label{ygreenformula}
\int_{Y^*}[\nabla_yv_0:\mathbf\Psi+v_0\cdot\diver_y\mathbf\Psi]\dd y=\int_{\partial T}\gamma_0^y v_0\cdot(\mathbf\Psi n)\dd\sigma_y.
\end{equation}
}
The contributions from opposite outer faces of $Y$ cancel because both $v_0$ and
$\mathbf\Psi$ are periodic.  Combining \eqref{surfacegreenlimit} and \eqref{ygreenformula},
we obtain
\begin{equation}
\label{surfaceidentityfluxclass}
\lim_{\eps\to0}
\eps\int_{\Gamma_\eps}
v_\eps\cdot(\mathbf\Psi n)^\eps\dd\sigma_x
=
\int_\Omega\int_{\partial T}
\gamma_0^y v_0\cdot(\mathbf\Psi n)\dd\sigma_y\dd x.
\end{equation}

On the other hand, \eqref{preliminarysurfacelimit} applied to the same test gives
\[
\lim_{\eps\to0}
\eps\int_{\Gamma_\eps}
v_\eps\cdot(\mathbf\Psi n)^\eps\dd\sigma_x
=
\int_\Omega\int_{\partial T}
z_0\cdot(\mathbf\Psi n)\dd\sigma_y\dd x.
\]
Hence
\begin{equation}
\label{z0traceannihilator}
\int_\Omega\int_{\partial T}
\big(z_0-\gamma_0^y v_0\big)
\cdot(\mathbf\Psi n)\dd\sigma_y\dd x=0
\end{equation}
for every admissible $\mathbf\Psi$.  Lemma~\ref{lem:normalfluxdense}, first for separated
variables and then by finite sums, shows that the class of surface tests
$(\mathbf\Psi n)|_{\partial T}$ is dense in
$C_c(\Omega;C(\partial T))^m$.  Since both
$z_0$ and $\gamma_0^y v_0$ belong to $L^{\Phi}(\Omega\times\partial T)^m$,
generalized H\"older's inequality extends \eqref{z0traceannihilator} to the whole test class.
Therefore
\[
z_0=\gamma_0^y v_0
\qquad\text{a.e. on }\Omega\times\partial T.
\]
The surface limit is unique; consequently the original subsequence satisfying
\eqref{tracevolumeconv1}--\eqref{tracevolumeconv2}, without any further extraction, satisfies
\eqref{darcyscaledtraceconv}--\eqref{darcyscaledtraceintegral}.
\end{proof}

\begin{remark}
\label{rem:unfolding}
Theorem~\ref{thm:darcyscaledtrace} is formulated with oscillating test functions, which are those of the Minty
argument of Subsection~\ref{subsec:minty}.  In the language of periodic unfolding in perforated domains
\cite{CioranescuDamlamianDonatoGrisoZaki2012}, with $\mathcal T_\eps(v)(x,y):=v(\eps k+\eps y)$ for
$x\in\eps(k+Y)$, $k\in K_\eps$, $y\in Y^*$, and $\mathcal T^b_\eps(v)$ defined by the same formula for $y\in\partial T$,
it asserts that the boundary unfolding of $v_\eps$ converges to the trace of the limit of the unfolded family.
The exact integration formulas
\[
 \int_\Omega\int_{Y^*}\Phi(|\mathcal T_\eps(v)|)\,\dd y\dd x=\int_{\widehat\Omega_\eps\cap\Omega_\eps}\Phi(|v|)\,\dd x,
 \qquad
 \int_\Omega\int_{\partial T}\Phi(|\mathcal T^b_\eps(v)|)\,\dd\sigma_y\dd x=\eps\int_{\Gamma_\eps}\Phi(|v|)\,\dd\sigma_x,
\]
where $\widehat\Omega_\eps:=\bigcup_{k\in K_\eps}\eps(k+Y)$, together with $\nabla_y\mathcal T_\eps(v)=\eps\,\mathcal T_\eps(\nabla v)$ and
$\mathcal T^b_\eps(v)=\gamma_0^y\mathcal T_\eps(v)$, link the two formulations in the Orlicz setting.
\end{remark}

\begin{corollary}
\label{cor:velocitysurfacecompactness}
Let $u_0$ be the limit obtained in Theorem~\ref{thm:darcyscaledcompactness}.  Then
\begin{equation}
\label{velocitytraceconv}
u_\eps
\overset{2s,\Gamma}{\rightharpoonup}
\gamma_0^y u_0
\quad\text{in }L^{\Phi}(\Omega\times\partial T)^N.
\end{equation}
Moreover,
\begin{equation}
\label{tangentialtraceconv}
(u_\eps)_\tau
\overset{2s,\Gamma}{\rightharpoonup}
P_\tau(y)\gamma_0^y u_0(x,y),
\end{equation}
and
\begin{equation}
\label{velocitytracepuretan}
(u_\eps)_\tau=u_\eps
\overset{2s,\Gamma}{\rightharpoonup}
\gamma_0^y u_0=(\gamma_0^y u_0)_\tau.
\end{equation}
\end{corollary}

\begin{proof}
The hypotheses of Theorem~\ref{thm:darcyscaledtrace} are supplied by the uniform bound
\eqref{uniformenergy} together with the convergences \eqref{u0twoscale} and
\eqref{grad0twoscale}; hence \eqref{velocitytraceconv} follows directly.  Multiplication by the continuous periodic matrix
$P_\tau(y)=I-n(y)\otimes n(y)$ preserves weak surface two-scale convergence, which gives
\eqref{tangentialtraceconv}.  Since $u_\eps\cdot n_\eps=0$ on $\Gamma_\eps$ and, by
Corollary~\ref{cor:limitincompressibility} and \eqref{fastdivnormal}, $\gamma_0^y u_0\cdot n=0$ on $\partial T$,
\eqref{tangentialtraceconv} reduces to \eqref{velocitytracepuretan}.
\end{proof}

Set
\begin{equation}
\label{actualMintyfields}
 S_\eps:=\mathbb S\!\left(\frac{x}{\eps},\eps\D u_\eps\right)\ \text{in }\Omega_\eps,
 \qquad
 B_\eps:=\mathcal B\!\left(\frac{x}{\eps},(u_\eps)_\tau\right)\ \text{on }\Gamma_\eps,
\end{equation}
and denote by $\widetilde S_\eps$ the zero extension of $S_\eps$ to $\Omega$.

\begin{corollary}
\label{cor:frictionprecompactness}
After extraction there exist $S_0\in L^{\widetilde\Phi}(\Omega\times Y^*)^{N\times N}_{\rm sym}$ and a
tangential field $B_0\in L^{\widetilde\Phi}(\Omega\times\partial T)^N$ such that
\begin{align}
\label{stressprelimit}
 \widetilde S_\eps&\overset{2s}{\rightharpoonup}
 \chi_{Y^*}(y)S_0(x,y)
 &&\text{weakly two-scale in }L^{\widetilde\Phi}(\Omega)^{N\times N},\\
\label{frictionprelimit}
 B_\eps&\overset{2s,\Gamma}{\rightharpoonup}B_0(x,y)
 &&\text{weakly in }L^{\widetilde\Phi}(\Omega\times\partial T)^N.
\end{align}
The identifications $S_0=\mathbb S(y,\D_yu_0)$ and $B_0=\mathcal B(y,\gamma_0^yu_0)$ are obtained in
Theorem~\ref{thm:energyclosure}.
\end{corollary}

\begin{proof}
The dual bound \eqref{stressdualbound} and Theorem~\ref{thm:volume-two-scale-compactness}, applied in
$L^{\widetilde\Phi}(\Omega)^{N\times N}$, give a weak two-scale limit of $\widetilde S_\eps$; it vanishes on $\Omega\times T$
because $\widetilde S_\eps=0$ in the holes (test with fields supported in $T$), and symmetry passes to the
limit.  The dual bound \eqref{frictiondualbound} and Theorem~\ref{thm:surface-two-scale-compactness},
applied with $\widetilde\Phi$ in place of $\Phi$, give \eqref{frictionprelimit}; tangentiality passes to
the limit by testing against continuous normal fields.
\end{proof}

\subsection{Vanishing inertia}
\label{subsec:inertia}

The nonlinear identification rests on a pressure-free limit identity, obtained by testing the weak
formulation with exactly solenoidal oscillating fields.  Since the energy controls only an Orlicz norm, whose
integrability on sets of finite measure is governed by the lower growth exponent at infinity $p^-_\infty$, the convective term vanishes against such tests only above the threshold
$\kappa_c$ of Assumption~\ref{ass:kappa}.  We first record the Lebesgue integrability supplied by the
energy.

\begin{lemma}
\label{lem:scaledL2}
Let $\theta(s):=N\left(1/\underline p-1/s\right)_+$, and let $1\le s\le \underline p^{\,*}$ if $\underline p<N$, and
$1\le s<\infty$ if $\underline p\ge N$.  Then there is $C_s>0$, independent of $\eps$, such that
\begin{equation}
\label{scaledLsinterpolation}
 \norm{\bar u_\eps}_{L^s(\Omega)}
 \le C_s\eps^{-\theta(s)}\qquad\forall\,\eps\in\mathscr E.
\end{equation}
\end{lemma}

\begin{proof}
By \eqref{OrliczLpembed} and \eqref{extensionfixedbound},
\begin{equation}
\label{Lpfilled}
 \norm{\bar u_\eps}_{L^{\underline p}(\Omega)}\le C,
 \qquad
 \norm{\nabla \bar u_\eps}_{L^{\underline p}(\Omega)}\le \frac{C}{\eps}.
\end{equation}
For $s\le \underline p$ the finite-measure embedding gives \eqref{scaledLsinterpolation} with $\theta(s)=0$.  For
$s>\underline p$ in the stated range, the Gagliardo--Nirenberg inequality on the fixed $C^{1,1}$ domain $\Omega$ gives
\[
 \norm{\bar u_\eps}_{L^s(\Omega)}
 \le C_s\norm{\nabla \bar u_\eps}_{L^{\underline p}(\Omega)}^{\vartheta}
          \norm{\bar u_\eps}_{L^{\underline p}(\Omega)}^{1-\vartheta}
      +C_s\norm{\bar u_\eps}_{L^{\underline p}(\Omega)},
 \qquad
 \vartheta=N\left(\frac1{\underline p}-\frac1s\right)\in[0,1],
\]
and \eqref{Lpfilled} yields \eqref{scaledLsinterpolation}.
\end{proof}

The next statement shows that inertia vanishes in the limit against solenoidal tests of the type constructed below.

\begin{proposition}
\label{prop:kappathreshold}
Under Assumption~\ref{ass:kappa}, there exists $q_\sharp>\max\{N,p^+\}$ such that the following holds for every
$q\in[q_\sharp,\infty)$.  If $\widehat V_\eps\in \mathbb V_\eps$ has the form
\begin{equation}
\label{testgradientstructure}
 \widehat V_\eps=V^\eps+r_\eps,
 \qquad
 \norm{\nabla V^\eps}_{L^q(\Omega_\eps)}\le C\eps^{-1},
 \qquad
 \norm{\nabla r_\eps}_{L^q(\Omega_\eps)}\le C.
\end{equation}
Then
\begin{equation}
\label{convectivevanishing}
 c_\eps(u_\eps,u_\eps,\widehat V_\eps)\longrightarrow0.
\end{equation}
\end{proposition}

\begin{proof}
Since $2q'\downarrow2<\underline p^{\,*}$ and $2N\big(1/\underline p-1/(2q')\big)_+\downarrow N(2/\underline p-1)_+$ as $q\to\infty$, and
$\kappa>1+N(2/\underline p-1)_+$ by \eqref{underlinep}, we may choose $q_\sharp>\max\{N,p^+\}$ such that, for every $q\ge q_\sharp$, $2q'<\underline p^{\,*}$ and
$\kappa>1+2N\big(1/\underline p -1/(2q')\big)_+$; fix such a $q$.
By \eqref{convectiveIBP}, $c_\eps(u_\eps,u_\eps,\widehat V_\eps)=-\eps^\kappa\int_{\Omega_\eps}(u_\eps\otimes u_\eps):\nabla \widehat V_\eps\,\dd x$.
Let $s_q:=2q'$ and $\theta_q:=\theta(s_q)$.  The choice of $q$ places $s_q$ in the admissible range of
Lemma~\ref{lem:scaledL2}; hence
\begin{equation}
\label{sqbound}
 \norm{u_\eps}_{L^{s_q}(\Omega_\eps)}
 \le \norm{\bar u_\eps}_{L^{s_q}(\Omega)}
 \le C_q\eps^{-\theta_q}.
\end{equation}
H\"older's inequality and \eqref{testgradientstructure} yield
\[
 \eps^\kappa\left|\int (u_\eps\otimes u_\eps):\nabla V^\eps\right|
 \le C\eps^{\kappa-1-2\theta_q}\longrightarrow0
\]
because $\kappa>1+2\theta_q$ by the preceding choice of $q$.  Similarly,
\[
 \eps^\kappa\left|\int (u_\eps\otimes u_\eps):\nabla r_\eps\right|
 \le C\eps^{\kappa-2\theta_q}\longrightarrow0.
\]
This proves \eqref{convectivevanishing}.
\end{proof}

\begin{remark}
\label{rem:kappasharp}
The exponent $\kappa_c$ is the threshold obtained when the Orlicz energy is used through the
embeddings $L^{\Phi}(\Omega_\eps)\hookrightarrow L^{s}(\Omega_\eps)$, $s<p^-_\infty$, that is, through the growth of
$\Phi$ at infinity.  For the power law $\Phi(t)=t^p/p$, and more generally when $\Phi(t)\le Ct^p$ for $t\ge1$
with $p=p^-_\infty$, the following example applies.  For $p<2$ choose a nonzero divergence-free
$\zeta\in C_c^\infty(Y^*)^N$ and one complete cell $k_\eps\in K_\eps$, and set
\[
 u_\eps(x)=\eps^{-N/p}\zeta\left(\frac{x}{\eps}-k_\eps\right).
\]
Then $\|u_\eps\|_{L^p(\Omega_\eps)}+\eps\|\nabla u_\eps\|_{L^p(\Omega_\eps)}\le C$, whereas
\[
 \|u_\eps\|_{L^2(\Omega_\eps)}^2\asymp \eps^{-N(2-p)/p},
\]
so the energy class allows an inertial size of order $\eps^{\kappa-1-N(2-p)/p}$ against a test carrying a
microscopic gradient of order $\eps^{-1}$.  This concentration profile belongs to the energy class without
solving the microscopic problem.  The example therefore shows that $\kappa_c$ is optimal for every argument
based on the Orlicz energy alone; the regime $\kappa\le\kappa_c$, where inertia may enter the limit, is thus the
natural setting for analyses using finer properties of the solutions.
\end{remark}

\subsection{Solenoidal recovery tests and the pressure-free identity}
\label{sec:preliminarylimit}

We introduce the solenoidal test space of the two-scale problem and a regular dense core.

For a cell field $V=V(x,y)$ set
\[
 \overline V(x):=\int_{Y^*}V(x,y)\,\dd y.
\]
Since the average is only known a priori in $L^\Phi(\Omega)^N$, we encode macroscopic incompressibility
and the zero normal flux without introducing a separate normal-trace space.  Define
\begin{equation}
\label{macrosolweakspace}
 \mathcal H_\sigma^\Phi(\Omega):=
 \left\{W\in L^\Phi(\Omega)^N:\
 \int_\Omega W\cdot\nabla\eta\,\dd x=0
 \quad\forall\,\eta\in C^\infty(\overline\Omega)\right\}.
\end{equation}
For sufficiently regular fields this condition is equivalent to
$\diver W=0$ in $\Omega$ and $W\cdot\nu=0$ on $\partial\Omega$.
We then set
\begin{equation}
\label{limitVspace}
\begin{aligned}
\mathbb V_\Phi:=\{V:\;&V\in L^{\Phi}(\Omega\times Y^*)^N,
\quad \nabla_yV\in L^{\Phi}(\Omega\times Y^*)^{N\times N},\\
&V(x,\cdot)\text{ is $Y$-periodic},\quad
\diver_yV=0\text{ in }Y^*,\quad
\gamma_0^y V\cdot n=0\text{ on }\partial T,\\
&\overline V\in\mathcal H_\sigma^\Phi(\Omega)\}.
\end{aligned}
\end{equation}
The trace term belongs to $L^{\Phi}(\Omega\times\partial T)$ by the fixed-cell trace inequality.
Corollary~\ref{cor:limitincompressibility} shows that
\begin{equation}
\label{u0inlimitV}
 u_0\in\mathbb V_\Phi.
\end{equation}

The recovery construction uses arbitrarily large but finite Lebesgue exponents in the cell variable.  We
therefore introduce
\[
 \mathcal W_{\rm reg}(Y^*):=
 \bigcap_{1<q<\infty}W^{1,q}_{\rm per}(Y^*)^N
\]
and the regular core
\begin{equation}
\label{smoothlimitcore}
\mathcal V_{\rm reg}:=\{V\in C_c^2(\Omega;\mathcal W_{\rm reg}(Y^*)):\;
\diver_yV=0,\ \gamma_0^y V\cdot n=0,\ \diver_x\overline V=0\}.
\end{equation}

Both density results below rest on the boundedness of Helmholtz projections in Orlicz spaces.

\begin{lemma}
\label{lem:helmholtz}
\begin{enumerate}[label=\textup{(\roman*)}]
\item The Helmholtz projections $\mathbb P_{\Omega,q}$ of $L^q(\Omega)^N$ onto $L^q_\sigma(\Omega)$, the
closure of $C_{c,\sigma}^\infty(\Omega)^N$ in $L^q(\Omega)^N$, $1<q<\infty$, extend to a bounded projection $\mathbb P_\Omega$ of
$L^\Phi(\Omega)^N$ with range $\mathcal H^\Phi_\sigma(\Omega)$.
\item For $a\in W^{1,q}_{\rm per}(Y^*)^N$, let $\pi$ be the zero-mean periodic solution of
$\Delta_y\pi=\diver_ya$ in $Y^*$, $\partial_n\pi=a\cdot n$ on $\partial T$, and set $\mathbb P_{Y}a:=a-\nabla_y\pi$.
Then $\mathbb P_Y$ extends to a bounded operator on $W^1_{\rm per}L^\Phi(Y^*)^N$, whose range consists of
divergence-free fields tangent to $\partial T$, and $\mathbb P_Yw=w$ for such fields.
\end{enumerate}
\end{lemma}

\begin{proof}
For (i), the projections are furnished by \cite[Theorem~1.4]{SimaderSohr1992} (see also
\cite[Chapter~III]{Galdi2011}); they are consistent on intersections of the spaces $L^q(\Omega)^N$, and their range is
characterized by annihilation of gradients.  For (ii), after identifying opposite faces, $Y^*$ is a compact
$C^{1,1}$ manifold with boundary $\partial T$, and the zero-mean Neumann solution satisfies
\[
 \norm{\pi}_{W^{2,q}(Y^*)}
 \le C_q\Big(
   \norm{\diver_ya}_{L^q(Y^*)}
   +\norm{a\cdot n}_{W^{1-1/q,q}(\partial T)}
   +\norm{\pi}_{L^q(Y^*)}
 \Big)
 \le C_q\norm a_{W^{1,q}(Y^*)},
\]
by the trace theorem and the zero-mean Poincar\'e inequality; hence
$\norm{\mathbb P_Ya}_{W^{1,q}(Y^*)}\le C_q\norm a_{W^{1,q}(Y^*)}$, again consistently in $q$.  These $L^q$ estimates, stated for bounded domains of $\R^N$, transfer to the
perforated torus by localization: the identified outer faces carry no boundary, so that interior estimates
apply near $\partial Y$ and boundary estimates near $\partial T$, and the lower-order terms are removed by the
usual compactness argument, using uniqueness of the zero-mean solution.  In both cases, choose
$1<q_0<p^-\le p^+<q_1<\infty$.  For $Q\in\{\Omega,Y^*\}$, the indices \eqref{indices} show that $\Phi$ satisfies the hypotheses of the
Marcinkiewicz-type theorem \cite[Theorem~VI.3.13]{RaoRen1991}, so that every linear operator bounded on
$L^{q_0}(Q)$ and on $L^{q_1}(Q)$ is bounded on $L^\Phi(Q)$; see also \cite[Theorems~III.5.16 and
IV.8.18]{BennettSharpley1988}, in whose normalization the Boyd indices of $L^\Phi$ lie in $[1/p^+,1/p^-]$.  For (i) this gives
the bounded extension of $\mathbb P_\Omega$; since $L^\Phi(\Omega)\hookrightarrow L^{q_0}(\Omega)$, the range is
$\mathcal H^\Phi_\sigma(\Omega)$.  For (ii), let $\mathcal G(f,F):=\langle f\rangle_{Y^*}+v$, where $v$ is the zero-mean periodic solution of
$\int_{Y^*}\nabla v\cdot\nabla\omega\,\dd y=\int_{Y^*}F\cdot\nabla\omega\,\dd y$ for all
$\omega\in W^{1,q'}_{\rm per}(Y^*)$.  By the $L^q(Y^*)$ theory of the weak Neumann problem, cf.\ \cite[Theorem~1.4]{SimaderSohr1992}, $\mathcal G$ is bounded
from $L^q(Y^*)\times L^q(Y^*)^N$ into $W^{1,q}_{\rm per}(Y^*)$ for every $1<q<\infty$, consistently in $q$, and
$\mathcal G(a,\nabla_ya)=a$.  Interpolating the operator
$(f,F)\mapsto(\mathbb P_Y\mathcal G(f,F),\nabla_y\mathbb P_Y\mathcal G(f,F))$, which is defined on the whole product
space, and restricting it to pairs $(a,\nabla_ya)$ gives the bound on $W^1_{\rm per}L^\Phi(Y^*)^N$.  The last assertion holds because $\pi=0$ when
$\diver_yw=0$ and $w\cdot n=0$.
\end{proof}

We first record the macroscopic density result.

\begin{lemma}
\label{lem:macroOrliczDensity}
Under the present $C^{1,1}$ assumption on $\Omega$,
\begin{equation}
\label{macroOrliczDensity}
 \overline{C_{c,\sigma}^\infty(\Omega)^N}^{\,L^\Phi(\Omega)^N}
 =\mathcal H_\sigma^\Phi(\Omega),
\end{equation}
where $C_{c,\sigma}^\infty(\Omega)^N$ denotes the compactly supported smooth divergence-free fields.
\end{lemma}

\begin{proof}
Let $W\in\mathcal H_\sigma^\Phi(\Omega)$ and let $T_MW$ be a bounded componentwise truncation.  Since
$\Phi\in\Delta_2$, $T_MW\to W$ in $L^\Phi(\Omega)^N$ as $M\to\infty$, and $W^M:=\mathbb P_\Omega(T_MW)\to\mathbb P_\Omega W=W$ by
Lemma~\ref{lem:helmholtz}(i).  For fixed $M$, $W^M\in L^{q_1}_\sigma(\Omega)$ is the limit in $L^{q_1}(\Omega)^N$ of fields
$W^M_j\in C_{c,\sigma}^\infty(\Omega)^N$; since $q_1>p^+$, convergence in $L^{q_1}(\Omega)^N$ implies convergence in
$L^\Phi(\Omega)^N$.  A diagonal choice in $(M,j)$ proves \eqref{macroOrliczDensity}.
\end{proof}

Define the zero-average cell space
\begin{equation}
\label{cellzeroavgspace}
\mathcal Z_0^\Phi(Y^*):=\left\{w\in W^1_{\rm per,0}L^\Phi(Y^*)^N:\
\diver_yw=0,\ \gamma_0 w\cdot n=0,\ \int_{Y^*}w\,\dd y=0\right\}.
\end{equation}

\begin{lemma}
\label{lem:cellaveragelifting}
There exists a linear map
\[
 \mathcal L:\mathbb R^N\longrightarrow
 \bigcap_{1<q<\infty}W^{1,q}_{\rm per}(Y^*)^N
\]
such that, for every $a\in\mathbb R^N$,
\begin{equation}
\label{cellaveragelift}
 \int_{Y^*}\mathcal L a\,\dd y=a,
 \qquad \diver_y\mathcal L a=0,
 \qquad \gamma_0 \mathcal L a\cdot n=0
 \quad\text{on }\partial T.
\end{equation}
Moreover, for every finite $1<q<\infty$ there is $C_q>0$ such that
\[
 \norm{\mathcal La}_{W^{1,q}(Y^*)}\le C_q|a|,
\]
and there is $C>0$ such that
\begin{equation}
\label{LiftingOrliczModular}
 \int_{Y^*}\left[\Phi(|\mathcal La|)+\Phi(|\nabla_y\mathcal La|)\right]\dd y
 \le C\Phi(C|a|).
\end{equation}
Consequently, if $a\in L^\Phi(\Omega)^N$, then
$\mathcal L(a(x))(y)$ and $\nabla_y\mathcal L(a(x))(y)$ belong to
$L^\Phi(\Omega\times Y^*)$, with the corresponding product modular bounded by
$C\int_\Omega\Phi(C|a(x)|)\,\dd x$.
\end{lemma}

\begin{proof}
For $\xi\in\mathbb R^N$, let $\phi_\xi$ be the unique zero-mean periodic solution of
\[
 \Delta_y\phi_\xi=0\quad\text{in }Y^*,\qquad
 \partial_n\phi_\xi=-\xi\cdot n\quad\text{on }\partial T,
\]
and set $h_\xi:=\xi+\nabla_y\phi_\xi$.  The compatibility condition follows from
$\int_{\partial T}\xi\cdot n\,\dd\sigma_y=0$.  Thus $h_\xi$ is periodic, divergence free and tangent
to $\partial T$.  Integration by parts, using $\diver_y h_\xi=0$, $h_\xi\cdot n=0$ on $\partial T$ and cancellation of opposite periodic faces, gives
\[
 \int_{Y^*}h_\xi\cdot\nabla_y\phi_\xi\,\dd y=0.
\]
Since $h_\xi=\xi+\nabla_y\phi_\xi$, it follows that
\[
 \xi\cdot\int_{Y^*}h_\xi\,\dd y
 =\int_{Y^*}|h_\xi|^2\,\dd y>0\qquad(\xi\ne0),
\]
so the linear map $M\xi:=\int_{Y^*}h_\xi\,\dd y$ satisfies
$\xi\cdot M\xi>0$ for $\xi\ne0$.  Hence $M$ is injective and, in finite dimension, invertible.  Define $\mathcal La:=h_{M^{-1}a}$.  Then \eqref{cellaveragelift} follows immediately.  Since $\partial T$ is $C^{1,1}$ and the outer
faces are identified periodically, standard Neumann regularity gives
$\mathcal La\in W^{1,q}_{\rm per}(Y^*)^N$ for every finite $q$, with a bound linear in $|a|$.
Because the range of $\mathcal L$ is finite dimensional, the index inequalities for $\Phi$ convert
these finite-$q$ bounds into \eqref{LiftingOrliczModular}.  Integrating the latter with
$a=a(x)$ gives the asserted product-space bound.
\end{proof}

In the cell, regular fields are dense in the zero-average solenoidal space.

\begin{lemma}
\label{lem:cellRegularDensity}
The subspace
\begin{equation}
\label{cellRegularZeroAvg}
 \mathcal Z_{0,\rm reg}(Y^*):=
 \mathcal Z_0^\Phi(Y^*)\cap
 \bigcap_{1<q<\infty}W^{1,q}_{\rm per}(Y^*)^N
\end{equation}
is dense in $\mathcal Z_0^\Phi(Y^*)$ with respect to the $W^1L^\Phi(Y^*)$ norm.
\end{lemma}

\begin{proof}
Let $w\in\mathcal Z_0^\Phi(Y^*)$.  By the ordinary density of smooth periodic functions in the reflexive
Orlicz--Sobolev space, choose $a_j$ smooth and periodic with
$a_j\to w$ in $W^1L^\Phi(Y^*)^N$.  Put $b_j:=\mathbb P_Ya_j$.  Then $b_j\to w$ in $W^1L^\Phi(Y^*)^N$ by Lemma~\ref{lem:helmholtz}(ii), because $\mathbb P_Yw=w$.  Elliptic
regularity applied to the smooth data gives $b_j\in W^{1,q}_{\rm per}(Y^*)^N$ for every finite $q$.  Let $\mathcal L$ be the cell-average lifting of Lemma~\ref{lem:cellaveragelifting} and define
\[
 w_j:=b_j-\mathcal L\left(\int_{Y^*}b_j\,\dd y\right).
\]
Since $w$ has zero cell average and the averaging map is continuous on $L^\Phi(Y^*)$, the vector in
parentheses tends to zero.  The finite-dimensional boundedness of $\mathcal L$ in every
$W^{1,q}(Y^*)^N$, $q<\infty$, and hence in $W^1L^\Phi(Y^*)^N$, shows that $w_j\to w$ in $W^1L^\Phi(Y^*)^N$.  By construction
$w_j\in\mathcal Z_{0,\rm reg}(Y^*)$, proving the lemma.
\end{proof}

Both density results combine into the density of the regular core in $\mathbb V_\Phi$.

\begin{lemma}
\label{lem:limitdensity}
Under the geometric assumptions on $\Omega$ and $Y^*$, $\mathcal V_{\rm reg}$ is dense in
$\mathbb V_\Phi$ for the norm
\begin{equation}
\label{limitVnorm}
 \norm{V}_{\mathbb V_\Phi}:=
 \norm{V}_{\Phi,\Omega\times Y^*}
 +\norm{\nabla_yV}_{\Phi,\Omega\times Y^*}.
\end{equation}
Consequently the traces converge in $L^{\Phi}(\Omega\times\partial T)$ along such an approximation.
\end{lemma}

\begin{proof}
We decompose into mean and zero-mean parts.  Let $V\in\mathbb V_\Phi$ and write
\begin{equation}
\label{Vmeanzerodecomp}
 V=\mathcal L\overline V+V^0,
 \qquad
 \int_{Y^*}V^0(x,y)\,\dd y=0\quad\text{for a.e. }x.
\end{equation}
Because $\overline V\in\mathcal H_\sigma^\Phi(\Omega)$,
Lemma~\ref{lem:macroOrliczDensity} gives
$\overline V_j\in C_{c,\sigma}^\infty(\Omega)^N$ with
$\overline V_j\to\overline V$ in $L^\Phi(\Omega)^N$.  Estimate
\eqref{LiftingOrliczModular} and Fubini then imply
\begin{equation}
\label{meanpartconv}
 \mathcal L\overline V_j\longrightarrow\mathcal L\overline V
 \quad\text{in the norm }\eqref{limitVnorm}.
\end{equation}

We next approximate the zero-mean part in the product modular.  For almost every $x$, the slice $V^0(x,\cdot)$ belongs to the separable Banach space
$\mathcal Z_0^\Phi(Y^*)$.  Joint measurability of $V^0$ and $\nabla_yV^0$, together with separability,
shows that the slice map is strongly measurable.  Standard simple-function approximation for strongly measurable maps, combined
with the absolute continuity of
\[
 E\longmapsto\int_E\int_{Y^*}
 \left[\Phi(|V^0|)+\Phi(|\nabla_yV^0|)\right]\dd y\dd x,
\]
produces finite-valued maps
\begin{equation}
\label{simplecellapprox}
 S_j(x,y)=\sum_{\ell=1}^{m_j}\mathbf 1_{E_{j\ell}}(x)w_{j\ell}(y),
 \qquad w_{j\ell}\in\mathcal Z_0^\Phi(Y^*),
\end{equation}
such that
\begin{equation}
\label{simplecellmodconv}
 \int_{\Omega\times Y^*}
 \left[\Phi(|S_j-V^0|)+\Phi(|\nabla_yS_j-\nabla_yV^0|)\right]
 \dd y\dd x\longrightarrow0.
\end{equation}
One first approximates the slice map pointwise in the $W^1L^\Phi(Y^*)$ norm by countably-valued simple
maps, then retains finitely many values on a set whose complement has arbitrarily small product modular;
the modular--norm relations on the fixed cell give \eqref{simplecellmodconv}.

By Lemma~\ref{lem:cellRegularDensity}, for each of the finitely many values in \eqref{simplecellapprox}
we may choose
$w_{j\ell}^{\rm reg}\in\mathcal Z_{0,\rm reg}(Y^*)$ so close that replacing $w_{j\ell}$ by
$w_{j\ell}^{\rm reg}$ changes the product modular by less than $1/j$.  Finally choose
$\eta_{j\ell}\in C_c^2(\Omega)$, $0\le\eta_{j\ell}\le1$, approximating
$\mathbf 1_{E_{j\ell}}$ almost everywhere and in measure.  Since the number of terms is finite and the
cell factors are fixed, dominated convergence with
$\Phi(C|w_{j\ell}^{\rm reg}|)+\Phi(C|\nabla_yw_{j\ell}^{\rm reg}|)$ gives
\begin{equation}
\label{smoothcoeffconv}
 V_j^0(x,y):=\sum_{\ell=1}^{m_j}\eta_{j\ell}(x)w_{j\ell}^{\rm reg}(y)
 \longrightarrow V^0
 \quad\text{in }\eqref{limitVnorm}.
\end{equation}
Every $V_j^0$ has zero cell average, is divergence free and tangent in $y$, and is compactly supported
and $C^2$ in $x$ with values in $\mathcal W_{\rm reg}(Y^*)$.

Finally set $V_j:=\mathcal L\overline V_j+V_j^0$.  Then $V_j\in\mathcal V_{\rm reg}$, because its cell average is exactly $\overline V_j$ and hence is
macroscopically divergence free.  Equations \eqref{meanpartconv} and
\eqref{smoothcoeffconv} show that $V_j\to V$ in \eqref{limitVnorm}.  The fixed-cell trace
inequality, integrated over $x$, then yields
\[
 \gamma_0^y V_j\longrightarrow\gamma_0^y V
 \quad\text{in }L^\Phi(\Omega\times\partial T)^N.
\]
This proves the density assertion.
\end{proof}

For a regular two-scale test field $V$, the raw field $V^\eps(x):=V(x,x/\eps)$ is tangent to every pore
wall because $V(x,y)\cdot n(y)=0$, but in general
\begin{equation}
\label{rawdiv}
 \diver V^\eps=(\diver_xV)^\eps.
\end{equation}
We correct this defect without changing the trace on $\Gamma_\eps$, using the following coarse
inverse-divergence estimate.  It is the $q$-independent form of the classical estimate
\cite{Masmoudi2002,Tartar1980} (see also \cite[Theorem~2.3]{Allaire1989}), and it follows from
Lemma~\ref{lem:tartarrestriction}.

\begin{lemma}
\label{lem:coarsebog}
There exists a linear map $\mathfrak B_\eps$, defined on $\bigcup_{1<q<\infty}L_0^q(\Omega_\eps)$ and
independent of $q$, such that, for every $1<q<\infty$,
\[
 \mathfrak B_\eps:L_0^q(\Omega_\eps)
 \longrightarrow W_0^{1,q}(\Omega_\eps)^N,
 \qquad \diver\mathfrak B_\eps g=g,
\]
and
\begin{equation}
\label{coarsebogestimate}
 \norm{\nabla\mathfrak B_\eps g}_{L^q(\Omega_\eps)}
 \le C_q\eps^{-1}\norm g_{L^q(\Omega_\eps)}.
\end{equation}
Moreover, for $w\in W_0^{1,q}(\Omega_\eps)^N$,
\begin{equation}
\label{scaledPoincareq}
 \norm w_{L^q(\Omega_\eps)}
 \le C_q\eps\norm{\nabla w}_{L^q(\Omega_\eps)}.
\end{equation}
The constants are independent of small $\eps$.
\end{lemma}

\begin{proof}
Let $g\in L_0^q(\Omega_\eps)$, let $\widetilde g\in L_0^q(\Omega)$ be its extension by zero through the holes, and let
$v:=\mathfrak B_\Omega\widetilde g\in W_0^{1,q}(\Omega)^N$ be given by the Bogovskii operator of
\cite[Theorem~6.6]{DieningRuzickaSchumacher2010} on $\Omega$, so that $\diver v=\widetilde g$ and
$\norm v_{W^{1,q}(\Omega)}\le C_q\norm g_{L^q(\Omega_\eps)}$.  Set $\mathfrak B_\eps g:=\mathcal R_\eps v$, with
$\mathcal R_\eps$ as in Lemma~\ref{lem:tartarrestriction}.  Since $\diver v=\widetilde g=0$ in every hole,
$m_{\eps,k}(v)=0$, and Lemma~\ref{lem:tartarrestriction}(ii) gives $\diver\mathfrak B_\eps g=g$.  The norm bound of
Lemma~\ref{lem:tartarrestriction} with $\Theta(t)=t^q$ gives
$\norm{\nabla\mathfrak B_\eps g}_{L^q(\Omega_\eps)}\le C_q\eps^{-1}\big(\norm v_{L^q(\Omega)}+\eps\norm{\nabla v}_{L^q(\Omega)}\big)
\le C_q\eps^{-1}\norm g_{L^q(\Omega_\eps)}$, which is \eqref{coarsebogestimate}.  The operators $\mathfrak B_\Omega$ and
$\mathcal R_\eps$ are linear and independent of $q$.

We now justify \eqref{scaledPoincareq}.  Decompose $\Omega_\eps$ into complete fluid cells and
the unperforated boundary remainder.  On each complete cell, scaling the fixed-cell Poincar\'e inequality
for functions with zero trace on $\partial T$ gives
\[
 \norm w_{L^q(\eps(k+Y^*))}
 \le C_q\eps\norm{\nabla w}_{L^q(\eps(k+Y^*))}.
\]
For the boundary remainder, extend $w$ by zero through the holes.  Since
$w\in W_0^{1,q}(\Omega_\eps)^N$, this zero extension belongs to $W_0^{1,q}(\Omega)^N$.
The power-law version of the thin-strip Poincar\'e inequality from
Lemma~\ref{lem:stripPoincare}, together with
$R_\eps\subset \Sigma_{\sqrt N\eps}$, gives
\[
 \norm w_{L^q(R_\eps)}
 \le C_q\eps\norm{\nabla w}_{L^q(\Omega_\eps)}.
\]
Summing the complete-cell and boundary-layer estimates proves
\eqref{scaledPoincareq}.
\end{proof}

Regular limit fields are recovered by exactly solenoidal microscopic fields with the same pore-wall trace.

\begin{theorem}
\label{thm:solenoidaltests}
Let $V\in\mathcal V_{\rm reg}$ and $\max\{N,p^+\}<q<\infty$.  Then, for all sufficiently small $\eps$,
there exists
$\widehat V_\eps\in \mathbb V_\eps$ such that
\begin{align}
\label{teststrongvalue}
 \norm{\widehat V_\eps-V^\eps}_{L^q(\Omega_\eps)}&\le C\eps,\\
\label{teststronggrad}
 \norm{\nabla \widehat V_\eps-\nabla V^\eps}_{L^q(\Omega_\eps)}&\le C,
\end{align}
and therefore, denoting zero extension to $\Omega$ by a tilde,
\begin{align}
\label{testtwoscalevalue}
 \widetilde{\widehat V}_\eps&\overset{2s}{\longrightarrow}
 \chi_{Y^*}(y)V(x,y)
 \quad\text{strongly in }L^{\Phi}(\Omega)^N,\\
\label{testtwoscalegrad}
 \eps\widetilde{\nabla \widehat V_\eps}&\overset{2s}{\longrightarrow}
 \chi_{Y^*}(y)\nabla_yV(x,y)
 \quad\text{strongly in }L^{\Phi}(\Omega)^{N\times N}.
\end{align}
Most importantly, the correction has zero trace on the pore walls, so
\begin{equation}
\label{testtraceexact}
 \widehat V_\eps(x)=V\left(x,\frac{x}{\eps}\right)
 \qquad\text{for }x\in\Gamma_\eps.
\end{equation}
\end{theorem}

\begin{proof}
Because the $x$-support of $V$ is compactly contained in $\Omega$, for small $\eps$ every cell
meeting that support is complete; we denote by
$K_\eps^V:=\{k\in K_\eps:\ \eps(k+\overline Y)\cap\supp_xV\ne\emptyset\}$ the corresponding indices. Put
\[
 g_\eps(x):=\diver V^\eps(x)=(\diver_xV)^\eps(x).
\]
The boundary flux of $V^\eps$ vanishes: the field is zero near $\partial\Omega$ and tangent on
$\Gamma_\eps$. Hence
\begin{equation}
\label{globalmeanzerog}
 \int_{\Omega_\eps}g_\eps\,\dd x=0.
\end{equation}

For $k\in K_\eps^V$ define
\[
 g_{\eps,k}(y):=(\diver_xV)(\eps(k+y),y),
 \qquad
 m_{\eps,k}:=\frac1{|Y^*|}\int_{Y^*}g_{\eps,k}(y)\,\dd y.
\]
Since $\diver_x\overline V=0$,
\[
 \int_{Y^*}(\diver_xV)(\eps k,y)\,\dd y=0.
\]
The $C^1$ dependence on $x$ therefore gives
\begin{equation}
\label{cellmeansmall}
 |m_{\eps,k}|\le C\eps
\end{equation}
uniformly in $k$.

Let $\mathfrak B_{Y^*}$ be the ordinary Bogovskii operator from
$L_0^q(Y^*)$ to $W_0^{1,q}(Y^*)^N$. Set
\[
 \rho_{\eps,k}(y):=
 \mathfrak B_{Y^*}(g_{\eps,k}-m_{\eps,k})(y)
\]
and define
\[
 r_\eps^{(0)}(x):=\eps \rho_{\eps,k}\left(\frac{x}{\eps}-k\right)\quad\text{on }\eps(k+Y^*),\ k\in K_\eps^V,
\]
with value zero elsewhere. Since $\rho_{\eps,k}$ has zero trace on the whole boundary of $Y^*$,
the pieces patch with zero trace across the cell faces and also vanish on every pore wall. We have
\begin{equation}
\label{localcorrdiv}
 \diver r_\eps^{(0)}=g_\eps-m_\eps,
\end{equation}
where $m_\eps=m_{\eps,k}$ on $\eps(k+Y^*)$ for $k\in K_\eps^V$ and $m_\eps=0$ elsewhere. Standard
cell estimates and scaling give
\begin{equation}
\label{localcorrest}
 \norm{r_\eps^{(0)}}_{L^q(\Omega_\eps)}\le C\eps,
 \qquad
 \norm{\nabla r_\eps^{(0)}}_{L^q(\Omega_\eps)}\le C.
\end{equation}
By \eqref{globalmeanzerog} and \eqref{localcorrdiv}, $m_\eps$ has zero global mean. By
\eqref{cellmeansmall},
\begin{equation}
\label{mglobest}
 \norm{m_\eps}_{L^q(\Omega_\eps)}\le C\eps.
\end{equation}
Set $r_\eps^{(1)}:=\mathfrak B_\eps m_\eps$.  Lemma~\ref{lem:coarsebog} and \eqref{mglobest} give
\[
 \norm{\nabla r_\eps^{(1)}}_{L^q(\Omega_\eps)}\le C,
 \qquad
 \norm{r_\eps^{(1)}}_{L^q(\Omega_\eps)}
 \le C\eps\norm{\nabla r_\eps^{(1)}}_{L^q(\Omega_\eps)}\le C\eps.
\]
Finally define
\begin{equation}
\label{Vepsconstruction}
 \widehat V_\eps:=V^\eps-r_\eps^{(0)}-r_\eps^{(1)}.
\end{equation}
Equations \eqref{rawdiv} and \eqref{localcorrdiv} show that
$\diver \widehat V_\eps=0$. Both corrections belong to $W_0^{1,q}(\Omega_\eps)^N$; hence they vanish
on $\partial\Omega$ and on $\Gamma_\eps$. The raw field is zero near $\partial\Omega$ and tangent
on the pore wall. Thus $\widehat V_\eps\in \mathbb V_\eps$ and \eqref{testtraceexact} holds.

The estimates above prove \eqref{teststrongvalue}--\eqref{teststronggrad}.  Since $q>N$, $V$ is continuous
on $\Omega\times\overline{Y^*}$, so that $V^\eps$ is measurable; approximating $V$ in
$C^1(\overline\Omega;W^{1,q}_{\rm per}(Y^*)^N)$ by fields of class $C^1$ in both variables, obtained by a bounded
linear smoothing in $y$, applying the classical chain rule to them and passing to the limit with the
cell-by-cell bound $\norm{W^\eps}_{L^q(\Omega_\eps)}\le C\norm W_{L^q_{\rm per}(Y;C(\overline\Omega))}$ and the
embedding used below, one obtains $V^\eps\in W^{1,q}(\Omega_\eps)^N$ and
\[
 \eps\nabla V^\eps
 =\nabla_yV\left(x,\frac{x}{\eps}\right)
 +\eps\nabla_xV\left(x,\frac{x}{\eps}\right).
\]
Since $q>p^+$, convergence to zero in $L^q(\Omega_\eps)$ implies convergence to zero in $L^{\Phi}(\Omega_\eps)$, by the index bounds for $\Phi$.  Moreover, $V\in C^2_c(\Omega;W^{1,q}_{\rm per}(Y^*)^N)$ gives, by Fubini's theorem,
$\chi_{Y^*}V,\chi_{Y^*}\nabla_yV\in L^q_{\rm per}(Y;W^{2,q}(\Omega))$, extended by zero to $T$, with $x$-support in
$\supp_xV$; since $q>N/2$, the embedding $W^{2,q}(\Omega)\hookrightarrow C(\overline\Omega)$ places both fields in
$L^q_{\rm per}(Y;C(\overline\Omega))$.  For small $\eps$ the zero extensions
satisfy $\widetilde{V^\eps}=(\chi_{Y^*}V)^\eps$ and $\widetilde{\nabla_yV(\cdot,\cdot/\eps)}=(\chi_{Y^*}\nabla_yV)^\eps$, because
$\supp_xV$ meets only complete cells.  Lemma~\ref{lem:Lqadmissible} therefore proves
\eqref{testtwoscalevalue}--\eqref{testtwoscalegrad}.  The two-stage correction is what preserves the leading
scaled gradient: the coarse operator, which loses a factor $\eps^{-1}$ in \eqref{coarsebogestimate}, is applied
only to the cell means $m_\eps$, which are of order $\eps$ because $\diver_x\overline V=0$.
\end{proof}

Testing the microscopic problem with these fields yields the pressure-free limit identity.

\begin{theorem}
\label{thm:preliminarylimit}
Let $S_0$ and $B_0$ be the weak two-scale limits in
\eqref{stressprelimit}--\eqref{frictionprelimit}. Then
\begin{equation}
\label{preliminarylimitidentity}
\int_\Omega\int_{Y^*}S_0:\D_yV\,\dd y\dd x+\int_\Omega\int_{\partial T}B_0\cdot\gamma_0^yV\,\dd\sigma_y\dd x=\int_\Omega\int_{Y^*}f\cdot V\,\dd y\dd x
\qquad\forall\,V\in\mathbb V_\Phi.
\end{equation}
\end{theorem}

\begin{proof}
We first take $V\in\mathcal V_{\rm reg}$, choose $q\ge q_\sharp$ with $q_\sharp$ as in
Proposition~\ref{prop:kappathreshold}, and let $(\widehat V_\eps)$ be the exactly solenoidal recovery sequence of
Theorem~\ref{thm:solenoidaltests}. Use $\widehat V_\eps$ in the microscopic weak formulation
\eqref{weakPeps}.

By \eqref{testtwoscalegrad},
\[
 \eps\widetilde{\D \widehat V_\eps}\overset{2s}{\longrightarrow}\chi_{Y^*}\D_yV
 \quad\text{strongly in }L^{\Phi}(\Omega)^{N\times N}.
\]
Together with the weak two-scale convergence \eqref{stressprelimit},
Lemma~\ref{lem:weak-strong-two-scale-pairing} gives
\begin{equation}
\label{bulkpass}
 \int_{\Omega_\eps}S_\eps:\eps\D \widehat V_\eps\dd x
 \longrightarrow
 \int_\Omega\int_{Y^*}S_0:\D_yV\dd y\dd x.
\end{equation}
Because the correction vanishes on $\Gamma_\eps$, \eqref{testtraceexact} and the weak
surface two-scale convergence of $B_\eps$ yield
\begin{equation}
\label{surfacepass}
 \eps\int_{\Gamma_\eps}B_\eps\cdot(\widehat V_\eps)_\tau\dd\sigma_x
 \longrightarrow
 \int_\Omega\int_{\partial T}B_0\cdot\gamma_0^y V\dd\sigma_y\dd x.
\end{equation}
Here $\gamma_0^y V$ is tangential by definition of $\mathcal V_{\rm reg}$, and it is an admissible surface
test: since $q>N$, the embedding $W^{1,q}(Y^*)\hookrightarrow C^{0,1-N/q}(\overline{Y^*})$ and the continuity of
$x\mapsto V(x,\cdot)$ in $W^{1,q}(Y^*)$ give $\gamma_0^yV\in C_c(\Omega;C(\partial T))^N$.

Proposition~\ref{prop:kappathreshold}, applied to the decomposition furnished by
Theorem~\ref{thm:solenoidaltests}, gives
\begin{equation}
\label{convpass}
 c_\eps(u_\eps,u_\eps,\widehat V_\eps)\longrightarrow0.
\end{equation}
Finally, $f$ is independent of $y$ and belongs to $L^{\widetilde\Phi}(\Omega)^N$, so it is a constant weakly
two-scale convergent family; \eqref{testtwoscalevalue} and Lemma~\ref{lem:weak-strong-two-scale-pairing}
give
\begin{equation}
\label{forcepass}
 \int_{\Omega_\eps}f\cdot \widehat V_\eps\dd x
 \longrightarrow
 \int_\Omega\int_{Y^*}f(x)\cdot V(x,y)\dd y\dd x.
\end{equation}
Combining \eqref{bulkpass}--\eqref{forcepass} proves
\eqref{preliminarylimitidentity} for $V\in\mathcal V_{\rm reg}$.

The three functionals in \eqref{preliminarylimitidentity} are continuous on
$\mathbb V_\Phi$: this follows from generalized H\"older in the bulk, generalized H\"older on the
surface, the fixed-cell trace inequality, and $f\in L^{\widetilde\Phi}(\Omega)^N$. The density statement
of Lemma~\ref{lem:limitdensity} therefore extends the identity to every
$V\in\mathbb V_\Phi$.
\end{proof}

\subsection{Coupled bulk--surface Minty identification}
\label{subsec:minty}

The trace theorem identifies the argument of the friction, but weak two-scale convergence alone does not
place $S_0$ and $B_0$ on their constitutive graphs.  Since bulk stress and wall friction enter the same
energy identity, both are identified by a single Minty argument whose closure condition is the energy
equality supplied by Theorem~\ref{thm:preliminarylimit}.  The $y$-dependence of $\mathbb S$ and
$\mathcal B$ is only measurable, so we first record the admissibility of the comparison fields.

\begin{lemma}
\label{lem:caratheodory}
Let $a:\Omega\times Y\to\mathbb R^m$, $Y$-periodic in $y$, and $h:\Omega\times\partial T\to\mathbb R^m$ be
bounded, measurable in $y$, continuous in $x$ for almost every $y$, and vanishing for $x$ outside a compact
set $K\Subset\Omega$.  If $z_\eps\overset{2s}{\rightharpoonup}z_0$ weakly two-scale in $L^\Phi(\Omega)^m$, then
$\int_\Omega z_\eps\cdot a(x,x/\eps)\,\dd x\to\int_\Omega\int_Y z_0\cdot a\,\dd y\dd x$; if
$z_\eps\overset{2s,\Gamma}{\rightharpoonup}z_0$ in $L^\Phi(\Omega\times\partial T)^m$, then
$\eps\int_{\Gamma_\eps}z_\eps\cdot h(x,x/\eps)\,\dd\sigma_x\to\int_\Omega\int_{\partial T}z_0\cdot h\,\dd\sigma_y\dd x$.
The same holds with $\Phi$ and $\widetilde\Phi$ interchanged.
In particular, if $H\in C_c(\Omega;C_{\rm per}(Y))^{N\times N}_{\rm sym}$ and
$w\in C_c(\Omega;C(\partial T))^N$ with $w\cdot n=0$, then $a=\mathbb S(y,H)$ and $h=\mathcal B(y,w)$
satisfy these hypotheses.
\end{lemma}

\begin{proof}
Since $a$ is bounded, continuous in $x$ for almost every $y$ and supported in $K$, it belongs to
$L^1_{\rm per}(Y;C_0(\Omega))^m$.  Approximating it in this space, extracting a subsequence, multiplying by a fixed
cut-off equal to one on $K$, and composing with the radial truncation $z\mapsto z\min\{1,\|a\|_\infty/|z|\}$, which
preserves continuity, we obtain $a_j\in C_c(\Omega;C_{\rm per}(Y))^m$, supported in a fixed compact subset of
$\Omega$, with $|a_j|\le\|a\|_\infty$ and $a_j\to a$ almost everywhere in $\Omega\times Y$.  Dominated convergence
and $\Delta_2$ give $a_j\to a$ in $L^{\widetilde\Phi}(\Omega\times Y)^m$.  For fixed $j$, the convergence with $a_j$
is the definition of two-scale convergence.  For fixed $A>0$, the function $\widetilde\Phi(|a-a_j|/A)$ is bounded,
continuous in $x$ for almost every $y$ and compactly supported in $x$, so Lemma~\ref{lem:riemann} gives the
convergence of the sampled modulars; as in the proof of Lemma~\ref{lem:Lqadmissible}, the sampled norm
$\|a^\eps-a_j^\eps\|_{\widetilde\Phi,\Omega}$ therefore converges to $\|a-a_j\|_{\widetilde\Phi,\Omega\times Y}$.  By \eqref{OrliczHolder} and the uniform bound of
$z_\eps$, replacing $a^\eps$ by $a_j^\eps$ produces an error which vanishes as $j\to\infty$, uniformly in
$\eps$, and the same holds in the limit term.  The surface case is identical.  Finally, $\mathbb S(y,H)$ and
$\mathcal B(y,w)$ are continuous in $x$ by continuity of the constitutive laws in their second variable,
and bounded by \eqref{Sgrowth} and \eqref{Bgrowth}, since $H$ and $w$ are bounded.
\end{proof}

Set $\mathcal Z_\Phi:=\{z\in L^{\Phi}(\Omega\times\partial T)^N:\ z\cdot n=0\}$.

\begin{lemma}
\label{lem:BNemytskii}
The maps $E\mapsto\mathbb S(\cdot,E)$ from $L^{\Phi}(\Omega\times Y^*)^{N\times N}_{\rm sym}$ into
$L^{\widetilde\Phi}(\Omega\times Y^*)^{N\times N}_{\rm sym}$ and $z\mapsto\mathcal B(\cdot,z)$ from
$\mathcal Z_\Phi$ into $L^{\widetilde\Phi}(\Omega\times\partial T)^N$ are strongly continuous.
\end{lemma}

\begin{proof}
We treat $\mathcal B$; the bulk case is identical.  Let $z_j\to z$ in $L^{\Phi}(\Omega\times\partial T)^N$.  After passage to a
subsequence, $z_j\to z$ almost everywhere and $\{\Phi(C|z_j|)\}_j$ is uniformly integrable for every
fixed $C>0$, by norm convergence and $\Phi\in\Delta_2$.  Continuity of $\mathcal B(y,\cdot)$ gives
$\mathcal B(y,z_j)\to\mathcal B(y,z)$ almost everywhere, and convexity, $\widetilde\Phi\in\Delta_2$ and
\eqref{Bgrowth} give, for a fixed large $C_1$,
\[
 \widetilde\Phi\!\left(
 \frac{|\mathcal B(y,z_j)-\mathcal B(y,z)|}{C_1}
 \right)
 \le C\big[\Phi(C|z_j|)+\Phi(C|z|)\big].
\]
Vitali's theorem yields modular, hence norm, convergence in $L^{\widetilde\Phi}(\Omega\times\partial T)^N$.  Since every
subsequence has a further subsequence with this property, the whole sequence converges.
\end{proof}

The identification rests on the following coupled Minty principle.

\begin{theorem}
\label{thm:coupledMinty}
Let $E_\eps\in L^{\Phi}(\Omega_\eps)_{\rm sym}^{N\times N}$ and $z_\eps\in L^{\Phi}(\Gamma_\eps,\mu_\eps)^N$ with
$z_\eps\cdot n_\eps=0$, and set $S_\eps(x)=\mathbb S(x/\eps,E_\eps(x))$ and
$B_\eps(x)=\mathcal B(x/\eps,z_\eps(x))$.  Assume, after zero extension in the volume,
\begin{align*}
 E_\eps&\overset{2s}{\rightharpoonup}\chi_{Y^*}E_0,
 &S_\eps&\overset{2s}{\rightharpoonup}\chi_{Y^*}S_0,\\
 z_\eps&\overset{2s,\Gamma}{\rightharpoonup}z_0,
 &B_\eps&\overset{2s,\Gamma}{\rightharpoonup}B_0,
\end{align*}
the volume convergences in $L^\Phi(\Omega)^{N\times N}$ and $L^{\widetilde\Phi}(\Omega)^{N\times N}$, the surface
limits $z_0\in L^\Phi(\Omega\times\partial T)^N$ and $B_0\in L^{\widetilde\Phi}(\Omega\times\partial T)^N$ being tangential, and
\begin{equation}
\label{coupledMintylimsup}
 \limsup_{\eps\to0}\left\{
 \int_{\Omega_\eps}S_\eps:E_\eps\dd x
 +\eps\int_{\Gamma_\eps}B_\eps\cdot z_\eps\dd\sigma_x
 \right\}
 \le
 \int_\Omega\int_{Y^*}S_0:E_0\dd y\dd x
 +\int_\Omega\int_{\partial T}B_0\cdot z_0\dd\sigma_y\dd x.
\end{equation}
Then
\begin{equation}
\label{coupledMintyidentification}
 S_0(x,y)=\mathbb S(y,E_0(x,y)),
 \qquad
 B_0(x,y)=\mathcal B(y,z_0(x,y))
 \quad\text{a.e.}
\end{equation}
\end{theorem}

\begin{proof}
Let $H\in C_c(\Omega;C_{\rm per}(Y))^{N\times N}_{\rm sym}$ and $w\in C_c(\Omega;C(\partial T))^N$ with
$w\cdot n=0$, and set $H_\eps(x)=H(x,x/\eps)$ and $w_\eps(x)=w(x,x/\eps)$.  Monotonicity of $\mathbb S$
and $\mathcal B$ gives
\begin{equation}
\label{microCoupledMinty}
0\le\int_{\Omega_\eps}
 [S_\eps-\mathbb S(x/\eps,H_\eps)]:(E_\eps-H_\eps)\dd x
+\eps\int_{\Gamma_\eps}
 [B_\eps-\mathcal B(x/\eps,w_\eps)]\cdot(z_\eps-w_\eps)\dd\sigma_x.
\end{equation}
In the expanded right-hand side, the sum of the two diagonal terms is controlled by
\eqref{coupledMintylimsup}, the cross terms with $H_\eps$ and $w_\eps$ converge by two-scale
convergence, and the terms containing $\mathbb S(x/\eps,H_\eps)$ and $\mathcal B(x/\eps,w_\eps)$ converge
by Lemma~\ref{lem:caratheodory} (the factor $\chi_\eps$ differs from $\chi_{Y^*}(x/\eps)$ only on $R_\eps$,
which does not meet the support of $H$ for small $\eps$).  Hence
\begin{equation}
\label{limitCoupledMinty}
0\le\int_\Omega\int_{Y^*}[S_0-\mathbb S(y,H)]:(E_0-H)\dd y\dd x+\int_\Omega\int_{\partial T}[B_0-\mathcal B(y,w)]\cdot(z_0-w)\dd\sigma_y\dd x.
\end{equation}
Continuous compactly supported fields are dense in $L^{\Phi}(\Omega\times Y^*)^{N\times N}_{\rm sym}$,
and projection of continuous fields by the continuous projector $P_\tau(y)=I-n(y)\otimes n(y)$ shows
that bounded continuous tangential fields with compact $x$-support are dense in $\mathcal Z_\Phi$.  By
Lemma~\ref{lem:BNemytskii} and \eqref{OrliczHolder}, \eqref{limitCoupledMinty} therefore holds
\begin{equation}
\label{limitCoupledMintyall}
 \forall\,H\in L^{\Phi}(\Omega\times Y^*)_{\rm sym}^{N\times N},\qquad
 \forall\,w\in\mathcal Z_\Phi.
\end{equation}
Take first $H=E_0$ and $w=z_0-t\eta$, with $t>0$ and $\eta\in L^\infty(\Omega\times\partial T)^N$
tangential.  Dividing by $t$ and letting $t\downarrow0$ by Lemma~\ref{lem:BNemytskii} gives
$\int_\Omega\int_{\partial T}[B_0-\mathcal B(y,z_0)]\cdot\eta\,\dd\sigma_y\dd x\ge0$; replacing $\eta$ by
$-\eta$ gives equality, and since bounded tangential fields separate points of the tangent bundle,
$B_0=\mathcal B(y,z_0)$.  Take next $w=z_0$ and $H=E_0-tK$, with bounded symmetric $K$; the same argument
gives $S_0=\mathbb S(y,E_0)$.
\end{proof}

For the actual solutions, the total microscopic dissipation equals the work of the force.

\begin{proposition}
\label{prop:totaldissipationlimit}
With $S_\eps$ and $B_\eps$ as in \eqref{actualMintyfields},
\begin{equation}
\label{micrototaldissipation}
D_\eps:=\int_{\Omega_\eps}S_\eps:\eps\D u_\eps\dd x+\eps\int_{\Gamma_\eps}B_\eps\cdot(u_\eps)_\tau\dd\sigma_x=\int_{\Omega_\eps}f\cdot u_\eps\dd x
\longrightarrow
 \int_\Omega\int_{Y^*}f(x)\cdot u_0(x,y)\dd y\dd x.
\end{equation}
\end{proposition}

\begin{proof}
The identity follows by taking $v=u_\eps$ in \eqref{weakPeps}, since
$c_\eps(u_\eps,u_\eps,u_\eps)=0$.  For the limit, write
$\int_{\Omega_\eps}f\cdot u_\eps\dd x=\int_\Omega f\cdot\widetilde u_\eps\dd x$ and use
\eqref{u0twoscale} with the test $f$, which is independent of $y$ and belongs to $L^{\widetilde\Phi}(\Omega)^N$.
\end{proof}

It remains to verify the closure condition with the limit energy.

\begin{theorem}
\label{thm:energyclosure}
The limits satisfy
\begin{equation}
\label{limitenergyequality}
 \int_\Omega\int_{Y^*}f\cdot u_0\,\dd y\dd x
=\int_\Omega\int_{Y^*}S_0:\D_yu_0\,\dd y\dd x
+\int_\Omega\int_{\partial T}B_0\cdot\gamma_0^y u_0\,\dd\sigma_y\dd x,
\end{equation}
and
\begin{equation}
\label{identifiednonlinearities}
 S_0(x,y)=\mathbb S(y,\D_yu_0(x,y)),
 \qquad
 B_0(x,y)=\mathcal B(y,\gamma_0^y u_0(x,y))
 \quad\text{a.e.}
\end{equation}
\end{theorem}

\begin{proof}
By \eqref{u0inlimitV}, $u_0\in\mathbb V_\Phi$, and $V=u_0$ in \eqref{preliminarylimitidentity}
gives \eqref{limitenergyequality}.  Set $E_\eps:=\eps\D u_\eps$ and $z_\eps:=(u_\eps)_\tau$.  By
Theorem~\ref{thm:darcyscaledcompactness}, Corollaries~\ref{cor:velocitysurfacecompactness}
and~\ref{cor:frictionprecompactness}, the hypotheses of Theorem~\ref{thm:coupledMinty} hold with
$E_0=\D_yu_0$ and $z_0=\gamma_0^yu_0$, and by Proposition~\ref{prop:totaldissipationlimit} and
\eqref{limitenergyequality} the left-hand side of \eqref{coupledMintylimsup} converges to its
right-hand side.  Theorem~\ref{thm:coupledMinty} gives \eqref{identifiednonlinearities}.
\end{proof}

Under a uniform form of monotonicity, the weak convergences of the scaled symmetric gradient and of the wall
trace become strong.

\begin{proposition}
\label{prop:corrector}
In addition to the hypotheses of Theorem~\ref{thm:energyclosure}, assume that there is $c_m>0$ such that, for
almost every $y\in Y$ and all $E,G\in\R^{N\times N}_{\rm sym}$, and for $\sigma_y$-almost every $y\in\partial T$ and all
$z,w\in\mathbb T_y$,
\begin{equation}
\label{uniformmonotone}
 [\mathbb S(y,E)-\mathbb S(y,G)]:(E-G)\ge c_m\Phi(|E-G|),\qquad
 [\mathcal B(y,z)-\mathcal B(y,w)]\cdot(z-w)\ge c_m\Phi(|z-w|).
\end{equation}
Then, along the sequence of Theorem~\ref{thm:energyclosure},
\begin{enumerate}[label=\textup{(\alph*)}]
\item $\chi_\eps\,\eps\D u_\eps\overset{2s}{\longrightarrow}\chi_{Y^*}\D_yu_0$ strongly two-scale in
$L^\Phi(\Omega)^{N\times N}$;
\item for every $\eta>0$ there is a tangential $w\in C_c(\Omega;C(\partial T))^N$ with
$\|\gamma_0^yu_0-w\|_{\Phi,\Omega\times\partial T}\le\eta$ and
$\limsup_{\eps\to0}\eps\int_{\Gamma_\eps}\Phi\big(|u_\eps-w(x,x/\eps)|\big)\,\dd\sigma_x\le\eta$; that is, the traces
$u_\eps|_{\Gamma_\eps}$ converge to $\gamma_0^yu_0$ in the sense of strong surface two-scale modular convergence
with respect to the measure $\eps\,\dd\sigma_x$.
\end{enumerate}
Since \eqref{uniformmonotone} implies Assumption~\ref{ass:strictmono}, under the assumptions of Theorem~A both
statements hold for the whole family.  Condition \eqref{uniformmonotone} holds for the benchmarks (i), (ii) and
(iii) of Subsection~\ref{subsec:benchmarks} when $p\ge2$, and for (v).
\end{proposition}

\begin{proof}
Write $E_\eps:=\eps\D u_\eps$, $z_\eps:=u_\eps=(u_\eps)_\tau$ on $\Gamma_\eps$, $E_0:=\D_yu_0$ and $z_0:=\gamma_0^yu_0$.  Let
$H\in C_c(\Omega;C_{\rm per}(Y))^{N\times N}_{\rm sym}$ and tangential $w\in C_c(\Omega;C(\partial T))^N$.  By
\eqref{uniformmonotone}, the right-hand side of \eqref{microCoupledMinty} is bounded below by
$c_m\big[\int_{\Omega_\eps}\Phi(|E_\eps-H_\eps|)\,\dd x+\eps\int_{\Gamma_\eps}\Phi(|z_\eps-w_\eps|)\,\dd\sigma_x\big]$.  The
passage to the limit leading to \eqref{limitCoupledMinty}, which uses \eqref{coupledMintylimsup} established in
the proof of Theorem~\ref{thm:energyclosure}, and the identification \eqref{identifiednonlinearities} give
\[
 c_m\limsup_{\eps\to0}\Big[\int_{\Omega_\eps}\Phi(|E_\eps-H_\eps|)\,\dd x+\eps\int_{\Gamma_\eps}\Phi(|z_\eps-w_\eps|)\,\dd\sigma_x\Big]
 \le R(H,w),
\]
where
\[
\begin{aligned}
 R(H,w):={}&\int_\Omega\int_{Y^*}[\mathbb S(y,E_0)-\mathbb S(y,H)]:(E_0-H)\,\dd y\dd x\\
 &+\int_\Omega\int_{\partial T}[\mathcal B(y,z_0)-\mathcal B(y,w)]\cdot(z_0-w)\,\dd\sigma_y\dd x.
\end{aligned}
\]
By \eqref{OrliczHolder} and the growth conditions, $R(H,w)\le C\big(\|E_0-H\|_{\Phi,\Omega\times Y^*}+\|z_0-w\|_{\Phi,\Omega\times\partial T}\big)$
for $H$ and $w$ in bounded subsets of $L^\Phi(\Omega\times Y^*)$ and $L^\Phi(\Omega\times\partial T)$.  By the density
statements used after \eqref{limitCoupledMinty}, for every $\delta\in(0,1)$ we may choose $H$ and $w$ with
$\|E_0-H\|_{\Phi,\Omega\times Y^*}+\|z_0-w\|_{\Phi,\Omega\times\partial T}\le\delta$, so that the limit superior above is at most
$C\delta$.  This gives (b).

For (a), let $\eta>0$ and $\psi\in L^\Phi(\Omega;C_{\rm per}(Y))^{N\times N}$ with
$\|\chi_{Y^*}E_0-\psi\|_{\Phi,\Omega\times Y}\le\eta/2$, and choose $H$ as above with $\delta$ small.  For small $\eps$,
$\chi_\eps H_\eps=(\chi_{Y^*}H)(x,x/\eps)$, because the $x$-support of $H$ does not meet $R_\eps$.  By
\eqref{modularnorm}, $\limsup_{\eps\to0}\|\chi_\eps(E_\eps-H_\eps)\|_{\Phi,\Omega}\le(C\delta)^{1/p^+}$.  The field
$\chi_{Y^*}H$ satisfies the hypotheses of Lemma~\ref{lem:Lqadmissible}, and the convergence of Luxemburg norms
established in its proof, combined with the approximation of $\psi$ in $L^\Phi(\Omega;C_{\rm per}(Y))$ by elements of
$C_c(\Omega;C_{\rm per}(Y))$ and the bound $\|\psi^\eps-\psi_0^\eps\|_{\Phi,\Omega}\le\|\sup_y|\psi-\psi_0|\|_{\Phi,\Omega}$, gives
$\limsup_{\eps\to0}\|(\chi_{Y^*}H)^\eps-\psi^\eps\|_{\Phi,\Omega}\le\|\chi_{Y^*}H-\psi\|_{\Phi,\Omega\times Y}\le\delta+\eta/2$.
Hence $\limsup_{\eps\to0}\|\chi_\eps E_\eps-\psi^\eps\|_{\Phi,\Omega}\le(C\delta)^{1/p^+}+\delta+\eta/2<\eta$ for $\delta$ small,
which is \eqref{strongTwoScaleDefinition}.

For the benchmarks, the radial maps $t\mapsto t^{p-1}$, $p\ge2$, satisfy
$(|E|^{p-2}E-|G|^{p-2}G):(E-G)\ge c_p|E-G|^p$; summing over the powers in (iii) and using the lower bounds of the
coefficients gives \eqref{uniformmonotone} with $\Phi$ as in each example, and in (v) the bulk map is the sum of
$2\mu(y)E$ and a monotone term.
\end{proof}

\section{Pressure reconstruction, cell problem, and effective nonlinear Darcy law}
\label{sec:effectivedarcy}

Throughout this section, Assumptions~\ref{ass:geometry}--\ref{ass:strictmono} hold, Assumption~\ref{ass:strictmono}
entering only through uniqueness statements (Subsection~\ref{subsec:relaxations}), and $u_0$, $S_0$, $B_0$
are the limits obtained in Section~\ref{sec:homogenization}.  With the nonlinearities identified, we
reconstruct the pressures, analyse the nonlinear cell problem and the effective permeability, prove the
convergence of the microscopic pressure and Theorem~A, and conclude with relaxations of the hypotheses.

\subsection{Microscopic and macroscopic pressure reconstruction}
\label{sec:microderham}

The identity of Theorem~\ref{thm:preliminarylimit} determines the limit equation only on fields which
are divergence-free in the fast variable and whose cell average is macroscopically solenoidal.  The
microscopic de Rham decomposition therefore produces a $y$-pressure $\pi(x,y)$ together with a cell
driving field $\xi(x)\in\R^N$, and the macroscopic de Rham decomposition then identifies $f-\xi$ as a
gradient in $x$.  Throughout this subsection we use the identification \eqref{identifiednonlinearities} of $S_0$ and $B_0$.
Define the cell test space
\begin{equation}
\label{Xcellmicroderham}
\begin{aligned}
\mathbb X_\Phi:=\{V:\;&V\in L^{\Phi}(\Omega\times Y^*)^N,
\quad \nabla_yV\in L^{\Phi}(\Omega\times Y^*)^{N\times N},
\quad V(x,\cdot)\text{ is $Y$-periodic},\\
&\gamma_0^y V(x,\cdot)\cdot n=0\text{ on }\partial T\text{ for a.e. }x\in\Omega\},
\end{aligned}
\end{equation}
its zero-cell-average subspace
\begin{equation}
\label{Xcellzeromicroderham}
 \mathbb X_{\Phi,0}:=\left\{V\in\mathbb X_\Phi:\
 \overline V(x)=0\text{ for a.e. }x\in\Omega\right\},
\end{equation}
and
\begin{equation}
\label{Gcellmicroderham}
 \mathbb G_\Phi:=\left\{g\in L^{\Phi}(\Omega\times Y^*):\
 \int_{Y^*}g(x,y)\,\dd y=0\text{ for a.e. }x\right\}.
\end{equation}
The fixed-cell trace inequality makes the surface trace of every $V\in\mathbb X_\Phi$ an element of
$L^{\Phi}(\Omega\times\partial T)^N$.  The pressure is reconstructed through a zero-average right inverse of
the divergence on the perforated periodic cell.

\begin{lemma}
\label{lem:microdivrightinverse}
There exists a bounded linear operator
$\mathfrak R_\Phi:\mathbb G_\Phi\to \mathbb X_{\Phi,0}$ such that
\begin{equation}
\label{microdivrightinverseproperties}
 \diver_y(\mathfrak R_\Phi g)=g
 \quad\text{a.e. in }\Omega\times Y^*,
\end{equation}
and
\begin{equation}
\label{microdivrightinversebound}
 \norm{\mathfrak R_\Phi g}_{\Phi,\Omega\times Y^*}
 +\norm{\nabla_y\mathfrak R_\Phi g}_{\Phi,\Omega\times Y^*}
 \le C\norm g_{\Phi,\Omega\times Y^*}
 \qquad\forall\,g\in\mathbb G_\Phi.
\end{equation}
\end{lemma}

\begin{proof}
Let first $g\in L^{\Phi}_0(Y^*)$. Since the bounded Lipschitz domain $Y^*$ is a John domain,
the Orlicz Bogovskii theorem \cite[Theorems~3.8 and~6.6]{DieningRuzickaSchumacher2010} (see also
\cite{BreitCianchi2015}) gives $\mathfrak B_{Y^*}g\in W_0^1L^{\Phi}(Y^*)^N$ with
$\diver_y\mathfrak B_{Y^*}g=g$ and the modular estimate
\[
 \int_{Y^*}\left[\Phi(|\mathfrak B_{Y^*}g|)
 +\Phi(|\nabla_y\mathfrak B_{Y^*}g|)\right]\dd y
 \le C\int_{Y^*}\Phi(C|g|)\dd y.
\]
Its zero trace on $\partial Y^*$ is compatible with periodic outer traces and gives zero normal trace on
$\partial T$.  Let $\mathcal L$ be the cell-average lifting of Lemma~\ref{lem:cellaveragelifting}, put
$a_g:=\int_{Y^*}\mathfrak B_{Y^*}g\,\dd y$ and $\mathfrak R_{Y^*}g:=\mathfrak B_{Y^*}g-\mathcal La_g$.  Then
$\mathfrak R_{Y^*}g$ is periodic, tangent to the hole, has zero cell average, and
$\diver_y\mathfrak R_{Y^*}g=g$.  Estimate \eqref{LiftingOrliczModular} and Jensen's inequality for
$a_g$ give
\[
 \int_{Y^*}\left[\Phi(|\mathfrak R_{Y^*}g|)
 +\Phi(|\nabla_y\mathfrak R_{Y^*}g|)\right]\dd y
 \le C\int_{Y^*}\Phi(C|g|)\dd y.
\]
Applying $\mathfrak R_{Y^*}$ pointwise in $x$ and integrating this modular inequality in $x$ gives the
operator $\mathfrak R_\Phi$ and \eqref{microdivrightinverseproperties}; strong measurability follows for
simple functions and then by density and boundedness.
\end{proof}

Define the continuous functional
\begin{equation}
\label{Amicroderham}
 \mathcal A_0(V):=
 \int_\Omega\int_{Y^*}S_0(x,y):\D_yV(x,y)\,\dd y\dd x
 +\int_\Omega\int_{\partial T}B_0(x,y)\cdot\gamma_0^y V(x,y)\,\dd\sigma_y\dd x,
 \qquad V\in\mathbb X_\Phi.
\end{equation}
The generalized H\"older inequality in the bulk and on the surface, together with the fixed-cell trace
estimate, gives
\begin{equation}
\label{Amicroderhambound}
 |\mathcal A_0(V)|\le C
 \left(\norm{S_0}_{\widetilde\Phi,\Omega\times Y^*}
       +\norm{B_0}_{\widetilde\Phi,\Omega\times\partial T}\right)
 \left(\norm V_{\Phi,\Omega\times Y^*}
       +\norm{\nabla_yV}_{\Phi,\Omega\times Y^*}\right).
\end{equation}

\begin{theorem}
\label{thm:microderham}
Assume the hypotheses of Theorem~\ref{thm:energyclosure}. Then there exists a unique pair
$\pi\in L^{\widetilde\Phi}(\Omega\times Y^*)$, $\xi\in L^{\widetilde\Phi}(\Omega)^N$ with
\begin{equation}
\label{micropressurenormalization}
 \int_{Y^*}\pi(x,y)\,\dd y=0
 \qquad\text{for a.e. }x\in\Omega,
\end{equation}
such that
\begin{equation}
\label{microderhamidentity}
\mathcal A_0(V)
 =\int_\Omega\int_{Y^*}\pi(x,y)\diver_yV(x,y)\,\dd y\dd x
 +\int_\Omega\xi(x)\cdot\overline V(x)\,\dd x
 \qquad\forall\,V\in\mathbb X_\Phi.
\end{equation}
Moreover,
\begin{equation}
\label{micropressureestimate}
 \norm\pi_{\widetilde\Phi,\Omega\times Y^*}
 +\norm\xi_{\widetilde\Phi,\Omega}
 \le C\left(
 \norm{S_0}_{\widetilde\Phi,\Omega\times Y^*}
 +\norm{B_0}_{\widetilde\Phi,\Omega\times\partial T}
 \right).
\end{equation}
\end{theorem}

\begin{proof}
We first show that $\mathcal A_0$ annihilates the zero-average solenoidal kernel.  Let
$V\in\mathbb X_{\Phi,0}$ satisfy $\diver_yV=0$.  Its cell average vanishes, hence $V\in\mathbb V_\Phi$,
and \eqref{preliminarylimitidentity} with \eqref{identifiednonlinearities} gives
$\mathcal A_0(V)=\int_\Omega f\cdot\overline V\,\dd x=0$.

To reconstruct the microscopic pressure, set $\Lambda(g):=\mathcal A_0(\mathfrak R_\Phi g)$ for
$g\in\mathbb G_\Phi$.  By Lemma~\ref{lem:microdivrightinverse} and \eqref{Amicroderhambound}, $\Lambda$ is
bounded on the closed subspace $\mathbb G_\Phi\subset L^{\Phi}(\Omega\times Y^*)$, and Hahn--Banach and
reflexive Orlicz duality provide $\varrho\in L^{\widetilde\Phi}(\Omega\times Y^*)$ with
$\Lambda(g)=\int_\Omega\int_{Y^*}\varrho\,g\,\dd y\dd x$ for all $g\in\mathbb G_\Phi$.  Since every
$g\in\mathbb G_\Phi$ has zero $y$-mean, the representative
\begin{equation}
\label{pidefinitionmicro}
 \pi(x,y):=\varrho(x,y)-\frac1{|Y^*|}\int_{Y^*}\varrho(x,z)\,\dd z
\end{equation}
represents the same functional, satisfies \eqref{micropressurenormalization}, and, by Jensen's
inequality, the bound for $\pi$ in \eqref{micropressureestimate}.  If $V\in\mathbb X_{\Phi,0}$ and
$g=\diver_yV$, then $V-\mathfrak R_\Phi g$ has zero cell average and zero fast divergence, so that
\begin{equation}
\label{microderhamzeromean}
 \mathcal A_0(V)=\mathcal A_0(\mathfrak R_\Phi g)=\int_\Omega\int_{Y^*}\pi\,\diver_yV\,\dd y\dd x
 \qquad\forall\,V\in\mathbb X_{\Phi,0}.
\end{equation}

The cell average gives a finite-dimensional obstruction, which defines the driving field.  For
$a\in L^{\Phi}(\Omega)^N$ put $(\mathcal La)(x,y):=\mathcal L(a(x))(y)$.  By
Lemma~\ref{lem:cellaveragelifting}, $\mathcal La\in\mathbb X_\Phi$, $\diver_y\mathcal La=0$ and
$\overline{\mathcal La}=a$, and \eqref{Amicroderhambound} shows that $a\mapsto\mathcal A_0(\mathcal La)$ is
bounded on $L^{\Phi}(\Omega)^N$.  Orlicz duality yields a unique $\xi\in L^{\widetilde\Phi}(\Omega)^N$ with
\begin{equation}
\label{xidefinitionmicro}
 \mathcal A_0(\mathcal La)=\int_\Omega\xi\cdot a\,\dd x
 \qquad\forall\,a\in L^{\Phi}(\Omega)^N,
\end{equation}
and the bound for $\xi$ in \eqref{micropressureestimate}.  For arbitrary $V\in\mathbb X_\Phi$,
$V^0:=V-\mathcal L\overline V$ belongs to $\mathbb X_{\Phi,0}$ and $\diver_yV^0=\diver_yV$; combining
\eqref{microderhamzeromean} and \eqref{xidefinitionmicro} yields \eqref{microderhamidentity}.

For uniqueness, let $(\pi_1,\xi_1)$ and $(\pi_2,\xi_2)$ satisfy \eqref{microderhamidentity} and
\eqref{micropressurenormalization}.  Testing with $V\in\mathbb X_{\Phi,0}$ and using the surjectivity of
$\diver_y:\mathbb X_{\Phi,0}\to\mathbb G_\Phi$ shows that $\pi_1-\pi_2$ depends only on $x$, hence vanishes;
testing with $V=\mathcal La$ gives $\xi_1=\xi_2$.
\end{proof}

Theorem~\ref{thm:microderham} and the pressure-free identity also yield the annihilation property
needed for the macroscopic pressure.  Let $W\in W^1_0L^{\Phi}(\Omega)^N$ with $\diver W=0$.  Then
$W\in\mathcal H_\sigma^\Phi(\Omega)$ and $V=\mathcal LW\in\mathbb V_\Phi$ with $\overline V=W$ and
$\diver_yV=0$, so that \eqref{microderhamidentity} and \eqref{preliminarylimitidentity} give
\begin{equation}
\label{ximacroannihilator2}
 \int_\Omega(f-\xi)\cdot W\,\dd x=0
 \qquad\forall\,W\in\mathcal V^{\Phi}_{0,\sigma}(\Omega)
 :=\{W\in W^1_0L^{\Phi}(\Omega)^N:\ \diver W=0\}.
\end{equation}

\begin{theorem}
\label{thm:macropressure}
Under the hypotheses of Theorem~\ref{thm:energyclosure}, there exists a unique
\begin{equation}
\label{p0space}
 p_0\in W^1L^{\widetilde\Phi}(\Omega)\quad\text{with}\quad\int_\Omega p_0\,\dd x=0,
\end{equation}
such that
\begin{equation}
\label{xifinalmacro}
 \xi(x)=f(x)-\nabla p_0(x)
 \quad\text{for a.e. }x\in\Omega,
\end{equation}
and
\begin{equation}
\label{p0estimate}
 \norm{p_0}_{\widetilde\Phi,\Omega}
 +\norm{\nabla p_0}_{\widetilde\Phi,\Omega}
 \le C\left(
 \norm f_{\widetilde\Phi,\Omega}
 +\norm\xi_{\widetilde\Phi,\Omega}
 \right).
\end{equation}
\end{theorem}

\begin{proof}
Set $F:=f-\xi\in L^{\widetilde\Phi}(\Omega)^N$.  Since the bounded $C^{1,1}$ domain $\Omega$ is a John domain,
\cite[Theorems~3.8 and~6.6]{DieningRuzickaSchumacher2010} provide a bounded right inverse
$\mathfrak B_\Omega:L^{\Phi}_0(\Omega)\to W^1_0L^{\Phi}(\Omega)^N$ of the divergence.  The functional
$g\mapsto-\int_\Omega F\cdot\mathfrak B_\Omega g\,\dd x$ is bounded on $L^{\Phi}_0(\Omega)$, whose dual is
$L^{\widetilde\Phi}(\Omega)$ modulo constants; it is therefore represented by a unique $p_0\in L^{\widetilde\Phi}(\Omega)$
with zero mean, and $\norm{p_0}_{\widetilde\Phi,\Omega}\le C\norm F_{\widetilde\Phi,\Omega}$.  For
$W\in W^1_0L^{\Phi}(\Omega)^N$, the field $W-\mathfrak B_\Omega(\diver W)$ belongs to
$\mathcal V^{\Phi}_{0,\sigma}(\Omega)$, so \eqref{ximacroannihilator2} gives
\[
 \int_\Omega F\cdot W\,\dd x=\int_\Omega F\cdot\mathfrak B_\Omega(\diver W)\,\dd x
 =-\int_\Omega p_0\diver W\,\dd x
 \qquad\forall\,W\in W^1_0L^{\Phi}(\Omega)^N.
\]
Hence $\nabla p_0=F$ in $\mathcal D'(\Omega)^N$, so $p_0\in W^1L^{\widetilde\Phi}(\Omega)$ and
\eqref{xifinalmacro}--\eqref{p0estimate} hold.  Two zero-mean functions with the same gradient coincide on
the connected domain $\Omega$.  This is the Orlicz form of the classical de Rham lemma
\cite{BreitCianchi2015,Galdi2011}.
\end{proof}

\begin{corollary}
\label{cor:twopressurecellbalance}
For almost every $x\in\Omega$, the pair $(u_0(x,\cdot),\pi(x,\cdot))$ satisfies
\begin{equation}
\label{cellbalancewithp0}
\begin{aligned}
&\int_{Y^*}\mathbb S\big(y,\D_yu_0(x,y)\big):\D_yv(y)\,\dd y
+\int_{\partial T}\mathcal B\big(y,\gamma_0^y u_0(x,y)\big)
 \cdot\gamma_0 v(y)\,\dd\sigma_y
-\int_{Y^*}\pi(x,y)\diver_yv(y)\,\dd y\\
&\qquad=\big(f(x)-\nabla p_0(x)\big)\cdot\int_{Y^*}v(y)\,\dd y
\qquad\forall\,v\in\mathcal W^{\rm per}_\Phi(Y^*),
\end{aligned}
\end{equation}
where $\mathcal W^{\rm per}_\Phi(Y^*):=\{v\in W^1_{\rm per}L^{\Phi}(Y^*)^N:\ \gamma_0 v\cdot n=0\text{ on }\partial T\}$.
Equivalently, $(u_0(x,\cdot),\pi(x,\cdot))$ solves the cell problem \eqref{introcellproblem} with
$\zeta=f(x)-\nabla p_0(x)$ in the distributional and weak-traction sense,
with $Y$-periodic velocity and $\int_{Y^*}\pi(x,y)\,\dd y=0$; $\pi(x,\cdot)\in L^{\widetilde\Phi}(Y^*)$ is regarded as a function on
the perforated torus obtained by identifying opposite faces of $Y$.  The macroscopic pressure $p_0$
determines the driving force $f-\nabla p_0$, while the cell pressure $\pi$ enforces the fast
incompressibility constraint.
\end{corollary}

\begin{proof}
Apply \eqref{microderhamidentity} to separated fields $V(x,y)=\eta(x)v(y)$, with $\eta\in C_c^\infty(\Omega)$
and $v$ in a countable dense subset of the admissible cell space, and use the fundamental lemma in $x$.
After removing a common null set, continuity of all terms in the cell norm extends the identity to every
admissible $v$; substituting \eqref{identifiednonlinearities} and \eqref{xifinalmacro} gives
\eqref{cellbalancewithp0}.  Taking $v\in C_c^\infty(Y^*)^N$ and using the symmetry of $\mathbb S$ gives the
first equation in \eqref{introcellproblem}; the second and third are
Corollary~\ref{cor:limitincompressibility}.  For admissible $v$, the Green formula on the perforated torus,
in which the outer faces cancel by periodicity and $v\cdot n=0$ on $\partial T$, shows that
\eqref{cellbalancewithp0} is the weak form of
$\int_{\partial T}\{[\mathbb S(y,\D_yu_0)n]_\tau+\mathcal B(y,\gamma_0^yu_0)\}\cdot v\,\dd\sigma_y=0$, which is the
last equation in \eqref{introcellproblem}.
\end{proof}

\subsection{Nonlinear periodic cell problem and cell pressure}

Having recovered both pressure levels, the limit system can be frozen at each macroscopic point.
Equation~\eqref{cellbalancewithp0} shows that all $x$-dependence enters the microscopic balance through
the driving vector $\zeta=f(x)-\nabla p_0(x)$.  We therefore solve the nonlinear solenoidal cell problem
for an arbitrary parameter $\zeta\in\mathbb R^N$.  Define
\begin{equation}
\label{cellsolspace}
 \mathcal V^\Phi_{\rm per}(Y^*)
 :=\left\{v\in W^1_{\rm per}L^\Phi(Y^*)^N:\
 \diver_yv=0\text{ in }Y^*,\quad
 \gamma_0 v\cdot n=0\text{ on }\partial T\right\}.
\end{equation}
This is a closed reflexive subspace of $W^1_{\rm per}L^\Phi(Y^*)^N$.  We use the bounded cell-average operator
\begin{equation}
\label{cellaverageoperator}
 J:\mathcal V^\Phi_{\rm per}(Y^*)\longrightarrow\R^N,
 \qquad Jv:=\int_{Y^*}v(y)\,\dd y.
\end{equation}
Lemma~\ref{lem:cellaveragelifting} shows that $J$ is onto and that $\mathcal L$ is a bounded
linear right inverse on the finite-dimensional space $\R^N$.

For $w,v\in\mathcal V^\Phi_{\rm per}(Y^*)$, set
{\small
\begin{equation}
\label{celloperatorform}
\langle\mathcal A_Y(w),v\rangle:=\int_{Y^*}\mathbb S(y,\D_yw):\D_yv\,\dd y+\int_{\partial T}\mathcal B(y,\gamma_0 w)\cdot\gamma_0 v\,\dd\sigma_y.
\end{equation}
}
Here the traces are tangential because the fields belong to \eqref{cellsolspace}.

\begin{theorem}
\label{thm:cellproblem}
For every
$\zeta\in\R^N$ there exists a unique
$w_\zeta\in\mathcal V^\Phi_{\rm per}(Y^*)$ satisfying
\begin{equation}
\label{cellweaksol}
 \langle\mathcal A_Y(w_\zeta),v\rangle
 =\zeta\cdot Jv
 \qquad\forall\,v\in\mathcal V^\Phi_{\rm per}(Y^*).
\end{equation}
Existence uses only the monotonicity conditions \eqref{Smonotone} and \eqref{Bmonotone}, and uniqueness uses
Assumption~\ref{ass:strictmono}.  Moreover, for each $\rho>0$ there is $C_\rho>0$ such
that
\begin{equation}
\label{celllocalbound}
 |\zeta|\le \rho
 \quad\Longrightarrow\quad
 \norm{w_\zeta}_{W^1L^\Phi(Y^*)}\le C_\rho.
\end{equation}
The exact energy identity is
\begin{equation}
\label{cellenergyidentity}
 \zeta\cdot Jw_\zeta
 =\int_{Y^*}\mathbb S(y,\D_yw_\zeta):\D_yw_\zeta\,\dd y
 +\int_{\partial T}\mathcal B(y,\gamma_0 w_\zeta)
 \cdot\gamma_0 w_\zeta\,\dd\sigma_y.
\end{equation}
\end{theorem}

\begin{proof}
The growth assumptions and the fixed-cell trace theorem show that $\mathcal A_Y$ maps bounded sets of
$\mathcal V^\Phi_{\rm per}(Y^*)$ into bounded sets of its dual.  Continuity of the two constitutive maps in their
second variables, together with the Orlicz growth bounds and Vitali's theorem, yields hemicontinuity.
Monotonicity is immediate from \eqref{Smonotone} and \eqref{Bmonotone}; existence uses only these
monotonicity conditions.

Introduce, for $v\in\mathcal V^\Phi_{\rm per}(Y^*)$, the cell dissipation
\[
 \mathcal E_Y(v):=
 \int_{Y^*}\Phi(|\D_yv|)\,\dd y
 +\int_{\partial T}\Phi(|\gamma_0 v|)\,\dd\sigma_y.
\]
With $a_0$ as in \eqref{Acoerciveenergy}, the constitutive coercivity gives
\begin{equation}
\label{cellAcoercive}
 \langle\mathcal A_Y(v),v\rangle\ge a_0\mathcal E_Y(v).
\end{equation}
The cell Orlicz--Korn inequality of Proposition~\ref{prop:cellkorn}, followed by the modular--norm
relations, implies
\[
 \frac{\langle\mathcal A_Y(v),v\rangle}
 {\norm v_{W^1L^\Phi(Y^*)}}
 \longrightarrow+\infty
 \qquad\text{as }\norm v_{W^1L^\Phi(Y^*)}\to\infty.
\]
Thus $\mathcal A_Y$ is bounded, monotone, hemicontinuous and coercive on the reflexive space
$\mathcal V^\Phi_{\rm per}(Y^*)$.  Since $v\mapsto\zeta\cdot Jv$ is a bounded linear functional, the
Browder--Minty theorem gives a solution of \eqref{cellweaksol}.

Testing \eqref{cellweaksol} by $v=w_\zeta$ gives \eqref{cellenergyidentity}.  Together with
\eqref{cellAcoercive}, the cell Korn inequality and Young's inequality, it also gives the local bound
\eqref{celllocalbound} for bounded sets of parameters.

For uniqueness, let $w_1,w_2$ solve \eqref{cellweaksol} for the same $\zeta$.
Subtract the two equations and test by $w_1-w_2$.  The resulting sum of two nonnegative terms is zero, and
strict monotonicity of $\mathbb S$ gives $\D_yw_1=\D_yw_2$ almost everywhere.  Lemma~\ref{lem:periodicKorn},
applied to $w_1-w_2$, yields $w_1=w_2$.
\end{proof}

The pressure associated with the solenoidal cell solution follows from the zero-average right inverse
of Lemma~\ref{lem:microdivrightinverse} with $x$ frozen.

\begin{corollary}
\label{cor:cellpressurezeta}
For every $\zeta\in\R^N$ and every cell solution $w_\zeta$, there exists a unique
$\pi_\zeta\in L^{\widetilde\Phi}(Y^*)$ with $\int_{Y^*}\pi_\zeta\,\dd y=0$ such that
\begin{equation}
\label{cellfullweak}
\begin{aligned}
&\int_{Y^*}\mathbb S(y,\D_yw_\zeta):\D_yv\,\dd y
 +\int_{\partial T}\mathcal B(y,\gamma_0 w_\zeta)
 \cdot\gamma_0 v\,\dd\sigma_y
 -\int_{Y^*}\pi_\zeta\diver_yv\,\dd y
 =\zeta\cdot\int_{Y^*}v\,\dd y\\
&\qquad\forall\,v\in\mathcal W^{\rm per}_\Phi(Y^*).
\end{aligned}
\end{equation}
Equivalently, $(w_\zeta,\pi_\zeta)$ solves \eqref{introcellproblem} in the distributional and weak-traction
sense on the perforated torus.
\end{corollary}

\begin{proof}
The functional $v\mapsto\langle\mathcal A_Y(w_\zeta),v\rangle-\zeta\cdot\int_{Y^*}v\,\dd y$ vanishes on
admissible fields with zero divergence by \eqref{cellweaksol}.  The argument of
Theorem~\ref{thm:microderham}, with $x$ frozen and the right inverse $\mathfrak R_{Y^*}$ of
Lemma~\ref{lem:microdivrightinverse}, gives the unique zero-mean $\pi_\zeta$ satisfying \eqref{cellfullweak};
here the mean obstruction is the prescribed vector $\zeta$.  The system \eqref{introcellproblem} follows by
the Green formula, as in Corollary~\ref{cor:twopressurecellbalance}.
\end{proof}

\subsection{Nonlinear permeability map}

The cell-averaged velocity defines the permeability map \eqref{introK}, $\mathcal K(\zeta)=Jw_\zeta$, which closes
the macroscopic constitutive law.

\begin{theorem}
\label{thm:Kproperties}
The permeability map $\mathcal K$ has the following properties.
\begin{enumerate}[label=\textup{(\roman*)}]
\item The map vanishes at the origin: $\mathcal K(0)=0$.
\item The map $\mathcal K$ is strictly monotone:
\begin{equation}
\label{Kstrictmono}
 (\zeta_1-\zeta_2)\cdot
 [\mathcal K(\zeta_1)-\mathcal K(\zeta_2)]>0
 \qquad\forall\,\zeta_1,\zeta_2\in\R^N,\ \zeta_1\ne\zeta_2.
\end{equation}
\item The map $\mathcal K$ is continuous on $\R^N$.
\item There exist constants $c_{\mathcal K},C_{\mathcal K}>0$, depending only on the cell and the structural data, such that
\begin{equation}
\label{KOrliczgrowth}
 \Phi\big(c_{\mathcal K}|\mathcal K(\zeta)|\big)
 +\mathcal E_Y(w_\zeta)
 \le C_{\mathcal K}\widetilde\Phi\big(C_{\mathcal K}|\zeta|\big)
 \qquad\forall\,\zeta\in\R^N.
\end{equation}
\item The cell dissipation is represented by
\begin{equation}
\label{Kenergy}
 \zeta\cdot\mathcal K(\zeta)
 =\langle\mathcal A_Y(w_\zeta),w_\zeta\rangle\ge0,
\end{equation}
and there is $c'_{\mathcal K}>0$, depending only on the cell and the structural data, such that
\begin{equation}
\label{KOrliczcoercive}
 \zeta\cdot\mathcal K(\zeta)\ge c'_{\mathcal K}\,\widetilde\Phi\big(c'_{\mathcal K}|\zeta|\big)
 \qquad\forall\,\zeta\in\R^N.
\end{equation}
In particular,
\begin{equation}
\label{Kenergyproper}
 \frac{\zeta\cdot\mathcal K(\zeta)}{|\zeta|}\longrightarrow+\infty
 \qquad\text{as }|\zeta|\to\infty.
\end{equation}
\end{enumerate}
\end{theorem}

\begin{proof}
For $\zeta=0$, the zero field solves the cell problem, and uniqueness yields $w_0=0$.

Let $w_i=w_{\zeta_i}$.  Subtract the two cell equations and test by $w_1-w_2$.  This gives
\begin{equation}
\label{Kmonotoneidentity}
\begin{aligned}
(\zeta_1-\zeta_2)\cdot[\mathcal K(\zeta_1)-\mathcal K(\zeta_2)]&=\int_{Y^*}[\mathbb S(y,\D_yw_1)-\mathbb S(y,\D_yw_2)]:(\D_yw_1-\D_yw_2)\dd y\\
&\quad+\int_{\partial T}[\mathcal B(y,\gamma_0 w_1)-\mathcal B(y,\gamma_0 w_2)]\cdot(\gamma_0 w_1-\gamma_0 w_2)\dd\sigma_y.
\end{aligned}
\end{equation}
The right-hand side is nonnegative.  If it vanishes, the argument of the uniqueness proof gives
$w_1=w_2$.  Subtracting the cell equations then yields
\[
 (\zeta_1-\zeta_2)\cdot Jv=0
 \qquad\forall\,v\in\mathcal V^\Phi_{\rm per}(Y^*).
\]
Since $J$ is onto, choose $v=\mathcal L(\zeta_1-\zeta_2)$ to conclude $\zeta_1=\zeta_2$.
Thus \eqref{Kstrictmono} is strict whenever the parameters differ.

For continuity, let $\zeta_n\to\zeta$.  The local estimate \eqref{celllocalbound} bounds $w_{\zeta_n}$ in
$\mathcal V^\Phi_{\rm per}(Y^*)$.  After extraction,
$w_{\zeta_n}\rightharpoonup w$ in that space and therefore
$Jw_{\zeta_n}\to Jw$ in $\R^N$.  Since
\[
 \mathcal A_Y(w_{\zeta_n})=J^*\zeta_n\longrightarrow J^*\zeta
 \qquad\text{strongly in }(\mathcal V^\Phi_{\rm per}(Y^*))^*,
\]
we also have
\[
 \langle\mathcal A_Y(w_{\zeta_n}),w_{\zeta_n}-w\rangle
 =\zeta_n\cdot J(w_{\zeta_n}-w)\longrightarrow0.
\]
The bounded monotone hemicontinuous operator $\mathcal A_Y$ is pseudomonotone.  Hence for every
$v\in\mathcal V^\Phi_{\rm per}(Y^*)$,
\[
 \langle\mathcal A_Y(w),w-v\rangle
 \le \liminf_{n\to\infty}
 \langle\mathcal A_Y(w_{\zeta_n}),w_{\zeta_n}-v\rangle
 =\langle J^*\zeta,w-v\rangle.
\]
Replacing $v$ by $w\pm t\eta$, dividing by $t>0$, and letting $t\downarrow0$ gives
$\langle\mathcal A_Y(w)-J^*\zeta,\eta\rangle=0$ for every
$\eta\in\mathcal V^\Phi_{\rm per}(Y^*)$.  Thus
$\mathcal A_Y(w)=J^*\zeta$, so $w$ solves the cell problem with parameter $\zeta$.  Uniqueness gives
$w=w_\zeta$.  Every subsequence has the same limit; therefore
\[
 \mathcal K(\zeta_n)=Jw_{\zeta_n}\longrightarrow Jw_\zeta=\mathcal K(\zeta).
\]

We next prove \eqref{KOrliczgrowth}; fix $\zeta\in\R^N$.  Jensen's inequality on the fixed cell and the cell
Orlicz--Korn inequality yield constants $c_1,C_1>0$ such that
\begin{equation}
\label{KJensenKorn}
 \Phi\big(c_1|\mathcal K(\zeta)|\big)
 \le C_1\mathcal E_Y(w_\zeta).
\end{equation}
On the other hand, coercivity and the energy identity give
\[
 a_0\mathcal E_Y(w_\zeta)
 \le \zeta\cdot\mathcal K(\zeta)
 \le |\zeta|\,|\mathcal K(\zeta)|.
\]
Choose $0<\eta<1$.  Young's inequality for the complementary pair
$(\Phi,\widetilde\Phi)$ gives
\[
 |\zeta|\,|\mathcal K(\zeta)|
 \le \Phi\big(\eta c_1|\mathcal K(\zeta)|\big)
 +\widetilde\Phi\left(\frac{|\zeta|}{\eta c_1}\right).
\]
By \eqref{scalingPhi}, $\Phi(\eta t)\le \eta^{p^-}\Phi(t)$ for $0<\eta\le1$; combine this estimate with
\eqref{KJensenKorn} and choose $\eta$ so small that the term containing
$\Phi(c_1|\mathcal K(\zeta)|)$ can be absorbed.  Using \eqref{scalingPhi} once more to adjust the constants inside the
$N$-functions, one obtains
\[
 \Phi\big(c_{\mathcal K}|\mathcal K(\zeta)|\big)
 \le C_{\mathcal K}\widetilde\Phi(C_{\mathcal K}|\zeta|).
\]
Returning to the coercivity estimate and using the same Young inequality yields the corresponding
bound for $\mathcal E_Y(w_\zeta)$ and completes the proof of \eqref{KOrliczgrowth}.

Identity \eqref{Kenergy} is the reformulation of \eqref{cellenergyidentity} obtained from the definitions
\eqref{celloperatorform} and \eqref{introK}.  To prove \eqref{KOrliczcoercive}, fix $\zeta,a\in\R^N$ and test
\eqref{cellweaksol} by $v=\mathcal La$, which belongs to $\mathcal V^\Phi_{\rm per}(Y^*)$ by
Lemma~\ref{lem:cellaveragelifting}.  Since $J\mathcal La=a$, Young's inequality \eqref{youngPhi}, the growth
conditions \eqref{Sgrowth} and \eqref{Bgrowth}, and the coercivity estimate above give
\[
\begin{aligned}
 \zeta\cdot a
 &=\int_{Y^*}\mathbb S(y,\D_yw_\zeta):\D_y\mathcal La\,\dd y
 +\int_{\partial T}\mathcal B(y,\gamma_0w_\zeta)\cdot\gamma_0\mathcal La\,\dd\sigma_y\\
 &\le C_0\,\mathcal E_Y(w_\zeta)
 +\int_{Y^*}\Phi(c_S|\D_y\mathcal La|)\,\dd y
 +\int_{\partial T}\Phi(c_B|\gamma_0\mathcal La|)\,\dd\sigma_y
 \le\frac{C_0}{a_0}\,\zeta\cdot\mathcal K(\zeta)+C_2\Phi(C_2|a|),
\end{aligned}
\]
where the last step uses \eqref{LiftingOrliczModular}, the bound $|\gamma_0\mathcal La|\le C|a|$ on $\partial T$,
which follows from the embedding $W^{1,q}(Y^*)\hookrightarrow C(\overline{Y^*})$ for $q>N$, and
\eqref{scalingPhi}.  Taking the supremum over $a\in\R^N$ and using
$\sup_{t\ge0}\{|\zeta|t-C_2\Phi(C_2t)\}=C_2\widetilde\Phi\big(|\zeta|/C_2^2\big)$, we obtain
$\zeta\cdot\mathcal K(\zeta)\ge (a_0C_2/C_0)\,\widetilde\Phi(|\zeta|/C_2^2)$, and \eqref{KOrliczcoercive} follows
by \eqref{scalingPhi} applied to $\widetilde\Phi$.  Finally, \eqref{Kenergyproper} follows from
\eqref{KOrliczcoercive}, because $\widetilde\Phi$ is an $N$-function, so that $\widetilde\Phi(t)/t\to\infty$ as
$t\to\infty$.
\end{proof}

\begin{corollary}
\label{cor:KNemytskii}
Let $D\subset\R^N$ have finite measure.  The pointwise realization
\[
 \widehat{\mathcal K}:L^{\widetilde\Phi}(D)^N\longrightarrow L^\Phi(D)^N,
 \qquad
 (\widehat{\mathcal K}\zeta)(x):=\mathcal K(\zeta(x)),
\]
is well defined and continuous.  In particular, for every $\zeta\in L^{\widetilde\Phi}(D)^N$,
\begin{equation}
\label{KNemytskiimodular}
 \int_D\Phi\big(c_{\mathcal K}|\mathcal K(\zeta(x))|\big)\,\dd x
 \le C_{\mathcal K}\int_D\widetilde\Phi\big(C_{\mathcal K}|\zeta(x)|\big)\,\dd x<\infty.
\end{equation}
\end{corollary}

\begin{proof}
By Theorem~\ref{thm:Kproperties}, $\mathcal K$ is continuous and satisfies \eqref{KOrliczgrowth}.  Hence measurability follows from continuity.  Estimate
\eqref{KNemytskiimodular} is the integrated form of \eqref{KOrliczgrowth}; the right-hand side is finite
because $\widetilde\Phi\in\Delta_2$.  If $\zeta_n\to\zeta$ strongly in
$L^{\widetilde\Phi}(D)^N$, then $\mathcal K(\zeta_n)\to\mathcal K(\zeta)$ in measure.  The growth estimate,
the $\Delta_2$ property and the modular convergence of $\zeta_n$ show that
$\{\Phi(c|\mathcal K(\zeta_n)|)\}_n$ is uniformly integrable for a suitable fixed $c>0$.
Using again $\Delta_2$ to control
$\Phi(c|\mathcal K(\zeta_n)-\mathcal K(\zeta)|)$ by the sum of the two individual modulars,
Vitali's theorem yields modular, hence norm, convergence in $L^\Phi(D)^N$.
\end{proof}

The two-sided bounds make the macroscopic problem well posed independently of the homogenization
process.

\begin{corollary}
\label{cor:Darcywellposed}
For every $f\in L^{\widetilde\Phi}(\Omega)^N$, problem \eqref{Darcyweakpressure} has a unique solution
$p_0\in\mathcal P^{\widetilde\Phi}_0(\Omega)$, and there are constants $c,C>0$, independent of $f$, such that
\[
 \int_\Omega\widetilde\Phi(c|\nabla p_0|)\,\dd x\le C\int_\Omega\widetilde\Phi(C|f|)\,\dd x.
\]
\end{corollary}

\begin{proof}
Define $\mathcal T:\mathcal P^{\widetilde\Phi}_0(\Omega)\to\mathcal P^{\widetilde\Phi}_0(\Omega)^*$ by
$\langle\mathcal T(p),\omega\rangle:=-\int_\Omega\mathcal K(f-\nabla p)\cdot\nabla\omega\,\dd x$.  By
Corollary~\ref{cor:KNemytskii}, $\mathcal T$ is bounded and continuous, and \eqref{Kstrictmono} makes it strictly
monotone.  For coercivity, set $\zeta:=f-\nabla p$.  Then
$\langle\mathcal T(p),p\rangle=\int_\Omega\mathcal K(\zeta)\cdot\zeta\,\dd x-\int_\Omega\mathcal K(\zeta)\cdot f\,\dd x$.
By \eqref{KOrliczcoercive}, Young's inequality, \eqref{KOrliczgrowth} and $\Delta_2$, for every $\delta>0$,
\[
 \langle\mathcal T(p),p\rangle
 \ge c'_{\mathcal K}\int_\Omega\widetilde\Phi(c'_{\mathcal K}|\zeta|)\,\dd x
 -\delta C\int_\Omega\widetilde\Phi(c'_{\mathcal K}|\zeta|)\,\dd x
 -C_\delta\int_\Omega\widetilde\Phi(|f|)\,\dd x .
\]
Choosing $\delta$ small and using the Poincar\'e inequality on $\mathcal P^{\widetilde\Phi}_0(\Omega)$ and the
modular--norm relations, we obtain $\langle\mathcal T(p),p\rangle/\norm p\to\infty$ as $\norm p\to\infty$.  The
Browder--Minty theorem on the reflexive space $\mathcal P^{\widetilde\Phi}_0(\Omega)$ gives a solution of
$\mathcal T(p_0)=0$, which is \eqref{Darcyweakpressure}, and strict monotonicity gives uniqueness.  Testing
\eqref{Darcyweakpressure} by $p_0$ and arguing as in the coercivity estimate gives the stated bound.
\end{proof}

\begin{remark}
\label{rem:Kcoercivity}
Estimates \eqref{KOrliczgrowth} and \eqref{KOrliczcoercive} are two-sided Orlicz bounds for the permeability
map: $\mathcal K$ maps $\widetilde\Phi$-growth of the driving force to $\Phi$-growth of the filtration velocity,
and its dissipation $\zeta\cdot\mathcal K(\zeta)$ grows like $\widetilde\Phi(|\zeta|)$.  No balance between the lower
and upper indices is needed; for the power law $\Phi(t)=t^p/p$ both bounds reduce to the scaling
$|\zeta|^{p'}$.
\end{remark}

\begin{remark}
\label{rem:Kpotential}
Suppose, in addition, that
\[
 \mathbb S(y,E)=D_EF(y,E),\qquad
 \mathcal B(y,z)=D_zG(y,z),
\]
where $F$ and $G$ are Carath\'eodory integrands, $C^1$ and strictly convex in their second variables,
and satisfy two-sided
Orlicz coercivity and growth estimates compatible with Assumptions~\ref{ass:stress}--\ref{ass:friction}.
Then the cell equation is the Euler equation of the strictly convex minimization problem associated with
\[
 v\longmapsto
 \int_{Y^*}F(y,\D_yv)\,\dd y
 +\int_{\partial T}G(y,\gamma_0 v)\,\dd\sigma_y
 -\zeta\cdot Jv.
\]
Standard convex duality therefore yields the effective potential
\[
 \Psi(\zeta):=
 \sup_{v\in\mathcal V^\Phi_{\rm per}(Y^*)}
 \left\{\zeta\cdot Jv-
 \int_{Y^*}F(y,\D_yv)\,\dd y-
 \int_{\partial T}G(y,\gamma_0 v)\,\dd\sigma_y\right\},
\]
with $\partial\Psi(\zeta)=\{\mathcal K(\zeta)\}$.  In particular, under these additional hypotheses,
$\Psi\in C^1(\R^N)$ and $\nabla\Psi=\mathcal K$; the homogenization theorem itself uses only the monotone
structure.
\end{remark}

\begin{remark}
\label{rem:Kfriction}
The structure of the cell problem gives the following qualitative properties of $\mathcal K$.
\begin{enumerate}[label=\textup{(\roman*)}]
\item If $\mathbb S(y,-E)=-\mathbb S(y,E)$ and $\mathcal B(y,-z)=-\mathcal B(y,z)$, then $-w_\zeta$ solves the cell problem
with parameter $-\zeta$, and uniqueness gives $\mathcal K(-\zeta)=-\mathcal K(\zeta)$.  More generally, let $R$ be an
orthogonal matrix with $R\Z^N=\Z^N$ which leaves the periodic fluid set $\R^N\setminus(\overline T+\Z^N)$
invariant, and assume $\mathbb S(Ry,RER^T)=R\,\mathbb S(y,E)R^T$ and $\mathcal B(Ry,Rz)=R\,\mathcal B(y,z)$.  The map
$v\mapsto Rv(R^T\cdot)$ preserves $\mathcal V^\Phi_{\rm per}(Y^*)$ and maps $Jv$ to $R\,Jv$; hence $Rw_\zeta(R^T\cdot)$
solves the cell problem with parameter $R\zeta$, and $\mathcal K(R\zeta)=R\,\mathcal K(\zeta)$.  In the Newtonian case, if this
holds for the whole group of signed permutation matrices (a hole with cubic symmetry and isotropic
coefficients), the permeability matrix $K_\beta$ commutes with all of them and is a multiple of the identity.
\item In the variational setting of Remark~\ref{rem:Kpotential}, normalize $F(y,0)=G(y,0)=0$; since
$\mathbb S(y,0)=0$ and $\mathcal B(y,0)=0$, convexity gives $F,G\ge0$.  Restricting the supremum defining $\Psi$ to
fields with zero trace on $\partial T$ shows that $\Psi\ge\Psi_{\rm ns}$, where $\Psi_{\rm ns}$ is the
potential of the no-slip cell problem; and if $G_1\le G_2$, the corresponding potentials satisfy
$\Psi_{G_1}\ge\Psi_{G_2}$.  This ordering of the effective potentials holds in general, and for homogeneous laws it
transfers to $\zeta\cdot\mathcal K(\zeta)$ as follows.  If $F(y,\cdot)$ and
$G(y,\cdot)$ are positively homogeneous of the same degree $p>1$, then
$\Psi$ is positively homogeneous of degree $p'$, as the substitution $v\mapsto t^{1/(p-1)}v$ shows, and Euler's
identity $\zeta\cdot\mathcal K(\zeta)=p'\Psi(\zeta)$ turns these inequalities into
\[
 p'\,\Psi_{\rm ns}(\zeta)\le\zeta\cdot\mathcal K_{G_2}(\zeta)\le\zeta\cdot\mathcal K_{G_1}(\zeta)
 \qquad\forall\,\zeta\in\R^N,\quad G_1\le G_2.
\]
In the Newtonian case $F(y,E)=\mu(y)|E|^2$, $G(y,z)=\beta(y)|z|^2/2$, this means that the permeability matrix
$K_\beta$ is nonincreasing in $\beta$ in the sense of quadratic forms and bounded below by the no-slip
permeability: in this homogeneous setting, increasing wall friction reduces the effective permeability, in
agreement with
\cite{LasseuxValdesParada2024}.
\item In the setting of (ii), assume moreover $F(y,E)\ge c_F\Phi(|E|)$ for some $c_F>0$, and let
\[
 \Psi_{\rm ps}(\zeta):=\sup_{v\in\mathcal V^\Phi_{\rm per}(Y^*)}
 \left\{\zeta\cdot Jv-\int_{Y^*}F(y,\D_yv)\,\dd y\right\}
\]
be the potential of the perfect-slip cell problem, that is, of \eqref{introcellproblem} with $\mathcal B\equiv0$.
It is finite for every hole $T$ allowed by Assumption~\ref{ass:geometry}, including axisymmetric ones, because
periodicity excludes rotations (Lemma~\ref{lem:periodicKorn}).  By \eqref{periodicKorn}, H\"older's inequality and \eqref{modularnorm},
$\zeta\cdot Jv-\int_{Y^*}F(y,\D_yv)\,\dd y\le C|\zeta|\,t-c_F\min\{t^{p^-},t^{p^+}\}$ with $t=\|\D_yv\|_{\Phi,Y^*}$,
so $\Psi_{\rm ps}(\zeta)<\infty$.  The same coercivity, weak upper semicontinuity and strict concavity in
$\D_yv$, combined with \eqref{periodicKorn}, give a unique maximizer, which solves the perfect-slip cell problem.
Since $G\ge0$, $\Psi_G\le\Psi_{\rm ps}$.  In the homogeneous case of (ii), $\Psi_{\rm ps}$ is positively homogeneous
of degree $p'$ and $\zeta\cdot\mathcal K_{\rm ps}(\zeta)=p'\Psi_{\rm ps}(\zeta)$, where $\mathcal K_{\rm ps}$ is the
perfect-slip permeability map; together with (ii),
\[
 p'\,\Psi_{\rm ns}(\zeta)\le\zeta\cdot\mathcal K_G(\zeta)\le\zeta\cdot\mathcal K_{\rm ps}(\zeta)
 \qquad\forall\,\zeta\in\R^N.
\]
In the Newtonian case, $K_{\rm ns}\le K_\beta\le K_{\rm ps}$ in the sense of quadratic forms: the slip-dependent
permeability lies between the no-slip and the perfect-slip permeabilities, as recalled in
Subsection~\ref{subsec:motivation-objectives}.  For the microscopic perfect-slip problem, see
Remark~\ref{rem:perfectslip}.
\end{enumerate}
\end{remark}

\subsection{Convergence of the microscopic pressure}
\label{subsec:pressureconv}

The homogenization was carried out in pressure-free form; we now show that the microscopic pressure itself
converges.  The pressure $p_\eps\in L^{\widetilde\Phi}(\Omega_\eps)$ of Remark~\ref{rem:fixedpressure} satisfies
\begin{equation}
\label{weakPepspressure}
\int_{\Omega_\eps}S_\eps:\eps\D v\,\dd x+\eps\int_{\Gamma_\eps}B_\eps\cdot v_\tau\,\dd\sigma_x
+c_\eps(u_\eps,u_\eps,v)-\int_{\Omega_\eps}p_\eps\diver v\,\dd x=\int_{\Omega_\eps}f\cdot v\,\dd x
\end{equation}
for every $v\in W^1L^\Phi(\Omega_\eps)^N$ with $\gamma_0v=0$ on $\partial\Omega$ and $\gamma_0v\cdot n_\eps=0$ on
$\Gamma_\eps$, with $S_\eps$ and $B_\eps$ as in \eqref{actualMintyfields}.  For such $v$, an integration by parts
in the first term of \eqref{convectiveform}, justified by density and using $\diver u_\eps=0$,
$u_\eps\cdot n_\eps=0$ on $\Gamma_\eps$ and $v=0$ on $\partial\Omega$, gives
\begin{equation}
\label{convectivegeneral}
c_\eps(u_\eps,u_\eps,v)=-\eps^\kappa\int_{\Omega_\eps}(u_\eps\otimes u_\eps):\nabla v\,\dd x.
\end{equation}
Since the Bogovskii constant of $\Omega_\eps$ degenerates as $\eps\to0$, the pressure is extended to the
holes following Tartar \cite{Tartar1980} and Allaire \cite{Allaire1989}, with the restriction operator
$\mathcal R_\eps$ of Lemma~\ref{lem:tartarrestriction} and the notation introduced there.

The extended pressure is defined by
\begin{equation}
\label{pressureextension}
 \widetilde P_\eps:=p_\eps\ \text{ in }\Omega_\eps,
 \qquad
 \widetilde P_\eps:=\frac1{|A_{\eps,k}|}\int_{A_{\eps,k}}p_\eps\,\dd x\ \text{ in }T_{\eps,k},\ k\in K_\eps,
 \qquad
 P_\eps:=\widetilde P_\eps-\frac1{|\Omega|}\int_\Omega\widetilde P_\eps\,\dd x.
\end{equation}
By Lemma~\ref{lem:tartarrestriction}(ii), for every $w\in W^1_0L^\Phi(\Omega)^N$,
\begin{equation}
\label{pressureextensionidentity}
\begin{aligned}
 \int_{\Omega_\eps}p_\eps\diver\mathcal R_\eps w\,\dd x
 &=\int_{\Omega_\eps}p_\eps\diver w\,\dd x+\sum_{k\in K_\eps}\Big(\int_{T_{\eps,k}}\diver w\,\dd x\Big)
 \frac1{|A_{\eps,k}|}\int_{A_{\eps,k}}p_\eps\,\dd x\\
 &=\int_\Omega\widetilde P_\eps\diver w\,\dd x.
\end{aligned}
\end{equation}

\begin{theorem}
\label{thm:pressureconv}
Let $(u_\eps)_{\eps\in\mathscr E}$, with pressures $p_\eps$, be a sequence along which the conclusions of
Section~\ref{sec:homogenization} hold, and let $p_0\in\mathcal P^{\widetilde\Phi}_0(\Omega)$ be the macroscopic pressure
of Theorem~\ref{thm:macropressure}.  Let
$q\ge q_\sharp$, with $q_\sharp$ as in Proposition~\ref{prop:kappathreshold}, and
$\theta_q:=N\big(1/\underline p-1/(2q')\big)_+$, so that $\kappa>1+2\theta_q$.  Then
$P_\eps=P^1_\eps+P^2_\eps$, where
\begin{equation}
\label{pressureconvergence}
 P^1_\eps\longrightarrow p_0\quad\text{strongly in }L^{\widetilde\Phi}(\Omega),
 \qquad
 \|P^2_\eps\|_{L^{q'}(\Omega)}\le C\eps^{\kappa-1-2\theta_q},
\end{equation}
as $\eps\to0$, $\eps\in\mathscr E$.  In particular, $P_\eps\to p_0$ strongly in $L^{q'}(\Omega)$.
\end{theorem}

\begin{proof}
Since $q>p^+$, the index bounds of $\widetilde\Phi$ give $W^{1,q}_0(\Omega)\hookrightarrow W^1_0L^\Phi(\Omega)$ and
$L^{\widetilde\Phi}(\Omega)\hookrightarrow L^{q'}(\Omega)$.

We first bound the pressure functional.  For $w\in W^{1,q}_0(\Omega)^N$, the field $v=\mathcal R_\eps w$ has zero
trace on $\partial\Omega\cup\Gamma_\eps$, so the friction term in \eqref{weakPepspressure} vanishes.  By
\eqref{pressureextensionidentity}, \eqref{weakPepspressure} and \eqref{convectivegeneral},
\[
 -\int_\Omega\widetilde P_\eps\diver w\,\dd x=\ell^1_\eps(w)+\ell^2_\eps(w),
\]
with
\[
 \ell^1_\eps(w):=\int_{\Omega_\eps}f\cdot\mathcal R_\eps w\,\dd x-\int_{\Omega_\eps}S_\eps:\eps\D\mathcal R_\eps w\,\dd x,
 \qquad
 \ell^2_\eps(w):=\eps^\kappa\int_{\Omega_\eps}(u_\eps\otimes u_\eps):\nabla\mathcal R_\eps w\,\dd x.
\]
By \eqref{OrliczHolder}, \eqref{stressdualbound} and Lemma~\ref{lem:tartarrestriction} with $\Theta=\Phi$,
\begin{equation}
\label{ell1bound}
 |\ell^1_\eps(w)|\le C\big(\|w\|_{\Phi,\Omega}+\eps\|\nabla w\|_{\Phi,\Omega}\big)
 \qquad\forall\,w\in W^1_0L^\Phi(\Omega)^N,
\end{equation}
so that $\ell^1_\eps$ extends to $W^1_0L^\Phi(\Omega)^N$.  By H\"older's inequality, \eqref{sqbound} and
Lemma~\ref{lem:tartarrestriction} with $\Theta(t)=t^q$, followed by the Poincar\'e inequality on $\Omega$,
\begin{equation}
\label{ell2bound}
 |\ell^2_\eps(w)|\le\eps^\kappa\|u_\eps\|^2_{L^{2q'}(\Omega_\eps)}\|\nabla\mathcal R_\eps w\|_{L^q(\Omega_\eps)}
 \le C\eps^{\kappa-1-2\theta_q}\|\nabla w\|_{L^q(\Omega)}.
\end{equation}

We next decompose $P_\eps$.  Let $\mathfrak B_\Omega$ be the Bogovskii operator of
\cite[Theorem~6.6]{DieningRuzickaSchumacher2010} on the John domain $\Omega$; it is bounded from $L^\Theta_0(\Omega)$
into $W^1_0L^\Theta(\Omega)^N$ for $\Theta=\Phi$ and $\Theta(t)=t^q$.  For $g\in L^q_0(\Omega)$, since $P_\eps$ and
$\widetilde P_\eps$ differ by a constant and $g$ has zero mean,
\[
 \int_\Omega P_\eps g\,\dd x=\int_\Omega\widetilde P_\eps\diver\mathfrak B_\Omega g\,\dd x
 =-\ell^1_\eps(\mathfrak B_\Omega g)-\ell^2_\eps(\mathfrak B_\Omega g).
\]
By \eqref{ell1bound}, $g\mapsto-\ell^1_\eps(\mathfrak B_\Omega g)$ is bounded on $L^\Phi_0(\Omega)$; since the dual
of $L^\Phi_0(\Omega)$ is $L^{\widetilde\Phi}(\Omega)$ modulo constants, it is represented by a unique
$P^1_\eps\in L^{\widetilde\Phi}(\Omega)$ with zero mean and $\|P^1_\eps\|_{\widetilde\Phi,\Omega}\le C$.  By
\eqref{ell2bound}, $g\mapsto-\ell^2_\eps(\mathfrak B_\Omega g)$ is represented by a unique
$P^2_\eps\in L^{q'}(\Omega)$ with zero mean and $\|P^2_\eps\|_{L^{q'}(\Omega)}\le C\eps^{\kappa-1-2\theta_q}$.  Thus
$P_\eps-P^1_\eps-P^2_\eps$ annihilates $L^q_0(\Omega)$, is therefore constant, and vanishes because all three
functions have zero mean.

To obtain the compactness of $P^1_\eps$, equip $L^\Phi(\Omega)^N\times L^\Phi(\Omega)^{N\times N}$ with the norm
$\|a\|_{\Phi,\Omega}+\eps\|b\|_{\Phi,\Omega}$.  The map $w\mapsto(w,\nabla w)$ is an isometry from
$W^1_0L^\Phi(\Omega)^N$, normed by $\|w\|_{\Phi,\Omega}+\eps\|\nabla w\|_{\Phi,\Omega}$, onto a subspace, on which
\eqref{ell1bound} bounds $\ell^1_\eps$.  The Hahn--Banach theorem and Orlicz duality give
$a_\eps\in L^{\widetilde\Phi}(\Omega)^N$ and $A_\eps\in L^{\widetilde\Phi}(\Omega)^{N\times N}$ with
$\|a_\eps\|_{\widetilde\Phi,\Omega}\le C$, $\|A_\eps\|_{\widetilde\Phi,\Omega}\le C\eps$ and
$\ell^1_\eps(w)=\int_\Omega a_\eps\cdot w\,\dd x+\int_\Omega A_\eps:\nabla w\,\dd x$.  Hence $P^1_\eps=P^{1a}_\eps+P^{1b}_\eps$,
where the zero-mean functions $P^{1a}_\eps$ and $P^{1b}_\eps$ represent $g\mapsto-\int_\Omega a_\eps\cdot\mathfrak B_\Omega g$
and $g\mapsto-\int_\Omega A_\eps:\nabla\mathfrak B_\Omega g$ on $L^\Phi_0(\Omega)$.  The second satisfies
$\|P^{1b}_\eps\|_{\widetilde\Phi,\Omega}\le C\eps$.  The first is the image of the bounded family $(a_\eps)$ under the
adjoint of $\mathfrak B_\Omega:L^\Phi_0(\Omega)\to L^\Phi(\Omega)^N$, which is compact as the composition of a bounded
operator into $W^1_0L^\Phi(\Omega)^N$ with the compact embedding \eqref{compactOrliczEmbedding}; by Schauder's
theorem, $(P^{1a}_\eps)$ is relatively compact in $L^{\widetilde\Phi}(\Omega)$.  Consequently, every sequence
$\eps_n\to0$ has a subsequence along which $P^1_{\eps_n}\to P$ strongly in $L^{\widetilde\Phi}(\Omega)$ for some
zero-mean $P$, and then $P_{\eps_n}\to P$ strongly in $L^{q'}(\Omega)$.

It remains to identify $P$.  Since $\overline T\subset Y$, there is $\delta\in(0,1/2)$ with
$\overline T\subset(\delta,1-\delta)^N$.  Fix $i\in\{1,\dots,N\}$, choose $j\ne i$ and a $1$-periodic
$a\in C^\infty(\R)$ supported in $[-\delta/2,\delta/2]+\Z$ with $\int_0^1a=1$, and set $\omega^{(i)}(y):=a(y_j)e_i$.  Then
$\omega^{(i)}\in C^\infty_{\rm per}(Y)^N$, $\diver_y\omega^{(i)}=0$, $\omega^{(i)}$ vanishes in a neighbourhood of
$\overline T+\Z^N$, and $\int_Y\omega^{(i)}\,\dd y=e_i$.  For $\psi\in C_c^\infty(\Omega)$ the field
$v_\eps(x):=\psi(x)\omega^{(i)}(x/\eps)$ belongs to $C_c^\infty(\Omega)^N$ and vanishes near the holes, so it is
admissible in \eqref{weakPepspressure}, its friction term vanishes, and $\diver v_\eps=\omega^{(i)}(x/\eps)\cdot\nabla\psi$.
Since $v_\eps$ vanishes in the holes and has compact support,
$\int_{\Omega_\eps}p_\eps\diver v_\eps\,\dd x=\int_\Omega P_{\eps}\,\omega^{(i)}(x/\eps)\cdot\nabla\psi\,\dd x$, which
converges along the subsequence to $\int_\Omega P\,\partial_i\psi\,\dd x$, because $P_{\eps_n}\to P$ strongly in
$L^{q'}(\Omega)$ and $\omega^{(i)}(\cdot/\eps)\rightharpoonup e_i$ weakly-$*$ in $L^\infty(\Omega)^N$.  Next,
$\eps\D v_\eps=\psi\,(\D_y\omega^{(i)})(x/\eps)+O(\eps)$ uniformly, so \eqref{stressprelimit} and
\eqref{identifiednonlinearities} give
$\int_{\Omega_\eps}S_\eps:\eps\D v_\eps\,\dd x\to\int_\Omega\int_{Y^*}\mathbb S(y,\D_yu_0):\psi\D_y\omega^{(i)}\,\dd y\dd x$.
By \eqref{convectivegeneral}, \eqref{sqbound} and $\|\nabla v_\eps\|_{L^q(\Omega)}\le C\eps^{-1}$,
$|c_\eps(u_\eps,u_\eps,v_\eps)|\le C\eps^{\kappa-1-2\theta_q}\to0$.  Finally
$\int_{\Omega_\eps}f\cdot v_\eps\,\dd x\to\int_\Omega\psi f_i\,\dd x$.  Passing to the limit in \eqref{weakPepspressure},
\[
 \int_\Omega\psi\int_{Y^*}\mathbb S(y,\D_yu_0):\D_y\omega^{(i)}\,\dd y\dd x-\int_\Omega P\,\partial_i\psi\,\dd x
 =\int_\Omega\psi f_i\,\dd x.
\]
Taking $v=\omega^{(i)}$ in \eqref{cellbalancewithp0}, whose friction and cell-pressure terms vanish, gives
$\int_{Y^*}\mathbb S(y,\D_yu_0):\D_y\omega^{(i)}\,\dd y=f_i-\partial_ip_0$ for almost every $x$.  Hence
$-\int_\Omega\psi\,\partial_ip_0\,\dd x=\int_\Omega P\,\partial_i\psi\,\dd x$ for every $\psi$ and every $i$, that is,
$\nabla(P-p_0)=0$ in $\mathcal D'(\Omega)$.  Since $\Omega$ is connected and both functions have zero mean,
$P=p_0$.  The limit being independent of the subsequence, \eqref{pressureconvergence} holds along
$\mathscr E$.
\end{proof}

\subsection{Proof of Theorem A}
\label{subsec:proofA}

Let $(u_\eps)_{0<\eps<\eps_0}$ be any family of weak solutions of \eqref{Peps}, and let $\mathscr E$ be an
arbitrary sequence tending to zero.  By Theorem~\ref{thm:darcyscaledcompactness} and
Corollaries~\ref{cor:velocitysurfacecompactness} and~\ref{cor:frictionprecompactness}, after extraction, the
convergences \eqref{finalconvU}--\eqref{finalconvtrace} hold, together with \eqref{stressprelimit} and
\eqref{frictionprelimit}.  Theorem~\ref{thm:energyclosure} identifies $S_0=\mathbb S(y,\D_yu_0)$ and
$B_0=\mathcal B(y,\gamma_0^yu_0)$, which gives \eqref{finalconvstress}--\eqref{finalconvfriction}.

Theorems~\ref{thm:microderham} and~\ref{thm:macropressure} provide the cell pressure $\pi$ and the
normalized macroscopic pressure $p_0\in\mathcal P^{\widetilde\Phi}_0(\Omega)$ with driving force $f-\nabla p_0$.  By
Corollary~\ref{cor:twopressurecellbalance}, for almost every $x\in\Omega$ the pair $(u_0(x,\cdot),\pi(x,\cdot))$
solves the cell problem \eqref{cellfullweak} with $\zeta=f(x)-\nabla p_0(x)$.  The uniqueness statements of
Theorem~\ref{thm:cellproblem} and Corollary~\ref{cor:cellpressurezeta} give
\begin{equation}
\label{u0cellrepresentation}
 u_0(x,\cdot)=w_{\,f(x)-\nabla p_0(x)},
 \qquad
 \pi(x,\cdot)=\pi_{\,f(x)-\nabla p_0(x)}
 \qquad\text{for a.e. }x\in\Omega,
\end{equation}
which is \eqref{finalu0cell}.  Integrating the first identity over $Y^*$ gives
$U=\mathcal K(f-\nabla p_0)$ almost everywhere, and $U\in L^\Phi(\Omega)^N$ by Corollary~\ref{cor:KNemytskii}.
Corollary~\ref{cor:limitincompressibility} gives $\diver U=0$ in $\Omega$ and $U\cdot\nu=0$ on $\partial\Omega$ in the
weak sense $\int_\Omega U\cdot\nabla\omega\,\dd x=0$ for all $\omega\in W^1L^{\widetilde\Phi}(\Omega)$; this is
\eqref{finalDarcyboxed}, equivalently \eqref{Darcyweakpressure}, the pairing being well defined by
Corollary~\ref{cor:KNemytskii} and \eqref{OrliczHolder}.  The properties of $\mathcal K$ are those of
Theorem~\ref{thm:Kproperties}.

Problem \eqref{Darcyweakpressure} has at most one solution in $\mathcal P^{\widetilde\Phi}_0(\Omega)$.  Indeed, let
$p_1,p_2$ be two solutions and put $\zeta_i=f-\nabla p_i$.  Subtracting the two identities and choosing
$\omega=p_1-p_2$ gives, since $\nabla(p_1-p_2)=-(\zeta_1-\zeta_2)$,
\[
 \int_\Omega
 [\mathcal K(\zeta_1)-\mathcal K(\zeta_2)]\cdot
 (\zeta_1-\zeta_2)\,\dd x=0.
\]
By \eqref{Kstrictmono} the integrand is nonnegative and vanishes only where $\zeta_1=\zeta_2$; hence
$\nabla p_1=\nabla p_2$, and connectedness of $\Omega$ with the zero-mean normalization gives $p_1=p_2$.

Consequently $p_0$ is uniquely determined, and so is $u_0$ by \eqref{u0cellrepresentation}; the limits of
\eqref{finalconvU}--\eqref{finalconvfriction} do not depend on the extracted subsequence.  Since every
sequence $\eps\to0$ has a subsequence along which these convergences hold with the same limits, they hold
for the whole family.  Finally, Theorem~\ref{thm:pressureconv}, applied along every such sequence, and the
uniqueness of $p_0$ give the convergence of the extended pressure for the whole family.  This proves
Theorem~A.

\subsection{Relaxations of the hypotheses}
\label{subsec:relaxations}

The arguments above remain valid, with the modifications described below, under weaker hypotheses.

\begin{remark}
\label{rem:nonstrict}
Assumption~\ref{ass:strictmono} is used only for uniqueness: in the proofs of
Theorems~\ref{thm:cellproblem} and~\ref{thm:Kproperties}, vanishing dissipation gives $\D_yw_1=\D_yw_2$, and
Lemma~\ref{lem:periodicKorn} then gives $w_1=w_2$ for every hole geometry, without any strict monotonicity of
$\mathcal B$.  Without Assumption~\ref{ass:strictmono}, Sections~\ref{sec:microscopic-estimates}
and~\ref{sec:homogenization} and Subsection~\ref{sec:microderham} remain valid, and so does
Theorem~\ref{thm:cellproblem} except for the uniqueness of $w_\zeta$, since its existence part, the local bound
\eqref{celllocalbound} and the energy identity do not use strict monotonicity.  For $\zeta\in\R^N$ let
$\mathcal S_\zeta\subset\mathcal V^\Phi_{\rm per}(Y^*)$ be the set of solutions of \eqref{cellweaksol} and
$\boldsymbol{\mathcal K}(\zeta):=J\mathcal S_\zeta\subset\R^N$.
\begin{enumerate}[label=\textup{(\alph*)}]
\item By Minty's lemma for the monotone hemicontinuous operator $\mathcal A_Y$,
$\mathcal S_\zeta=\bigcap_{v\in\mathcal V^\Phi_{\rm per}(Y^*)}\{w:\ \langle J^*\zeta-\mathcal A_Y(v),w-v\rangle\ge0\}$, an
intersection of closed half-spaces; $\mathcal S_\zeta$ is therefore closed and convex, and bounded by
\eqref{celllocalbound}, hence weakly compact.  Thus $\boldsymbol{\mathcal K}(\zeta)$ is a nonempty convex compact set.
\item Identity \eqref{Kmonotoneidentity} holds for any $w_i\in\mathcal S_{\zeta_i}$, so
$(\zeta_1-\zeta_2)\cdot(k_1-k_2)\ge0$ for $k_i\in\boldsymbol{\mathcal K}(\zeta_i)$.
\item The graph of $\boldsymbol{\mathcal K}$ is closed: if $\zeta_n\to\zeta$ and $k_n=Jw_n\to k$ with
$w_n\in\mathcal S_{\zeta_n}$, the pseudomonotonicity argument in the proof of Theorem~\ref{thm:Kproperties} shows
that every weak limit $w$ of $(w_n)$ belongs to $\mathcal S_\zeta$, and $Jw=k$.
\item $\boldsymbol{\mathcal K}$ is maximal monotone.  Let $(\zeta_0,k_0)$ satisfy $(\zeta-\zeta_0)\cdot(k-k_0)\ge0$ for all
$\zeta$ and $k\in\boldsymbol{\mathcal K}(\zeta)$.  For $\eta\in\R^N$ and $t>0$ choose $k_t\in\boldsymbol{\mathcal K}(\zeta_0+t\eta)$;
then $\eta\cdot(k_t-k_0)\ge0$, and, by \eqref{celllocalbound} and (c), a subsequence of $k_t$ converges as
$t\downarrow0$ to some $k_\eta\in\boldsymbol{\mathcal K}(\zeta_0)$ with $\eta\cdot(k_\eta-k_0)\ge0$.  If
$k_0\notin\boldsymbol{\mathcal K}(\zeta_0)$, the separation theorem for the compact convex set $\boldsymbol{\mathcal K}(\zeta_0)$
gives $\eta$ with $\eta\cdot(k-k_0)<0$ for all $k\in\boldsymbol{\mathcal K}(\zeta_0)$, a contradiction.
\item The proofs of \eqref{KOrliczgrowth} and \eqref{KOrliczcoercive} apply to every element of $\mathcal S_\zeta$, so
$\Phi(c_{\mathcal K}|k|)\le C_{\mathcal K}\widetilde\Phi(C_{\mathcal K}|\zeta|)$ and
$\zeta\cdot k\ge c'_{\mathcal K}\widetilde\Phi(c'_{\mathcal K}|\zeta|)$ for all $k\in\boldsymbol{\mathcal K}(\zeta)$.
\end{enumerate}
Consequently, for every sequence $\eps\to0$ there are a subsequence and a pair
$(u_0,p_0)\in\mathbb W^\Phi_{\rm per}(\Omega\times Y^*)^N\times\mathcal P^{\widetilde\Phi}_0(\Omega)$ such that
\eqref{finalconvU}--\eqref{finalconvfriction} hold along the subsequence, $u_0(x,\cdot)\in\mathcal S_{f(x)-\nabla p_0(x)}$
for almost every $x$, and the filtration velocity satisfies the Darcy inclusion
\[
 U\in\boldsymbol{\mathcal K}(f-\nabla p_0),\qquad\diver U=0\ \text{ in }\Omega,\qquad U\cdot\nu=0\ \text{ on }\partial\Omega;
\]
the arguments are those of Subsection~\ref{subsec:proofA}, in which uniqueness enters only through the
identification of the limit.  The limit is then characterized along subsequences, and
Theorem~\ref{thm:pressureconv} holds along the same subsequences.
\end{remark}

\begin{remark}
\label{rem:perfectslip}
At the microscopic level the rescaled cells are not periodic, and the rigid motions in
Proposition~\ref{prop:cellkorn} are removed by friction and impermeability jointly.  Assume that
\[
 \{r\in\mathcal R:\ r\cdot n=0\text{ on }\partial T\}=\{0\},
\]
a condition which excludes the axisymmetric hole geometries supporting tangential rigid rotations.  Then
impermeability alone removes the rigid motions in Proposition~\ref{prop:cellkorn}, the coercivity
\eqref{Bcoercive} may be replaced by $\mathcal B(y,z)\cdot z\ge0$, and perfect slip $\mathcal B\equiv0$ becomes
admissible in Theorem~A.  At the cell level no such condition is needed, by Lemma~\ref{lem:periodicKorn}; the
perfect-slip permeability then bounds the slip-dependent one from above (Remark~\ref{rem:Kfriction}).
\end{remark}

\begin{remark}
\label{rem:energyinequality}
Condition \eqref{NSexponent} is used only for microscopic solvability and for the energy equality at fixed
$\eps$, and the passage to the limit uses this equality only through an upper bound.  Hence, if
$p^-_\infty>2N/(N+2)$, so that $\underline p$ may be chosen with $\underline p^{\,*}>2$ in place of
$\underline p>3N/(N+2)$ in \eqref{underlinep}, the conclusions of Theorem~A remain valid for every family of
solutions of \eqref{Peps} satisfying the energy inequality
\[
 \int_{\Omega_\eps}\mathbb S\!\left(\frac{x}{\eps},\eps\D u_\eps\right):\eps\D u_\eps\,\dd x
 +\eps\int_{\Gamma_\eps}\mathcal B\!\left(\frac{x}{\eps},(u_\eps)_\tau\right)\cdot(u_\eps)_\tau\,\dd\sigma_x
 \le\int_{\Omega_\eps}f\cdot u_\eps\,\dd x:
\]
Proposition~\ref{prop:kappathreshold} requires only $2<\underline p^{\,*}$, and the closure in
Theorem~\ref{thm:energyclosure} uses only an upper bound for the limit superior of the microscopic dissipation.
This isolates the energy inequality as the only property of the microscopic solutions used in the limit; for
steady no-slip flows of generalized Newtonian fluids, weak solutions in the range
$2N/(N+2)<p^-_\infty\le3N/(N+2)$ are constructed by Lipschitz truncation \cite{FrehseMalekSteinhauer2003}, and the
above conclusions apply to those solutions which satisfy the energy inequality.
\end{remark}

\section*{Statements and declarations}

\textit{Funding.} This research did not receive any specific grant from funding agencies in the public,
commercial, or not-for-profit sectors.

\textit{Declaration of competing interest.} The authors declare that they have no known competing financial
interests or personal relationships that could have appeared to influence the work reported in this paper.

\textit{Data availability.} No data were used for the research described in the article.

\textit{Author contributions.} All authors contributed to the conception, analysis and writing of this work,
and read and approved the final manuscript.

\textit{Declaration of generative AI and AI-assisted technologies.} During the preparation of this work the authors used Claude (Anthropic) in order to assist
with the literature review.  The author(s) reviewed and edited the output as needed and take full responsibility for the content of the published article.


\end{document}